\documentclass[10pt,a4paper,notitlepage]{article}
\usepackage[latin1]{inputenc}
\usepackage{graphicx,amsmath,amsfonts,amssymb,amsthm}
\usepackage[mathscr]{euscript}
\usepackage{rotating}

\date{January 23, 2006\footnote{Corrected May 2, 2006.}} 
\author{Douglas Lundholm\\ \scriptsize{F01, KTH}}
\title{Geometric (Clifford) algebra and its applications}

\numberwithin{equation}{section}
\numberwithin{table}{section}

\theoremstyle{plain}
\newtheorem{thm}{Theorem}[section]
\newtheorem{lem}[thm]{Lemma}
\newtheorem{prop}[thm]{Proposition}
\newtheorem*{cor}{Corollary}

\theoremstyle{definition}
\newtheorem{defn}{Definition}[section]

\newtheorem{exmp}{Example}[section]

\theoremstyle{remark}
\newtheorem*{rem}{Remark}

\newcommand{\cl}{\mathcal{C}l}
\newcommand{\liprod}{\ \raisebox{.2ex}{$\llcorner$}\ }
\newcommand{\riprod}{\ \raisebox{.2ex}{$\lrcorner$}\ }
\newcommand{\iprod}{\ \raisebox{.3ex}{$\scriptscriptstyle \bullet$}\ }
\newcommand{\cliffconj}{{\scriptscriptstyle \square}}
\newcommand{\symdiff}{\!\bigtriangleup\!}
\newcommand{\dual}{^\mathbf{c}}

\begin{document}

\pagestyle{empty}

\maketitle
\thispagestyle{empty}

\begin{abstract}
	In this Master of Science Thesis I introduce geometric algebra both
	from the traditional geometric setting of vector spaces, and also from
	a more combinatorial view which simplifies common relations and operations.
	This view enables us to define Clifford algebras with scalars in arbitrary rings
	and provides new suggestions for an infinite-dimensional approach.
	
	Furthermore, I give a quick review of classic results regarding
	geometric algebras, such as their classification in terms of
	matrix algebras, the connection to orthogonal and Spin groups,
	and their representation theory. A number of lower-dimensional examples
	are worked out in a systematic way using so called norm functions, 
	while general applications of representation theory
	include normed division algebras and vector fields on spheres.
	
	I also consider examples in relativistic physics, where reformulations in terms of 
	geometric algebra give rise to both computational and conceptual
	simplifications.
	
\end{abstract}


\newpage
\tableofcontents
\newpage

\pagestyle{plain}
\setcounter{page}{1}

\section{Introduction}

	The foundations of geometric algebra, or what today is more commonly known
	as Clifford algebra, were put forward already in 1844 by Grassmann. He 
	introduced vectors, scalar products and extensive quantities such as exterior
	products. His ideas were far ahead of his time and formulated in an abstract
	and rather philosophical form which was hard to follow for contemporary
	mathematicians. Because of this, his work was largely ignored until around
	1876, when Clifford took up Grassmann's ideas and formulated a natural
	algebra on vectors with combined interior and exterior products.
	He referred to this as an application of Grassmann's geometric algebra.
	
	Due to unfortunate historic events, such as Clifford's early death in 1879,
	his ideas did not reach the wider part of the mathematics community.
	Hamilton had independently invented the quaternion algebra which
	was a special case of Grassmann's constructions, a fact Hamilton quickly realized
	himself. Gibbs reformulated, largely due to a misinterpretation, the
	quaternion algebra to a system for calculating with vectors in three
	dimensions with scalar and cross products.
	This system, which today is taught at an elementary academic level,
	found immediate applications in physics, which at that time circled
	around Newton's mechanics and Maxwell's electrodynamics.
	Clifford's algebra only continued to be employed within small
	mathematical circles, while physicists struggled to transfer the
	three-dimensional concepts in Gibbs' formulation to special relativity
	and quantum mechanics.
	Contributions and independent reinventions of Grassmann's and Clifford's constructions
	were made along the way by Cartan, Lipschitz, Chevalley, Riesz, Atiyah, Bott, Shapiro,
	and others.
	
	Only in 1966 did Hestenes identify the Dirac algebra, which had
	been constructed for relativistic quantum mechanics, as the geometric
	algebra of spacetime. This spawned new interest in geometric algebra,
	and led, though with a certain reluctance in the scientific community,
	to applications and reformulations in a wide range of fields in
	mathematics and physics.
	More recent applications include image analysis, computer vision,
	robotic control and electromagnetic field simulations.
	Geometric algebra is even finding its way into the computer game
	industry.
	
	There are a number of aims of this Master of Science Thesis. Firstly, I want to give
	a compact introduction to geometric algebra which sums up classic results
	regarding its basic structure and operations, the relations between different geometric
	algebras, and the important connection to orthogonal groups via Spin groups.
	I also clarify a number of statements which have been used in a rather sloppy,
	and sometimes incorrect, manner in the literature. All stated theorems are accompanied 
	by proofs, or references to where a strict proof can be found.
	
	Secondly, I want to show why I think that geometric and Clifford algebras
	are important, by giving examples of applications in mathematics and physics.
	The applications chosen cover a wide range of topics, some with no
	direct connection to geometry. The applications in physics serve to illustrate 
	the computational and, most importantly, conceptual simplifications that
	the language of geometric algebra can provide.
	
	Another aim of the thesis is to present some of the ideas of my supervisor
	Lars Svensson in the subject of generalizing Clifford algebra in the algebraic
	direction. I also present some of my own ideas regarding norm functions
	on geometric algebras.
	
	The reader will be assumed to be familiar with basic algebraic concepts such
	as tensors, fields, rings and homomorphisms. Some basics in topology are also helpful.
	To really appreciate the examples in physics, the reader should be
	familiar with special relativity and preferably also relativistic electrodynamics
	and quantum mechanics.
	For some motivation and a picture of where we are heading, it could be
	helpful to have seen some examples of geometric algebras before.
	For a quick 10-page introduction with some applications in physics,
	see \cite{lundholm}.
	
	Throughout, we will use the name \emph{geometric algebra} in the context
	of vector spaces, partly in honor of Grassmann's contributions, but
	mainly for the direct and natural connection to geometry that this algebra admits.
	In a more general algebraic setting, where a combinatorial rather than 
	geometric interpretation exists, we call the corresponding construction
	\emph{Clifford algebra}.
	
	\newpage
	
\section{Foundations}

	In this section we define geometric algebra and work out a number of
	its basic properties. We consider the definition that is most common in the
	mathematical literature, namely as a quotient space 
	on the tensor algebra of a vector space with a quadratic form.
	We see that this leads, in the finite-dimensional case, to the
	equivalent definition as an algebra with generators $\{e_i\}$ satisfying
	$e_i e_j + e_j e_i = 2 g_{ij}$ for some metric $g$.
	This is perhaps the most well-known definition.
	
	We go on to consider an alternative definition of geometric algebra
	based on its algebraic and combinatorial features. 
	The resulting algebra, here called \emph{Clifford algebra} due to
	its higher generality but less direct connection to geometry, allows us to define common operations and
	prove fundamental identities in a remarkably simple way compared
	to traditional fomulations. 
	
	Returning to the vector space setting, we go on to study some of the
	geometric features from which geometric algebra earns its name.
	We also consider parts of the extensive linear function theory 
	which exists for geometric algebras.
	
	Finally, we note that the generalized Clifford algebra offers
	interesting views regarding the infinite-dimensional case.

\subsection{Geometric algebra $\mathcal{G}(\mathcal{V},q)$}

	The traditional definition of geometric algebra is carried out in the context
	of vector spaces with an inner product, or more generally a quadratic form.
	We consider here a vector space $\mathcal{V}$ of arbitrary dimension over 
	some field $\mathbb{F}$.
	
	\begin{defn} \label{def_quad_form}
		A \emph{quadratic form} $q$ on a vector space $\mathcal{V}$ is a map 
		$q\!: \mathcal{V} \to \mathbb{F}$ such that
		\begin{displaymath}
		\begin{array}{rl}
			i)  & q(\alpha v) = \alpha^2 q(v) \quad \forall \ \alpha \in \mathbb{F}, v \in \mathcal{V}\\[5pt]
			ii) & q(v+w) - q(v) - q(w) \quad \textrm{is linear in both $v$ and $w$.}
		\end{array}
		\end{displaymath}
		The bilinear form $\beta_q (v,w) := \frac{1}{2} \big( q(v+w) - q(v) - q(w) \big)$ is called the \emph{polarization} of $q$.
	\end{defn}

	\begin{exmp}
		If $\mathcal{V}$ has a bilinear form $\langle \cdot , \cdot \rangle$
		then $q(v) := \langle v, v \rangle$ is a quadratic form
		and $\beta_q$ is the symmetrization of $\langle \cdot , \cdot \rangle$.
		This could be positive definite (an inner product), 
		or indefinite (a \emph{metric} of arbitrary signature).
	\end{exmp}
	
	\begin{exmp}
		If $\mathcal{V}$ is a normed vector space over $\mathbb{R}$, 
		with norm denoted by $| \cdot |$, where the	parallelogram 
		identity $|x+y|^2 + |x-y|^2 = 2|x|^2 + 2|y|^2$ holds
		then $q(v) := |v|^2$ is a quadratic form and $\beta_q$ is an inner 
		product on $\mathcal{V}$. This is a classic result, sometimes called the Jordan-von Neumann theorem.
	\end{exmp}
	
	Let $\mathcal{T(V)} := \bigoplus_{k=0}^{\infty} \bigotimes^{k} \mathcal{V}$ 
	denote the tensor algebra on $\mathcal{V}$, the elements of which are
	finite sums of tensors of different grades on $\mathcal{V}$.
	Consider the ideal generated by all elements of the form\footnote{Mathematicians 
	often choose a different sign convention here, resulting in reversed signature
	in many of the following results. The convention used here seems more natural in my opinion,
	since e.g. squares of vectors in euclidean spaces become positive instead of negative.}
	$v \otimes v - q(v)$ for vectors $v$,
	\begin{equation}
		\mathcal{I}_q(\mathcal{V}) := \Big\{ A \otimes \big(v \otimes v - q(v)\big) \otimes B \quad : \quad v \in V,\ A,B \in \mathcal{T(V)} \Big\}.
	\end{equation}
	We define the geometric algebra over $\mathcal{V}$ by quoting out this ideal from $\mathcal{T(V)}$.
	
	\begin{defn} \label{def_g}
		The \emph{geometric algebra} $\mathcal{G}(\mathcal{V},q)$ over the 
		vector space $\mathcal{V}$ with quadratic form $q$ is defined by
		\begin{displaymath}
			\mathcal{G}(\mathcal{V},q) := \mathcal{T(V)} / \mathcal{I}_q(\mathcal{V}).
		\end{displaymath}
		When it is clear from the context what vector space or quadratic form we are working with,
		we will often denote $\mathcal{G}(\mathcal{V},q)$ by $\mathcal{G}(\mathcal{V})$, or just $\mathcal{G}$.
	\end{defn}
	
	\noindent
	The product in $\mathcal{G}$, called the \emph{geometric} or \emph{Clifford product}, 
	is inherited from the tensor product in $\mathcal{T(V)}$ and
	we denote it by juxtaposition (or $\cdot$ if absolutely necessary),
	\begin{displaymath}
	\setlength\arraycolsep{2pt}
	\begin{array}{ccl}
		\mathcal{G} \times \mathcal{G} 	& \to 		& \mathcal{G}, \\
		(A,B) 							& \mapsto 	& AB := [ A \otimes B ].
	\end{array}
	\end{displaymath}
	Note that this product is bilinear and associative. We immediately find the following 
	identities on $\mathcal{G}$ for $v,w \in \mathcal{V}$:
	\begin{equation} \label{generator_relations}
		v^2 = q(v) \quad \Rightarrow \quad vw + wv = 2\beta_q(v,w).
	\end{equation}
	
	One of the most important consequences of this definition of the 
	geometric algebra is the following 
	
	\begin{prop}[Universality] \label{prop_universality}
		Let $\mathcal{A}$ be an associative algebra over $\mathbb{F}$ with a unit denoted by $1_\mathcal{A}$. 
		If $f\!: \mathcal{V} \to \mathcal{A}$ is linear and
		\begin{equation} \label{universality_property}
			f(v)^2 = q(v)1_\mathcal{A} \quad \forall\ v \in \mathcal{V}
		\end{equation}
		then $f$ extends uniquely to an $\mathbb{F}$-algebra homomorphism
		$F\!: \mathcal{G}(\mathcal{V},q) \to \mathcal{A}$, i.e.
		\begin{displaymath}
		\setlength\arraycolsep{2pt}
		\begin{array}{rcll}
			F(\alpha) &=& \alpha 1_\mathcal{A},	&\quad \forall\ \alpha \in \mathbb{F}, \\
			F(v) &=& f(v),  					&\quad \forall\ v \in \mathcal{V}, \\
			F(xy) &=& F(x)F(y),					\\
			F(x+y) &=& F(x)+F(y),				&\quad \forall\ x,y \in \mathcal{G}.
		\end{array}
		\end{displaymath}
		Furthermore, $\mathcal{G}$ is the unique associative $\mathbb{F}$-algebra with this property.
	\end{prop}
	\begin{proof}
		Any linear map $f\!: \mathcal{V} \to \mathcal{A}$ extends to a unique
		algebra homomorphism $\hat{f}\!: \mathcal{T}(\mathcal{V}) \to \mathcal{A}$
		defined by $\hat{f}(u \otimes v) := f(u)f(v)$ etc.
		Property \eqref{universality_property} implies that $\hat{f}=0$ on the
		ideal $\mathcal{I}_q(\mathcal{V})$ and so $\hat{f}$ descends to a well-defined map $F$ 
		on $\mathcal{G}(\mathcal{V},q)$ which has the required properties.
		Suppose now that $\mathcal{C}$ is an associative $\mathbb{F}$-algebra
		with unit and that $i\!: \mathcal{V} \hookrightarrow \mathcal{C}$ is an embedding
		with the property that any linear map $f\!: \mathcal{V} \to \mathcal{A}$
		with property \eqref{universality_property} extends uniquely to an
		algebra homomorphism $F\!: \mathcal{C} \to \mathcal{A}$. Then the isomorphism
		from $\mathcal{V} \subseteq \mathcal{G}$ to $i(\mathcal{V}) \subseteq \mathcal{C}$
		clearly induces an algebra isomorphism $\mathcal{G} \to \mathcal{C}$.
	\end{proof}

	So far we have not made any assumptions on the dimension of $\mathcal{V}$.
	We will come back to the infinite-dimensional case when discussing
	the more general Clifford algebra. 
	Here we will familiarize ourselves with the properties of quadratic forms
	on finite-dimensional spaces.
	For the remainder of this subsection
	we will therefore assume that $\dim \mathcal{V} = n < \infty$.

	\begin{defn} \label{def_o_basis}
		A basis $\{e_1,\ldots,e_n\}$ of $(\mathcal{V},q)$ is said to be \emph{orthogonal} or \emph{canonical}
		if $\beta_q(e_i,e_j) = 0$ for all $i \neq j$. 
		The basis is called \emph{orthonormal} if we also have that $q(e_i) \in \{-1,0,1\}$ for all $i$.
	\end{defn}
	
	We have a number of classical theorems regarding orthogonal bases.
	Proofs of these can be found e.g. in \cite{svensson}.
	
	\begin{thm} \label{thm_o_basis}
		If $\dim \mathcal{V} < \infty$ and \emph{char} $\mathbb{F} \neq 2$ then 
		there exists an orthogonal basis of $(\mathcal{V},q)$.
	\end{thm}
	
	\noindent
	Because this rather fundamental theorem breaks down for fields of characteristic
	two (such as $\mathbb{Z}_2)$, we will always assume that char $\mathbb{F} \neq 2$ when
	talking about geometric algebra. General fields and rings will be treated by
	the general Clifford algebra, however.
	
	
	\begin{thm}[Sylvester's Law of Inertia] \label{thm_sylvester}
		Assume that $\dim \mathcal{V} < \infty$ and $\mathbb{F}=\mathbb{R}$.
		If $E_1$ and $E_2$ are two orthogonal bases of $\mathcal{V}$ and
		\begin{displaymath}
		\setlength\arraycolsep{2pt}
		\begin{array}{lcl}
			E_i^+ &:=& \{ e \in E_i : q(e)>0 \}, \\
			E_i^- &:=& \{ e \in E_i : q(e)<0 \}, \\
			E_i^0 &:=& \{ e \in E_i : q(e)=0 \}
		\end{array}
		\end{displaymath}
		then
		\begin{displaymath}
		\setlength\arraycolsep{2pt}
		\begin{array}{c}
			|E_1^+| = |E_2^+|, \\
			|E_1^-| = |E_2^-|, \\
			\text{\emph{Span }} E_1^0 = \text{\emph{Span }} E_2^0.
		\end{array}
		\end{displaymath}
	\end{thm}
	
	\noindent
	This means that there is a unique \emph{signature} $(s,t,u) := (|E_i^+|,|E_i^-|,|E_i^0|)$ associated to $(\mathcal{V},q)$.
	For the complex case we have the following simpler result:
	
	\begin{thm} \label{thm_complex_sylvester}
		If $E_1$ and $E_2$ are orthogonal bases of $\mathcal{V}$ with $\mathbb{F}=\mathbb{C}$ and
		\begin{displaymath}
		\setlength\arraycolsep{2pt}
		\begin{array}{lcl}
			E_i^\times &:=& \{ e \in E_i : q(e) \neq 0 \}, \\
			E_i^0 &:=& \{ e \in E_i : q(e)=0 \}
		\end{array}
		\end{displaymath}
		then
		\begin{displaymath}
		\setlength\arraycolsep{2pt}
		\begin{array}{c}
			\text{\emph{Span }} E_1^\times = \text{\emph{Span }} E_2^\times, \\
			\text{\emph{Span }} E_1^0 = \text{\emph{Span }} E_2^0.
		\end{array}
		\end{displaymath}
		If $E_i^0 = \varnothing$ ($q$ is nondegenerate) then there exists a basis
		$E$ with $q(e)=1 \ \forall\ e \in E$.
	\end{thm}
	
	\noindent
	From the above follows that we can talk about \emph{the} signature of a quadratic form 
	or a metric without ambiguity.
	We use the short-hand notation $\mathbb{R}^{s,t,u}$ to denote the $(s+t+u)$-dimensional real vector space
	with a quadratic form of signature $(s,t,u)$, while $\mathbb{C}^n$ is
	understood to be the complex $n$-dimensional space with a nondegenerate quadratic form.
	When $u=0$ or $t=u=0$ we may simply write $\mathbb{R}^{s,t}$ or $\mathbb{R}^{s}$.
	A space of type $\mathbb{R}^{n,0}$ is called \emph{euclidean} and $\mathbb{R}^{0,n}$
	\emph{anti-euclidean}, while the spaces $\mathbb{R}^{1,n}$ ($\mathbb{R}^{n,1}$) are called
	(\emph{anti-})\emph{lorentzian}.
	Within real and complex spaces we can always find bases that are orthonormal.

	\begin{rem}
		The general condition for orthonormal bases to exist 
		is that the field $\mathbb{F}$ is a so called
		\emph{spin field}. This means that every $\alpha \in \mathbb{F}$ 
		can be written as $\alpha = \beta^2$ or $-\beta^2$ for some $\beta \in \mathbb{F}$.
		The fields $\mathbb{R}$, $\mathbb{C}$ and $\mathbb{Z}_p$ for $p$ a prime with $p \equiv 3 \pmod{4}$, are spin,
		but e.g. $\mathbb{Q}$ is not.
	\end{rem}
	
	Consider now the geometric algebra $\mathcal{G}$ over a real or complex space $\mathcal{V}$.
	If we pick an orthonormal basis $E = \{e_1,\ldots,e_n\}$ of $\mathcal{V}$ it follows from Definition \ref{def_g} and \eqref{generator_relations}
	that $\mathcal{G}$ is the free associative algebra generated by $E$ modulo the relations
	\begin{equation}
		e_i^2 = q(e_i) = \pm 1 \ \textrm{or} \ 0 \qquad \textrm{and} \qquad e_i e_j = - e_j e_i,\ i \neq j.
	\end{equation}
	We also observe that $\mathcal{G}$ is spanned by $\{E_{i_1 i_2 \ldots i_k}\}_{i_1<i_2< \ldots <i_k}$, where $E_{i_1 i_2 \ldots i_k} := \\ e_{i_1} e_{i_2} \ldots e_{i_k}$.
	Thus, one can view $\mathcal{G}$ as vector space isomorphic to $\wedge^* \mathcal{V}$, the exterior algebra of $\mathcal{V}$.
	This is a description of geometric algebra (Clifford algebra) which 
	may be more familiar to e.g. physicists.
	
	\begin{rem}
		If we take $q=0$ we actually obtain an \emph{algebra isomorphism} $\mathcal{G} \cong \wedge^* \mathcal{V}$.
		In this case $\mathcal{G}$ is called a \emph{Grassmann algebra}.
	\end{rem}
	
	One element in $\mathcal{G}$ deserves special attention, namely the so called \emph{pseudoscalar} 
	\begin{equation}
		I := e_1 e_2 \ldots e_n.
	\end{equation}
	Note that this definition is basis independent up to orientation when $q$ is nondegenerate.
	Indeed, let $\{Re_1,\ldots,Re_n\}$ be another orthonormal basis with the same orientation, where $R \in O(\mathcal{V},q)$,
	the group of linear transformations which leave $q$ invariant\footnote{The details
	surrounding such transformations will be discussed in Section \ref{sec_groups}}.
	Then $Re_1 Re_2 \ldots Re_n = \sum_{\pi \in S_n} \textrm{sign}(\pi)\ R_{\pi(1)1} \ldots R_{\pi(n)n} \cdot \\ e_1 e_2 \ldots e_n = \det R \ e_1 e_2 \ldots e_n = I$ 
	due to the anticommutativity of the $e_i$:s.
	Note that, by selecting a certain pseudoscalar for $\mathcal{G}$ we also 
	impose a certain orientation on $\mathcal{V}$. There is no such thing as an absolute orientation;
	instead all statements concerning orientation will be made relative to the chosen one.
	
	The square of the pseudoscalar is given by (and gives information about) the signature and dimension of $(\mathcal{V},q)$.
	For $\mathcal{G}(\mathbb{R}^{s,t,u})$ we have that
	\begin{equation}
		I^2 = (-1)^{\frac{1}{2}n(n-1) + t} \delta_{u,0}, \quad \textrm{where} \ n=s+t+u.
	\end{equation}
	We say that $\mathcal{G}$ is \emph{degenerate} if the quadratic form 
	is degenerate, or equivalently if $I^2=0$.
	For odd $n$, $I$ commutes with all elements in $\mathcal{G}$ 
	and the center of $\mathcal{G}$ is $Z(\mathcal{G}) = \ \textrm{Span}_\mathbb{F} \ \{1,I\}$.
	For even $n$, the center consists of the scalars $\mathbb{F}$ only.
	
\subsection{Combinatorial Clifford algebra $\cl(X,R,r)$}
	
	We now take a temporary step away from the comfort of fields and vector spaces
	and instead consider the purely algebraic features of geometric algebra that were
	uncovered in the previous subsection.
	Note that we could roughly write 
	\begin{equation}
		\mathcal{G}(\mathcal{V}) = \textrm{Span}_\mathbb{F} \ \{E_A\}_{A \subseteq \{1,2,\ldots,n\}}
	\end{equation}
	for an $n$-dimensional space $\mathcal{V}$ over $\mathbb{F}$, 
	and that the geometric product of these basis elements behaves as
	\begin{equation}
		E_A E_B = \tau(A,B) \ E_{A \bigtriangleup B}, \quad \textrm{where} \quad \tau(A,B)=1,-1 \ \textrm{or}\ 0,
	\end{equation}
	and $A \symdiff B := (A \cup B) \!\smallsetminus\! (A \cap B)$ is 
	the symmetric difference between the sets $A$ and $B$.
	Motivated by this we consider the following generalization.
	
	\begin{defn} \label{def_cl}
		Let $X$ be a finite set and $R$ a commutative ring with unit. Let
		$r\!: X \to R$ be some function which is to be thought of as a \emph{signature} on $X$.
		The \emph{Clifford algebra} over $(X,R,r)$ is defined as the set
		\begin{displaymath}
			\cl(X,R,r) := \bigoplus_{\mathscr{P}(X)} R,
		\end{displaymath}
		i.e. the free $R$-module generated by $\mathscr{P}(X)$, the set of all subsets of $X$.
		We may use the shorter notation $\cl(X)$, or just $\cl$, when 
		the current choice of $X$, $R$ and $r$ is clear from the context.
		We call $R$ the \emph{scalars} of $\cl$.
	\end{defn}
	
	\begin{exmp}
		A typical element of $\cl(\{x,y,z\},\mathbb{Z},r)$ could for example look like
		\begin{equation}
			5 \varnothing + 3 \{x\} + 2 \{y\} - \{x,y\} + 12 \{x,y,z\}.
		\end{equation}
	\end{exmp}
	
	We have not yet defined a product on $\cl$.
	In addition to being $R$-bilinear and associative, we would like the product to satisfy
	$\{x\}^2 = r(x) \varnothing$ for $x \in X$, 
	$\{x\}\{y\} = -\{y\}\{x\}$ for $x \neq y \in X$ and 
	$\varnothing A = A \varnothing = A$	for all $A \in \mathscr{P}(X)$.
	In order to arrive at such a product we make use of the following

	\begin{lem} \label{lem_finite_tau}
		There exists a map $\tau\!: \mathscr{P}(X) \times \mathscr{P}(X) \to R$ such that
		\begin{displaymath}
		\begin{array}{rl}
			i)  & \tau(\{x\},\{x\}) = r(x) \quad \forall\ x \in X, \\[5pt]
			ii) & \tau(\{x\},\{y\}) = -\tau(\{y\},\{x\}) \quad \forall\ x,y \in X : x \neq y, \\[5pt]
			iii)& \tau(\varnothing,A) = \tau(A,\varnothing) = 1 \quad \forall\ A \in \mathscr{P}(X), \\[5pt]
			iv) & \tau(A,B) \tau(A \symdiff B,C) = \tau(A,B \symdiff C) \tau(B,C) \quad \forall\ A,B,C \in \mathscr{P}(X), \\[5pt]
			v)  & \tau(A,B) \in \{-1,1\} \quad \textrm{if} \quad A \cap B = \varnothing.
		\end{array}
		\end{displaymath}
	\end{lem}
	\begin{proof}
		We proceed by induction on the cardinality $|X|$ of $X$.
		For $X=\varnothing$ the lemma is trivial, so let $z \in X$
		and assume the lemma holds for $Y := X \!\smallsetminus\! \{z\}$.
		Hence, there is a $\tau'\!: \mathscr{P}(Y) \times \mathscr{P}(Y) \to R$ 
		which has the properties (\emph{i})-(\emph{v}) above. If $A \subseteq Y$ we
		write $A' = A \cup \{z\}$ and, for $A,B$ in $\mathscr{P}(Y)$ we extend
		$\tau'$ to $\tau\!: \mathscr{P}(X) \times \mathscr{P}(X) \to R$ in the
		following way:
		\begin{displaymath}
		\setlength\arraycolsep{2pt}
		\begin{array}{rcl}
			\tau(A,B)	&:=& \tau'(A,B) \\[3pt]
			\tau(A',B)	&:=& (-1)^{|B|} \tau'(A,B) \\[3pt]
			\tau(A,B')	&:=& \tau'(A,B) \\[3pt]
			\tau(A',B') &:=& r(z) (-1)^{|B|} \tau'(A,B)
		\end{array}
		\end{displaymath}
		Now it is straightforward to verify that (\emph{i})-(\emph{v}) holds for
		$\tau$, which completes the proof.
	\end{proof}

	\begin{defn} \label{def_cl_prod}
		Define the \emph{Clifford product}
		\begin{displaymath}
		\setlength\arraycolsep{2pt}
		\begin{array}{ccc}
			\cl(X) \times \cl(X) & \to & \cl(X) \\
			(A,B) & \mapsto & AB
		\end{array}
		\end{displaymath}
		by taking $AB := \tau(A,B) A \symdiff B$ for $A,B \in \mathscr{P}(X)$ and extending linearly.
		We choose to use the $\tau$ which is constructed as in the proof of Lemma \ref{lem_finite_tau} 
		by consecutively adding elements from the set $X$.
		A unique such $\tau$ may only be selected after imposing a certain order (orientation) on the set $X$.
	\end{defn}
	
	Using Lemma \ref{lem_finite_tau} one easily verifies that this product has all
	the properties that we asked for above. For example, in order to verify associativity
	we note that 
	\begin{equation}
		A(BC) = A\big(\tau(B,C) B \symdiff C\big) = \tau(A, B \symdiff C) \tau(B,C) A \symdiff (B \symdiff C),
	\end{equation}
	while 
	\begin{equation}
		(AB)C = \tau(A,B) (A \symdiff B)C = \tau(A,B) \tau(A \symdiff B,C) (A \symdiff B) \symdiff C.
	\end{equation}
	Associativity now follows from (\emph{iv}) and the associativity of the symmetric difference.
	As is expected from the analogy with $\mathcal{G}$, we also have the property that
	different basis elements of $\cl$ commute up to a sign.

	\begin{prop} \label{prop_reverse}
		If $A,B \in \mathscr{P}(X)$ then
		\begin{displaymath}
			AB = (-1)^{\frac{1}{2}|A|(|A|-1) \ + \ \frac{1}{2}|B|(|B|-1) \ + \ \frac{1}{2}|A \bigtriangleup B|(|A \bigtriangleup B|-1)} BA.
		\end{displaymath}
	\end{prop}
	\begin{proof}
		By the property (\emph{v}) in Lemma \ref{lem_finite_tau} it is sufficient
		to prove this for $A = \{a_1\} \{a_2\} \ldots \{a_k\}$, $B = \{b_1\} \{b_2\} \ldots \{b_l\}$,
		where $a_i$ are disjoint elements in $X$ and likewise for $b_i$.
		If $A$ and $B$ have $m$ elements in common then $AB = (-1)^{(k-m)l\ +\ m(l-1)}BA = (-1)^{kl-m}BA$
		by property (\emph{ii}).
		But then we are done, since $\frac{1}{2} \big(-k(k-1) - l(l-1) + (k+l-2m)(k+l-2m-1) \big) \equiv kl + m \pmod{2}$.
	\end{proof}
	
	We are now ready to make the formal connection between $\mathcal{G}$ and $\cl$.
	Let $(\mathcal{V},q)$ be a vector space over $\mathbb{F}$ with a quadratic form.
	Pick an orthogonal basis $E = \{e_1,\ldots,e_n\}$ of $\mathcal{V}$ 
	and consider the Clifford algebra $Cl(E,\mathbb{F},q|_E)$.
	Define $f\!: \mathcal{V} \to \cl$ by $f(e_i) := \{e_i\}$ for $i=1,\ldots,n$ and extend linearly.
	We then have 
	\begin{equation}
	\setlength\arraycolsep{2pt}
	\begin{array}{rl}
		f(v)^2 	&= f(\sum_i v_i e_i) f(\sum_j v_j e_j) = \sum_{i,j} v_i v_j f(e_i)f(e_j) \\[10pt]
				&= \sum_{i,j} v_i v_j \{e_i\}\{e_j\} = \sum_i v_i^2 \{e_i\}^2 \\[10pt]
				&= \sum_i v_i^2 q|_E(e_i) \varnothing = \sum_i v_i^2 q(e_i) \varnothing = q(\sum_i v_i e_i) \varnothing = q(v) \varnothing.
	\end{array}
	\end{equation}
	By Proposition \ref{prop_universality}, $f$ extends uniquely to a homomorphism $F\!: \mathcal{G} \to \cl$.
	Since $\dim \mathcal{G} = \dim \cl = 2^n$ and $F$ is easily seen to be surjective 
	from the property (\emph{v}),
	we arrive at an isomorphism
	\begin{equation} \label{g_isomorphic_to_cl}
		\mathcal{G}(\mathcal{V},q) \cong \cl(E,\mathbb{F},q|_E).
	\end{equation}
	We make this equivalence between $\mathcal{G}$ and $\cl$ even more transparent by 
	suppressing the unit $\varnothing$ in expressions and writing simply $e$ instead
	of $\{e\}$ for singletons $e \in E$. For example, with an orthonormal
	basis $\{e_1,e_2,e_3\}$ in $\mathbb{R}^3$, both $\mathcal{G}$ and $\cl$ are then spanned by
	\begin{equation}
		\{ 1, \ e_1, e_2, e_3, \ e_1 e_2, e_1 e_3, e_2 e_3, \ e_1 e_2 e_3 \}.
	\end{equation}
	
	There is a natural grade structure on $Cl$ given by the cardinality of the subsets of $X$.
	Consider the following

	\begin{defn} \label{def_grade_space}
		The \emph{subspace of $k$-vectors} in $\cl$, or the \emph{grade-$k$ part} of $\cl$, is defined by
		\begin{displaymath}
			\cl^k(X,R,r) := \bigoplus_{A \in \mathscr{P}(X)\ :\ |A| = k} R.
		\end{displaymath}
		Of special importance are the \emph{even} and \emph{odd subspaces},
		\begin{displaymath}
			\cl^{\pm}(X,R,r) := \bigoplus_{k \ \text{is}\ \substack{\text{\tiny even}\\ \text{\tiny odd}}} \cl^k(X,R,r).
		\end{displaymath}
		This notation carries over to the corresponding subspaces of $\mathcal{G}$ 
		and we write $\mathcal{G}^k$, $\mathcal{G}^\pm$ etc. where for example 
		$\mathcal{G}^0 = \mathbb{F}$ and $\mathcal{G}^1 = \mathcal{V}$.
		The elements of $\mathcal{G}^2$ are also called \emph{bivectors}, while
		arbitrary elements of $\mathcal{G}$ are traditionally called \emph{multivectors}.
	\end{defn}
	
	\noindent
	We then have a split of $\cl$ into graded subspaces as
	\begin{equation} \label{cl_grades}
	\setlength\arraycolsep{2pt}
	\begin{array}{rcl}
		\cl(X)
			&=& \cl^{+} \!\oplus \cl^{-} \\[5pt]
			&=& \cl^0 \oplus \cl^1 \oplus \cl^2 \oplus \ldots \oplus \cl^{|X|}.
	\end{array}
	\end{equation}
	Note that, under the Clifford product, $\cl^{\pm} \cdot \cl^{\pm} \subseteq \cl^{+}$ and $\cl^{\pm} \cdot \cl^{\mp} \subseteq \cl^{-}$.
	Hence, the even-grade elements $\cl^{+}$ form a subalgebra of $\cl$.
	
	In $\cl(X,R,r)$ we have the possibility of defining a unique pseudoscalar
	independently of the signature $r$, namely the set $X$ itself.
	Note, however, that it can only be normalized if $X^2 = \tau(X,X) \in R$ is invertible,
	which requires that $r$ is nondegenerate. We will almost always talk about
	pseudoscalars in the setting of nondegenerate vector spaces, so this will not be a problem.
	
\subsection{Standard operations}
	
	
	A key feature of Clifford algebras is that they contain a surprisingly large
	amount of structure. In order to really be able to harness the power of this structure
	we need to introduce powerful notation. 
	Most of the following definitions will be made on $\cl$ for simplicity,
	but because of the equivalence between $\mathcal{G}$ and $\cl$ they
	carry over to $\mathcal{G}$ in a straightforward manner.
	
	We will find it convenient to introduce the notation that for any proposition $P$, 
	$(P)$ will denote the number $1$ if $P$ is true and $0$ if $P$ is false.
	
	\begin{defn} \label{def_operations}
		For $A,B \in \mathscr{P}(X)$ define
		\begin{displaymath}
		\setlength\arraycolsep{2pt}
		\begin{array}{ccll}
			A \wedge B &:=& (A \cap B = \varnothing) \ AB	& \quad\textrm{\emph{outer product}} \\[5pt]
			A \liprod B &:=& (A \subseteq B) \ AB			& \quad\textrm{\emph{left inner product}} \\[5pt]
			A \riprod B &:=& (A \supseteq B) \ AB			& \quad\textrm{\emph{right inner product}} \\[5pt]
			A * B &:=& (A = B) \ AB							& \quad\textrm{\emph{scalar product}} \\[5pt]
			\langle A \rangle_n &:=& (|A|=n) \ A			& \quad\textrm{\emph{projection on grade $n$}} \\[5pt]
			A^\star &:=& (-1)^{|A|} \ A						& \quad\textrm{\emph{grade involution}} \\[3pt]
			A^\dagger &:=& (-1)^{\binom{|A|}{2}} \ A		& \quad\textrm{\emph{reversion}}
		\end{array}
		\end{displaymath}
		and extend linearly to $\cl(X,R,r)$.
	\end{defn}
	
	\noindent
	The grade involution is also called the (\emph{first}) \emph{main involution}. It has the property
	\begin{equation} \label{grade_inv_property}
		(xy)^\star = x^\star y^\star, \quad v^\star = -v
	\end{equation}
	for all $x,y \in Cl(X)$ and $v \in Cl^1(X)$, as is easily verified 
	by expanding in linear combinations of elements in $\mathscr{P}(X)$
	and using that $|A \symdiff B| \equiv |A| + |B|\ (\textrm{mod}\ 2)$.
	The reversion earns its name from the property
	\begin{equation} \label{reversion_property}
		(xy)^\dagger = y^\dagger x^\dagger, \quad v^\dagger = v,
	\end{equation}
	and it is sometimes called the \emph{second main involution} or the \emph{principal antiautomorphism}.
	This reversing behaviour follows directly from Proposition \ref{prop_reverse}.
	We will find it convenient to have a name for the composition of these two involutions.
	Hence, we define the \emph{Clifford conjugate} $x^\cliffconj$ of $x \in Cl(X)$ by $x^\cliffconj := x^{\star\dagger}$
	and observe the property
	\begin{equation} \label{cliffconj_property}
		(xy)^\cliffconj = y^\cliffconj x^\cliffconj, \quad v^\cliffconj = -v.
	\end{equation}
	Note that all the above involutions act by changing sign on some of the graded subspaces.
	We can define general involutions of that kind which will come in handy later.
	
	\begin{defn} \label{def_grade_involutions}
		For $A \in \mathscr{P}(X)$ define
		\begin{displaymath}
		\setlength\arraycolsep{2pt}
		\begin{array}{rcl}
			[A]					&:=& (-1)^{(A \neq \varnothing)} A, \\[5pt]
			[A]_{p,q,\ldots,r} 	&:=& (-1)^{(|A| = p,q,\ldots,\ \textrm{or}\ r)} A.
		\end{array}
		\end{displaymath}
		and extend linearly to $\cl(X,R,r)$.
	\end{defn}
	
	\noindent
	We summarize the action of these involutions in Table \ref{table_involutions}.
	Note the periodicity.
	
	\begin{table}[ht]
		\begin{displaymath}
		\begin{array}{c|cccccccc}
				& \cl^0 & \cl^1 & \cl^2 & \cl^3 & \cl^4 & \cl^5 & \cl^6 & \cl^7 \\
			\hline
			\\[-1.5ex]
			\star				& + & - & + & - & + & - & + & - \\
			\dagger				& + & + & - & - & + & + & - & - \\
			\cliffconj			& + & - & - & + & + & - & - & + \\
			\lbrack\ \rbrack	& + & - & - & - & - & - & - & -
		\end{array}
		\end{displaymath}
		\caption{The action of involutions on graded subspaces of $\cl$. \label{table_involutions}}
	\end{table}
	
	The scalar product has the symmetric property $x * y = y * x$ for all $x,y \in \cl$.
	Therefore, it forms a symmetric bilinear map $\cl \times \cl \to R$
	which is degenerate if and only if $\cl$ (i.e. the signature $r$) is degenerate.
	This map coincides with the bilinear form $\beta_q$ when restricted to $\mathcal{V} \subseteq \mathcal{G}(\mathcal{V},q)$.
	Note also that subspaces of different grade are orthogonal with respect to the scalar product.
	
	Another product that is often seen in the context of geometric algebra is
	the \emph{inner product}, defined by $A \iprod B := (A \subseteq B \ \textrm{or} \ A \supseteq B)\ AB = A \liprod B + A \riprod B - A * B$.
	We will stick to the left and right inner products, however, 
	because they admit a simpler handling of grades, something which is illustrated\footnote{The
	corresponding identities with $\!\iprod\!$ instead of $\llcorner,\lrcorner$ need to be supplied with grade restrictions.}
	by the following
	
	\begin{prop} \label{prop_trix}
		For all $x,y,z \in \cl(X)$ we have
		\begin{displaymath}
		\setlength\arraycolsep{2pt}
		\begin{array}{ccc}
			x \wedge (y \wedge z) &=& (x \wedge y) \wedge z, \\[2pt]
			x \liprod (y \riprod z) &=& (x \liprod y) \riprod z, \\[2pt]
			x \liprod (y \liprod z) &=& (x \wedge y) \liprod z, \\[2pt]
			x * (y \liprod z) &=& (x \wedge y) * z.
		\end{array}
		\end{displaymath}
	\end{prop}
	\begin{proof}
		This follows directly from Definition \ref{def_operations} and basic set logic.
		For example, taking $A,B,C \in \mathscr{P}(X)$ we have
		\begin{equation}
		\setlength\arraycolsep{2pt}
		\begin{array}{rl}
			A \liprod (B \liprod C)
				&= (B \subseteq C)(A \subseteq B \symdiff C)ABC \\[5pt]
				&= (B \subseteq C \ \textrm{and} \ A \subseteq C \!\smallsetminus\! B)ABC \\[5pt]
				&= (A \cap B = \varnothing \ \textrm{and} \ A \cup B \subseteq C)ABC \\[5pt]
				&= (A \cap B = \varnothing)(A \symdiff B \subseteq C)ABC \\[5pt]
				&= (A \wedge B) \liprod C.
		\end{array}
		\end{equation}
		The other identities are proven in an equally simple way.
	\end{proof}
	
	To work efficiently with geometric algebra it is crucial to understand
	how vectors behave under these operations.
	
	\begin{prop} \label{prop_trix2}
		For all $x,y \in \cl(X)$ and $v \in \cl^1(X)$ we have
		\begin{displaymath}
		\setlength\arraycolsep{2pt}
		\begin{array}{rcl}
			vx				&=& v \liprod x + v \wedge x, \\[3pt]
			v \liprod x		&=& \frac{1}{2}(vx - x^\star v) = - x^\star \!\riprod v, \\[3pt]
			v \wedge x		&=& \frac{1}{2}(vx + x^\star v) = \phantom{-} x^\star \!\wedge v, \\[3pt]
			v \liprod (xy)	&=& (v \liprod x)y + x^\star (v \liprod y).
		\end{array}
		\end{displaymath}
	\end{prop}
	
	\noindent
	The first three identities are shown simply by using linearity and set relations,
	while the fourth follows immediately from the second.
	Note that for 1-vectors $u,v \in \cl^1$ we have the basic relations
	\begin{equation}
		u \liprod v = v \liprod u = u \riprod v = u * v = \frac{1}{2}(uv + vu)
	\end{equation}
	and
	\begin{equation} \label{vector_wedge}
		u \wedge v = - v \wedge u  = \frac{1}{2}(uv - vu).
	\end{equation}
	
	
	It is often useful to expand the various products and involutions 
	in terms of the grades involved.
	The following identities are left as exercises.
	
	\begin{prop} \label{prop_grade_def}
		For all $x,y \in Cl(X)$ we have
		\begin{displaymath}
		\setlength\arraycolsep{2pt}
		\begin{array}{ccl}
			x \wedge y		&=& \sum_{n,m \geq 0} \big\langle \langle x \rangle_n \langle y \rangle_m \big\rangle_{n+m}, \\[5pt]
			x \liprod y		&=& \sum_{0 \leq n \leq m} \big\langle \langle x \rangle_n \langle y \rangle_m \big\rangle_{m-n}, \\[5pt]
			x \riprod y		&=& \sum_{n \geq m \geq 0} \big\langle \langle x \rangle_n \langle y \rangle_m \big\rangle_{n-m}, \\[5pt]
			x \iprod y		&=& \sum_{n,m \geq 0} \big\langle \langle x \rangle_n \langle y \rangle_m \big\rangle_{|n-m|}, \\[5pt]
			x * y			&=& \langle xy \rangle_0, \\[5pt]
			x^\star			&=& \sum_{n \geq 0} (-1)^n \langle x \rangle_n, \\[3pt]
			x^\dagger		&=& \sum_{n \geq 0} (-1)^{\binom{n}{2}} \langle x \rangle_n.
		\end{array}
		\end{displaymath}
	\end{prop}
	
	In the general setting of a Clifford algebra with scalars in a ring $R$, we need
	to be careful about the notion of linear (in-)dependence. A subset
	$\{x_1,x_2,\ldots,x_m\}$ of $\cl$ is called \emph{linearly dependent} iff there
	exist $r_1,\ldots,r_m \in R$, not all zero, such that
	\begin{equation}
		r_1 x_1 + r_2 x_2 + \ldots + r_m x_m = 0.
	\end{equation}
	Note that a single nonzero 1-vector could be linearly dependent in this context.
	We will prove an important theorem concerning linear dependence where we need the following
	
	\begin{lem} \label{lem_blade_det}
		If $u_1,u_2,\ldots,u_k$ and $v_1,v_2,\ldots,v_k$ are 1-vectors then
		\begin{displaymath}
			(u_1 \wedge u_2 \wedge \cdots \wedge u_k) * (v_k \wedge v_{k-1} \wedge \cdots \wedge v_1)
				= \det\ [u_i*v_j]_{1 \leq i,j \leq k}.
		\end{displaymath}
	\end{lem}
	\begin{proof}
		Since both sides of the expression are multilinear and alternating in both
		the $u_i$:s and the $v_i$:s, we need only consider ordered disjoint elements 
		$\{e_i\}$ in the basis of singleton sets in $X$.
		Both sides are zero, except in the case
		\begin{equation}
		\setlength\arraycolsep{2pt}
		\begin{array}{l}
			(e_{i_1} e_{i_2} \ldots e_{i_k}) * (e_{i_k} e_{i_{k-1}} \ldots e_{i_1}) = \\[5pt]
			\qquad = r(e_{i_1}) r(e_{i_2}) \ldots r(e_{i_k}) = \det\ [r(e_{i_p}) \delta_{p,q}]_{1 \leq p,q \leq k} \\[5pt]
			\qquad = \det\ [e_{i_p}*e_{i_q}]_{1 \leq p,q \leq k},
		\end{array}
		\end{equation}
		so we are done.
	\end{proof}
	
	\begin{thm} \label{thm_linear_indep}
		The 1-vectors $\{x_1,x_2,\ldots,x_m\}$ are linearly independent iff
		the m-vector $\{x_1 \wedge x_2 \wedge \cdots \wedge x_m\}$ is linearly independent.
	\end{thm}

	\begin{proof}
		Assume that $r_1 x_1 + \ldots + r_m x_m = 0$, where, say, $r_1 \neq 0$.
		Then
		\begin{equation}
		\setlength\arraycolsep{2pt}
		\begin{array}{l}
			r_1 (x_1 \wedge \cdots \wedge x_m) = (r_1 x_1) \wedge x_2 \wedge \cdots \wedge x_m \\[5pt]
			\qquad	= (r_1 x_1 + \ldots + r_m x_m) \wedge x_2 \wedge \cdots \wedge x_m = 0,
		\end{array}
		\end{equation}
		since $x_i \wedge x_i = 0$.
		
		Conversely, assume that $rX = 0$ for
		$r \neq 0$ in $R$ and $X = x_1 \wedge \cdots \wedge x_m$.
		We will use the basis minor theorem for arbitrary rings which
		can be found in the appendix.
		Assume that $x_j = x_{1j} e_1 + \ldots + x_{nj} e_n$, where $x_{ij} \in R$
		and $e_i \in X$ are basis elements such that $e_i^2 = 1$. This assumption on
		the signature is no loss in generality, since this theorem only concerns
		the exterior algebra associated to the outer product. It will only serve to
		simplify our reasoning below. Collect the coordinates in a matrix
		\begin{equation}
			A := \left[
			\begin{array}{cccc}
				rx_{11} & x_{12} & \cdots & x_{1m} \\
				rx_{21} & x_{22} & \cdots & x_{2m} \\
				\vdots  & \vdots & 		  & \vdots \\
				rx_{n1} & x_{n2} & \cdots & x_{nm} \\
			\end{array}
			\right] \in R^{n \times m},\ m \leq n
		\end{equation}
		and note that we can expand $rX$ in a grade-$m$ basis as
		\begin{equation}
			rX = \sum_{E \subseteq X : |E| = m} (rX * E^\dagger) E = \sum_{E \subseteq X : |E| = m} (\det A_{E,\{1,\ldots,m\}}) E,
		\end{equation}
		where we used Lemma \ref{lem_blade_det}.
		We find that the determinant of each $m \times m$ minor of $A$ is zero.
		
		Now, let $k$ be the rank of $A$. Then we must have $k<m$, and if $k=0$ then
		$rx_{i1}=0$ and $x_{ij}=0$ for all $i=1,\ldots,n$, $j>1$.
		But that would mean that $\{x_1,\ldots,x_m\}$ are linearly dependent.
		Therefore we assume that $k>0$ and, without loss of generality, that
		\begin{equation}
			d := \det \left[
			\begin{array}{cccc}
				rx_{11} & x_{12} & \cdots & x_{1k} \\
				\vdots  & \vdots & 		  & \vdots \\
				rx_{k1} & x_{k2} & \cdots & x_{kk} \\
			\end{array}
			\right] \neq 0.
		\end{equation}
		By the basis minor theorem there exist $r_1,\ldots,r_k \in R$
		such that
		\begin{equation}
			r_1 rx_1 + r_2 x_2 + \ldots + r_k x_k + dx_m = 0.
		\end{equation}
		Hence, $\{x_1,\ldots,x_m\}$ are linearly dependent.
	\end{proof}
	
	For our final set of operations,
	we will consider a nondegenerate geometric algebra $\mathcal{G}$ with
	pseudoscalar $I$. The nondegeneracy implies that there exists a natural
	duality between the inner and outer products.
	
	\begin{defn} \label{def_dual}
		We define the \emph{dual} of $x \in \mathcal{G}$ by $x\dual := xI^{-1}$.
		The \emph{dual outer product} or \emph{meet} $\vee$ is defined such that the diagram
		\begin{displaymath}
		\setlength\arraycolsep{2pt}
		\begin{array}{rcccccl}
						& \mathcal{G} 	& \times 	& \mathcal{G} 	& \xrightarrow{\vee} 	& \mathcal{G} \\
			(\ )\dual 	& \downarrow 	& 			& \downarrow	&						& \downarrow	& (\ )\dual \\
						& \mathcal{G} 	& \times	& \mathcal{G}	& \xrightarrow{\wedge}	& \mathcal{G}
		\end{array}
		\end{displaymath}
		commutes, i.e. $(x \vee y)\dual := x\dual \wedge y\dual \ \Rightarrow\ x \vee y = ((xI^{-1}) \wedge (yI^{-1}))I$.
	\end{defn}
	
	\begin{rem}
		In $\cl(X)$, the corresponding dual of $A \in \mathscr{P}(X)$ is 
		$A\dual = AX^{-1} = \tau(X,X)^{-1} \tau(A,X) A \symdiff X \propto A^c$, 
		the complement of the set $A$.
		Hence, we really find that the dual is the linearization of a 
		sign (or orientation) -respecting complement. This motivates our choice of notation.
	\end{rem}
	
	\begin{prop} \label{prop_dual}
		For all $x,y \in G$ we have
		\begin{equation} \label{dual_products}
		\setlength\arraycolsep{2pt}
		\begin{array}{ccc}
			x \liprod y\dual &=& (x \wedge y)\dual, \\
			x \wedge y\dual &=& (x \liprod y)\dual. \\
		\end{array}
		\end{equation}
	\end{prop}
	\begin{proof}
		Using Proposition \ref{prop_trix} and the fact that $xI = x \liprod I$, we obtain
		\begin{equation}
			x \liprod (yI^{-1}) = x \liprod (y \liprod I^{-1}) = (x \wedge y) \liprod I^{-1} = (x \wedge y)I^{-1},
		\end{equation}
		and from this follows also the second identity
		\begin{equation}
			(x \wedge y\dual)I^{-1}I = (x \liprod y^\mathbf{cc})I = (x \liprod y)I^{-2}I.
			\qedhere
		\end{equation}
	\end{proof}
	
	\noindent
	It is instructive to compare these results with those in the language of differential
	forms and Hodge duality, which are completely equivalent. In that setting one
	often starts with an outer product and then uses a metric to define a dual.
	The inner product is then defined from the outer product and dual according to \eqref{dual_products}.
	
\subsection{Vector space geometry}
	
	We will now leave the general setting of Clifford algebra for a moment and instead
	focus on the geometric properties of $\mathcal{G}$ and its newly
	defined operations.
	
	\begin{defn} \label{def_blades}
		A \emph{blade} is an outer product of 1-vectors. We define the following:
		\begin{displaymath}
		\setlength\arraycolsep{2pt}
		\begin{array}{lcll}
			\mathcal{B}_k &:=& \{ v_1 \wedge v_2 \wedge \cdots \wedge v_k \in \mathcal{G} : v_i \in \mathcal{V} \}
				& \quad\textrm{\emph{the set of $k$-blades}} \\[5pt]
			\mathcal{B} &:=& \bigcup_{k=0}^{\infty} \mathcal{B}_k
				& \quad\textrm{\emph{the set of all blades}} \\[5pt]
			\mathcal{B}^* &:=& \mathcal{B} \!\smallsetminus\! \{0\}
				& \quad\textrm{\emph{the nonzero blades}} \\[5pt]
			\mathcal{B}^\times &:=& \{ B \in \mathcal{B} : B^2 \neq 0\}
				& \quad\textrm{\emph{the invertible blades}}
		\end{array}
		\end{displaymath}
		The \emph{basis blades} associated to an orthogonal basis $E = \{e_i\}_{i=1}^{\dim \mathcal{V}}$ is 
		the basis of $\mathcal{G}$ generated by $E$, i.e.
		\begin{displaymath}
			\mathcal{B}_E := \{ e_{i_1} \wedge e_{i_2} \wedge \cdots \wedge e_{i_k} \in \mathcal{G} : i_1 < i_2 < \ldots < i_k \} \simeq \mathscr{P}(E)\ \textrm{in}\ \cl.
		\end{displaymath}
	\end{defn}
	
	\noindent
	We include the unit 1 among the blades and call it the $0$-blade.
	Note that $\mathcal{B}_k \subseteq \mathcal{G}^k$ and that by applying Proposition \ref{prop_trix2}
	recursively we can expand a blade as a sum of geometric products,
	\begin{equation} \label{blade_expansion}
	\setlength\arraycolsep{2pt}
	\begin{array}{rcl}
		a_1 \wedge a_2 \wedge \cdots \wedge a_k 
			&=& \frac{1}{2} a_1(a_2 \wedge \cdots \wedge a_k) + \frac{1}{2} (-1)^{k-1} (a_2 \wedge \cdots \wedge a_k)a_1 \\[5pt]
			&=& \ldots = \frac{1}{k!} \sum_{\pi \in S_k} \textrm{sign}(\pi)\ a_{\pi(1)} a_{\pi(2)} \ldots a_{\pi(k)}.
	\end{array}
	\end{equation}
	This expression is clearly similar to a determinant, except that this is a product of \emph{vectors}
	instead of scalars.
	
	The key property of blades is that they represent linear subspaces of $\mathcal{V}$.
	This is made precise by the following 
	
%
	
	\begin{prop} \label{prop_blade_subspace}
		If $A = a_1 \wedge a_2 \wedge \cdots \wedge a_k \neq 0$ is a 
		nonzero $k$-blade and $a \in \mathcal{V}$ then
		\begin{displaymath}
			a \wedge A = 0 \quad \Leftrightarrow \quad a \in \textrm{\emph{Span}} \{a_1,a_2,\ldots,a_k\}.
		\end{displaymath}
	\end{prop}
	\begin{proof}
		This follows directly from Theorem \ref{thm_linear_indep} since $\{a_i\}$
		are linearly independent and $a \wedge A=0 \Leftrightarrow \{a,a_i\}$ are linearly dependent.
	\end{proof}
	
	\noindent
	Hence, to every nonzero $k$-blade $A = a_1 \wedge a_2 \wedge \cdots \wedge a_k$ 
	there corresponds a unique $k$-dimensional subspace
	\begin{equation} \label{blade_subspace}
		\bar{A} := \{ a \in \mathcal{V} : a \wedge A = 0 \} = \textrm{Span} \{a_1,a_2,\ldots,a_k\}.
	\end{equation}
	Conversely, if $\bar{A} \subseteq \mathcal{V}$ is a $k$-dimensional subspace of $\mathcal{V}$, 
	then we can find a nonzero $k$-blade $A$ representing $\bar{A}$ by simply taking a basis
	$\{a_i\}_{i=1}^k$ of $\bar{A}$ and forming
	\begin{equation} \label{subspace_blade}
		A := a_1 \wedge a_2 \wedge \cdots \wedge a_k.
	\end{equation}
	We thus have the geometric interpretation of blades as subspaces with an associated
	orientation (sign) and magnitude.
	Since every element in $\mathcal{G}$ is a linear combination of basis blades,
	we can think of every element as representing a linear combination of subspaces.
	In the case of a nondegenerate algebra these basis subspaces are 
	nondegenerate as well.
	On the other hand, any blade which represents a nondegenerate subspace can also be 
	treated as a basis blade associated to an orthogonal basis.
	This will follow in the discussion below.
	
	\begin{prop} \label{prop_blade_product}
		Every $k$-blade can be written as a geometric product of $k$ vectors.
	\end{prop}
	\begin{proof}
		Take a nonzero $A = a_1 \wedge \cdots \wedge a_k \in \mathcal{B}^*$. 
		Pick an orthogonal basis $\{e_i\}_{i=1}^k$ of the subspace $(\bar{A},q|_{\bar{A}})$.
		Then we can write $a_i = \sum_j \beta_{ij} e_j$ for some $\beta_{ij} \in \mathbb{F}$,
		and $A = \det\ [\beta_{ij}]\ e_1 e_2 \ldots e_k$ by \eqref{blade_expansion}.
	\end{proof}
	
	\noindent
	There are a number of useful consequences of this result.
	
	\begin{cor}
		If $A \in \mathcal{B}$ then $A^2$ is a scalar.
	\end{cor}
	\begin{proof}
		Use the expansion of $A$ above to obtain 
		\begin{equation} \label{blade_squared}
			A^2 = (\det\ [\beta_{ij}])^2\ (-1)^{\frac{1}{2}k(k-1)} q(e_1) q(e_2) \ldots q(e_k) \in \mathbb{F}. \qedhere
		\end{equation}
	\end{proof}
	
	\begin{cor}
		If $A \in \mathcal{B}^\times$ then $A$ has an inverse $A^{-1} = \frac{1}{A^2}A$.
	\end{cor}
	
	\begin{cor}
		If $A \in \mathcal{B}^\times$ then $q$ is nondegenerate on $\bar{A}$
		and there exists an
		orthogonal basis $E$ of $\mathcal{V}$ such that $A \in \mathcal{B}_E$.
	\end{cor}
	\begin{proof}
		The first statement follows directly from \eqref{blade_squared}.
		For the second statement note that, since $q$ is nondegenerate on $\bar{A}$,
		we have $\bar{A} \cap \bar{A}^\perp = 0$. 
		Take an orthogonal basis $\{e_i\}$ of $\bar{A}$.
		For any $v \in \mathcal{V}$ we have that $v - \sum_i \beta_q(v,e_i) q(e_i)^{-1} e_i \in \bar{A}^\perp$.
		Thus, $\mathcal{V} = \bar{A} + \bar{A}^\perp$ and
		we can extend $\{e_i\}$ to an orthogonal basis of $\mathcal{V}$ consisting of 
		one part in $\bar{A}$ and one part in $\bar{A}^\perp$.
		By rescaling this basis we have $A = e_1 \wedge \cdots \wedge e_k$.
	\end{proof}
	
	\begin{rem}
		Note that if we have an orthogonal basis of a subspace of $\mathcal{V}$ where $q$ is degenerate, 
		then it may not be possible to extend this basis to an orthogonal basis for all of $\mathcal{V}$.
		$\mathbb{R}^{1,1}$ for example has two null-spaces, but these are not orthogonal.
		If the space is euclidean or anti-euclidean, though, orthogonal bases can always be extended
		(e.g. using the Gram-Schmidt algorithm).
	\end{rem}
	
	It is useful to be able to work efficiently with general bases of $\mathcal{V}$ 
	and $\mathcal{G}$ which	need not be orthogonal.
	Let $\{e_1,\ldots,e_n\}$ be \emph{any} basis of $\mathcal{V}$.
	Then $\{e_\mathbf{i}\}$ is a basis of $\mathcal{G}(\mathcal{V})$, where we
	use a multi-index notation
	\begin{equation} \label{genbase_index}
		\mathbf{i} = (i_1,i_2,\ldots,i_k), \quad i_1 < i_2 < \ldots < i_k, \quad 0 \leq k \leq n
	\end{equation}
	and
	\begin{equation} \label{genbase_blade}
		e_{()} := 1, \quad e_{(i_1,i_2,\ldots,i_k)} := e_{i_1} \wedge e_{i_2} \wedge \cdots \wedge e_{i_k}.
	\end{equation}
	Sums over $\mathbf{i}$ are understood to be performed over all allowed such indices.
	If $\mathcal{G}$ is nondegenerate then the scalar product $(A,B) \mapsto A * B$
	is also nondegenerate and we can find a so called \emph{reciprocal basis} $\{e^1,\ldots,e^n\}$
	of $\mathcal{V}$ such that
	\begin{equation} \label{reciprocal_def}
		e^i * e_j = \delta^i_j.
	\end{equation}
	The reciprocal basis is easily verified to be given by
	\begin{equation} \label{reciprocal_exp}
		e^i = (-1)^{i-1} (e_1 \wedge \cdots \wedge \check{e}_i \wedge \cdots \wedge e_n) e_{(1,\ldots,n)}^{-1},
	\end{equation}
	where $\check{}$ denotes a deletion.
	Furthermore, we have that $\{e^\mathbf{i}\}$ is a reciprocal basis of $\mathcal{G}$,
	where $e^{(i_1,\ldots,i_k)} := e^{i_k} \wedge \cdots \wedge e^{i_1}$.
	This follows since by Lemma \ref{lem_blade_det} and \eqref{reciprocal_def},
	\begin{equation} \label{reciprocal_g}
		e^\mathbf{i} * e_\mathbf{j} 
			= (e^{i_k} \wedge \cdots \wedge e^{i_1}) * (e_{j_1} \wedge \cdots \wedge e_{j_l})
			= \delta^k_l \det\ [e^{i_p} * e_{j_q}]_{p,q} = \delta^\mathbf{i}_\mathbf{j}.
	\end{equation}
	We now have the coordinate expansions
	\begin{equation} \label{coordinate_expansions}
	\setlength\arraycolsep{2pt}
	\begin{array}{rcll}
		v &=& \sum_i v * e^i e_i = \sum_i v * e_i e^i & \quad \forall\ v \in \mathcal{V}, \\[5pt]
		x &=& \sum_\mathbf{i} x * e^\mathbf{i} e_\mathbf{i} = \sum_\mathbf{i} x * e_\mathbf{i} e^\mathbf{i} & \quad \forall\ x \in \mathcal{G}(\mathcal{V}).
	\end{array}
	\end{equation}
	
	In addition to being useful in coordinate expansions, the general and reciprocal bases
	also provide a geometric understanding of the dual operation because of the following
	
	\begin{thm}
		Assume that $\mathcal{G}$ is nondegenerate.
		If $A = a_1 \wedge \cdots \wedge a_k \in \mathcal{B}^*$ and we extend $\{a_i\}_{i=1}^k$
		to a basis $\{a_i\}_{i=1}^n$ of $\mathcal{V}$ then
		\begin{displaymath}
			A\dual \propto a^{k+1} \wedge a^{k+2} \wedge \cdots \wedge a^n,
		\end{displaymath}
		where $\{a^i\}$ is the reciprocal basis of $\{a_i\}$.
	\end{thm}
	\begin{proof}
		Using an expansion of the inner product into sub-blades (this will not be
		explained in detail here, see \cite{hestenes_sobczyk} or \cite{svensson}) 
		plus orthogonality \eqref{reciprocal_def}, we obtain
		\begin{equation}
		\setlength\arraycolsep{2pt}
		\begin{array}{rcl}
			A\dual
				&=& A \liprod I^{-1} \propto (a_k \wedge \cdots \wedge a_1) \liprod (a^1 \wedge \cdots \wedge a^k \wedge a^{k+1} \wedge \cdots \wedge a^n) \\[5pt]
				&=& (a_k \wedge \cdots \wedge a_1)*(a^1 \wedge \cdots \wedge a^k)\ a^{k+1} \wedge \cdots \wedge a^n \\[5pt]
				&=& a^{k+1} \wedge \cdots \wedge a^n.
		\end{array} \qedhere
		\end{equation}
	\end{proof}
	
	\begin{cor}
		If $A$ and $B$ are blades then $A\dual$, $A \wedge B$, $A \vee B$ and $A \liprod B$ are blades as well.
	\end{cor}
	
	\noindent
	The blade-subspace correspondence then gives us a geometric interpretation of
	these operations.
	
	\begin{prop} \label{prop_subspace_ops}
		If $A,B \in \mathcal{B}^*$ are nonzero blades then $\overline{A\dual} = \bar{A}^\perp$ and
		\begin{displaymath}
		\setlength\arraycolsep{2pt}
		\begin{array}{rcl}
			A \wedge B \neq 0 
				&\Rightarrow& 
				\overline{A \wedge B} = \bar{A} + \bar{B} \ \textrm{and}\ \bar{A} \cap \bar{B} = 0, \\[5pt]
			\bar{A} + \bar{B} = \mathcal{V}
				 &\Rightarrow& 
				\overline{A \vee B} = \bar{A} \cap \bar{B}, \\[5pt]
			A \liprod B \neq 0 
				&\Rightarrow& 
				\overline{A \liprod B} = \bar{A}^\perp \cap \bar{B}, \\[5pt]
			\bar{A} \subseteq \bar{B} 
				&\Rightarrow& 
				A \liprod B = AB, \\[5pt]
			\bar{A} \cap \bar{B}^\perp \neq 0 
				&\Rightarrow& 
				A \liprod B = 0.
		\end{array}
		\end{displaymath}
	\end{prop}
	
	\noindent
	The proofs of the statements in the above corollary and proposition are left as exercises. 
	Some of them can be found in \cite{svensson} and \cite{hestenes_sobczyk}.

\subsection{Linear functions}
	
	Since $\mathcal{G}$ is itself a vector space which embeds $\mathcal{V}$, it is natural to consider the
	properties of linear functions on $\mathcal{G}$. There is a special class of such
	functions, called outermorphisms, which can be said to respect the structure
	of $\mathcal{G}$ in a natural way. We will see that, just as the geometric algebra $\mathcal{G}(\mathcal{V},q)$
	is completely determined by the underlying vector space $(\mathcal{V},q)$, an outermorphism
	is completely determined by its behaviour on $\mathcal{V}$.
	
	\begin{defn} \label{def_outermorphism}
		A linear map $F\!: \mathcal{G} \to \mathcal{G}'$ 
		is called an \emph{outermorphism} or \emph{$\wedge$-morphism} if 
		\begin{displaymath}
		\setlength\arraycolsep{2pt}
		\begin{array}{rl}
			i)   & F(1) = 1, \\[5pt]
			ii)  & F(\mathcal{G}^m) \subseteq \mathcal{G}'^m \quad \forall\ m \geq 0, \quad \textrm{(grade preserving)} \\[5pt]
			iii) & F(x \wedge y) = F(x) \wedge F(y) \quad \forall\ x,y \in \mathcal{G}. \phantom{lagg pa samma tab}
		\end{array}
		\end{displaymath}
		A linear transformation $F\!: \mathcal{G} \to \mathcal{G}$ 
		is called a \emph{dual outermorphism} or \emph{$\vee$-morphism} if 
		\begin{displaymath}
		\setlength\arraycolsep{2pt}
		\begin{array}{rl}
			i)   & F(I) = I, \\[5pt]
			ii)  & F(\mathcal{G}^m) \subseteq \mathcal{G}^m \quad \forall\ m \geq 0, \\[5pt]
			iii) & F(x \vee y) = F(x) \vee F(y) \quad \forall\ x,y \in \mathcal{G}. \phantom{lagg pa samma tab}
		\end{array}
		\end{displaymath}
	\end{defn}
	
	\begin{thm} \label{thm_outermorphism}
		For every linear map $f\!: \mathcal{V} \to \mathcal{W}$ there exists a unique outermorphism
		$f_\wedge\!: \mathcal{G}(\mathcal{V}) \to \mathcal{G}(\mathcal{W})$ such that $f_\wedge(v) = f(v) \ \forall\ v \in \mathcal{V}$.
	\end{thm}
	\begin{proof}
		Take a general basis $\{e_1,\ldots,e_n\}$ of $\mathcal{V}$ and define,
		for $1 \leq i_1 < i_2 < \ldots < i_m \leq n$,
		\begin{equation}
			f_\wedge(e_{i_1} \wedge \cdots \wedge e_{i_m}) := f(e_{i_1}) \wedge \cdots \wedge f(e_{i_m}),
		\end{equation}
		and extend $f_\wedge$ to the whole of $\mathcal{G}(\mathcal{V})$ by linearity.
		We also define $f_\wedge(\alpha) := \alpha$ for $\alpha \in \mathbb{F}$.
		Hence, (\emph{i}) and (\emph{ii}) are satisfied. (\emph{iii}) is easily
		verified by expanding in the induced basis $\{e_\mathbf{i}\}$ of $\mathcal{G}(\mathcal{V})$.
		Unicity is obvious since our definition was necessary.
	\end{proof}
	
	\noindent
	Uniqueness immediately implies the following.
	
	\begin{cor}
		If $f\!: \mathcal{V} \to \mathcal{V}'$ and $g\!: \mathcal{V}' \to \mathcal{V}''$
		are linear then $(g \circ f)_\wedge = g_\wedge \circ f_\wedge$.
	\end{cor}

	\begin{cor}
		If $F\!: \mathcal{G}(\mathcal{V}) \to \mathcal{G}(\mathcal{W})$ is an outermorphism
		then $F = (F|_\mathcal{V})_\wedge$.
	\end{cor}
	
	\noindent
	In the setting of $\cl$ this means that an outermorphism $F\!: \cl(X,R,r) \to \cl(X',R,r')$ 
	is completely determined by its values on $X$.
	
	We have noted that a nondegenerate $\mathcal{G}$ results in a nondegenerate bilinear form $x*y$.
	This gives us a canonical isomorphism $\theta\!: \mathcal{G} \to \mathcal{G}^* = \textrm{Lin}(\mathcal{G},\mathbb{F})$
	between the elements of $\mathcal{G}$ and the linear functionals on	$\mathcal{G}$ as follows.
	For every $x \in \mathcal{G}$ we define a linear functional $\theta(x)$ by $\theta(x)(y) := x*y$.
	Taking a general basis $\{e_\mathbf{i}\}$ of $\mathcal{G}$ and using \eqref{reciprocal_g} we obtain a 
	dual basis $\{\theta(e^\mathbf{i})\}$ such that $\theta(e^\mathbf{i})(e_\mathbf{j}) = \delta_\mathbf{j}^\mathbf{i}$.
	Now that we have a canonical way of moving between $\mathcal{G}$ and its dual space $\mathcal{G}^*$, 
	we can for every linear map $F\!: \mathcal{G} \to \mathcal{G}$ define an \emph{adjoint map}
	$F^*\!: \mathcal{G} \to \mathcal{G}$ by
	\begin{equation} \label{adjoint_def}
		F^*(x) := \theta^{-1}\big( \theta(x) \circ F \big).
	\end{equation}
	Per definition, this has the expected and unique property
	\begin{equation} \label{adjoint_property}
		F^*(x) * y = x * F(y)
	\end{equation}
	for all $x,y \in \mathcal{G}$.
	Note that if we restrict our attention to $\mathcal{V}$ this construction results 
	in the usual adjoint $f^*$ of a linear map $f\!: \mathcal{V} \to \mathcal{V}$.
	
	\begin{thm}[Hestenes' Theorem] \label{thm_hestenes}
		Assume that $\mathcal{G}$ is nondegenerate and let $F\!: \mathcal{G} \to \mathcal{G}$ be an outermorphism. 
		Then the adjoint $F^*$ is also an outermorphism and
		\begin{displaymath}
		\setlength\arraycolsep{2pt}
		\begin{array}{ccc}
			x \liprod F(y) 	&=& F\big( F^*(x) \liprod y \big), \\[5pt]
			F(x) \riprod y 	&=& F\big( x \riprod F^*(y) \big),
		\end{array}
		\end{displaymath}
		for all $x,y \in \mathcal{G}$.
	\end{thm}
	\begin{proof}
		We first prove that $F^*$ is an outermorphism. 
		The fact that $F^*$ is grade preserving follows from \eqref{adjoint_property} and
		the grade preserving property of $F$. Now take basis blades
		$x = x_1 \wedge \cdots \wedge x_m$ and $y = y_m \wedge \cdots \wedge y_1$ 
		with $x_i,y_j \in \mathcal{V}$. Then
		\begin{displaymath}
		\setlength\arraycolsep{2pt}
		\begin{array}{rcl}
			F^*(x_1 \wedge \cdots \wedge x_m) * y
				&=& (x_1 \wedge \cdots \wedge x_m) * F(y_m \wedge \cdots \wedge y_1) \\[5pt]
				&=& (x_1 \wedge \cdots \wedge x_m) * \big(F(y_m) \wedge \cdots \wedge F(y_1)\big) \\[5pt]
				&=& \det\ [x_i * F(y_j)]_{i,j} = \det\ [F^*(x_i) * y_j]_{i,j} \\[5pt]
				&=& \big(F^*(x_1) \wedge \cdots \wedge F^*(x_m)\big) * (y_m \wedge \cdots \wedge y_1) \\[5pt]
				&=& \big(F^*(x_1) \wedge \cdots \wedge F^*(x_m)\big) * y,
		\end{array}
		\end{displaymath}
		where we have used Lemma \ref{lem_blade_det}.
		By linearity and nondegeneracy it follows that $F^*$ is an outermorphism.
		The first identity stated in the therorem now follows quite easily from 
		Proposition \ref{prop_trix}. For any $z \in \mathcal{G}$ we have
		\begin{displaymath}
		\setlength\arraycolsep{2pt}
		\begin{array}{rcl}
			z * \big(x \liprod F(y)\big)
				&=& (z \wedge x) * F(y) = F^*(z \wedge x) * y \\[5pt]
				&=& \big(F^*(z) \wedge F^*(x)\big) * y = F^*(z) * \big(F^*(x) \liprod y\big) \\[5pt]
				&=& z * F\big(F^*(x) \liprod y\big).
		\end{array}
		\end{displaymath}
		The nondegeneracy of the scalar product then gives the first identity.
		The second identity is proven similarly, using that $(x \riprod y) * z = x * (y \wedge z)$.
	\end{proof}
	
	\noindent
	From uniqueness of outermorphisms we also obtain the following
	
	\begin{cor}
		If $f\!: \mathcal{V} \to \mathcal{V}$ is a linear transformation then
		$(f^*)_\wedge = (f_\wedge)^*$.
	\end{cor}

	\noindent
	This means that we can simply write $f^*_\wedge$ for the adjoint outermorphism of $f$.
	
	A powerful concept in geometric algebra (or exterior algebra) is the generalization of eigenvectors
	to so called \emph{eigenblades}. For a function $f\!: \mathcal{V} \to \mathcal{V}$, a $k$-eigenblade
	with eigenvalue $\lambda \in \mathbb{F}$ is a blade $A \in \mathcal{B}_k$ such that
	\begin{equation} \label{eigenblade_def}
		f_\wedge(A) = \lambda A.
	\end{equation}
	Just as eigenvectors can be said to represent 
	invariant 1-dimensional subspaces of a function, a $k$-blade with nonzero
	eigenvalue represents an invariant $k$-dimensional subspace.
	One important example of an eigenblade is the pseudoscalar $I$, which represents the whole 
	invariant vector space $\mathcal{V}$.
	Since $f_\wedge$ is
	grade preserving, we must have $f_\wedge(I) = \lambda I$ for some $\lambda \in \mathbb{F}$
	which we call the \emph{determinant} of $f$, i.e.
	\begin{equation} \label{determinant_def}
		f_\wedge(I) = \det f\ I.
	\end{equation}
	Expanding $\det f = f(I) * I^{-1}$ in a basis using Lemma \ref{lem_blade_det}, one finds
	that this agrees with the usual definition of the determinant of a linear function.
	
	In the following we assume that $\mathcal{G}$ is nondegenerate, so that $I^2 \neq 0$.
	
	\begin{defn} \label{def_dual_map}
		For linear $F\!: \mathcal{G} \to \mathcal{G}$ we define the \emph{dual map} 
		$F\dual\!: \mathcal{G} \to \mathcal{G}$ by
		$F\dual(x) := F(xI)I^{-1}$, so that the following diagram commutes:
		\begin{displaymath}
		\setlength\arraycolsep{2pt}
		\begin{array}{rcccl}
						& \mathcal{G} 	& \xrightarrow{F} 		& \mathcal{G} \\
			(\ )\dual 	& \downarrow	&						& \downarrow	& (\ )\dual \\
						& \mathcal{G}	& \xrightarrow{F\dual}	& \mathcal{G}
		\end{array}
		\end{displaymath}
	\end{defn}
	
	\begin{prop}
		We have the following properties of the dual map:
		\begin{displaymath}
		\setlength\arraycolsep{2pt}
		\begin{array}{rrcl}
			i) 	& F^\mathbf{cc} &=& F, \\[5pt]
			ii) & (F \circ G)\dual &=& F\dual \circ G\dual, \\[5pt]
			iii)& id\dual &=& id, \\[5pt]
			iv) & F(\mathcal{G}^s) \subseteq \mathcal{G}^t &\Rightarrow& F\dual(\mathcal{G}^{\dim \mathcal{V} -s}) \subseteq \mathcal{G}^{\dim \mathcal{V} -t}, \\[5pt]
			v) 	& F \ \textrm{$\wedge$-morphism} &\Rightarrow& F\dual \ \textrm{$\vee$-morphism},
		\end{array}
		\end{displaymath}
		for all linear $F,G\!: \mathcal{G} \to \mathcal{G}$.
	\end{prop}
	
	\noindent
	The proofs are simple and left as exercises to the reader.
	As a special case of Theorem \ref{thm_hestenes} we obtain, with $y=I$
	and a linear map $f\!: \mathcal{V} \to \mathcal{V}$,
	\begin{equation}
		\det f\ xI = f_\wedge \big( f_\wedge^*(x)I \big),
	\end{equation}
	so that
	\begin{equation}
		\det f\ id = f_\wedge\dual \circ f_\wedge^* = f_\wedge \circ f_\wedge^{*\mathbf{c}}.
	\end{equation}
	If $\det f \neq 0$ we then have a simple expression for the inverse;
	\begin{equation}
		f_\wedge^{-1} = (\det f)^{-1} f_\wedge^{*\mathbf{c}},
	\end{equation}
	which is essentially the dual of the adjoint. $f^{-1}$ is obtained 
	by restricting to $\mathcal{V}$.
	An orthogonal transformation $f$ has $f^{-1} = f^*$ and $\det f = 1$, 
	so in that case $f_\wedge = f_\wedge\dual$.

\subsection{Infinite-dimensional Clifford algebra}
	
	This far we have only defined the Clifford algebra 
	$\cl(X,R,r)$ of a \emph{finite} set $X$, resulting
	in a \emph{finite-dimen\-sional} algebra $\mathcal{G}(\mathcal{V})$ whenever $R$ is a field.
	In order for this combinatorial construction to qualify as a complete generalization of $\mathcal{G}$,
	we would at least like to be able to define the corresponding
	Clifford algebra of an infinite-dimensional vector space, something which
	was possible for $\mathcal{G}$ in Definition \ref{def_g}.
	
	The treatment of $\cl$ in the previous subsections has been 
	deliberately put in a form which eases the generalization to an infinite $X$.
	Reconsidering Definition \ref{def_cl}, we now have two possibilites; 
	either we consider the set $\mathscr{P}(X)$ of all subsets of $X$,
	or the set $\mathscr{F}(X)$ of all \emph{finite} subsets.
	We therefore define, for an arbitrary set $X$, ring $R$, and signature $r \!: X \to R$,
	\begin{equation}
		\cl(X,R,r) := \bigoplus_{\mathscr{P}(X)} R
		\qquad \textrm{and} \qquad
		\cl_\mathscr{F}(X,R,r) := \bigoplus_{\mathscr{F}(X)} R.
	\end{equation}
	Elements in $\cl$ ($\cl_\mathscr{F}$) are finite linear combinations
	of (finite) subsets of $X$.
	
	Our problem now is to define a Clifford product for $\cl$ and $\cl_\mathscr{F}$.
	This can be achieved just as in the finite case if only we can find 
	a map $\tau \!: \mathscr{P}(X) \times \mathscr{P}(X) \to R$
	satisfying the conditions in Lemma \ref{lem_finite_tau}.
	This is certainly not a trivial task. Starting with the case $\cl_\mathscr{F}$ it
	is sufficient to construct such a map on $\mathscr{F}(X)$.
	
	We call a map $\tau \!: \mathscr{F}(X) \times \mathscr{F}(X) \to R$ \emph{grassmannian on} $X$ if it
	satisfies (\emph{i})-(\emph{v}) in Lemma \ref{lem_finite_tau}, with 
	$\mathscr{P}(X)$ replaced by $\mathscr{F}(X)$.
	
	\begin{thm}
		For any $X,R,r$ there exists a grassmannian map on $\mathscr{F}(X)$.
	\end{thm}
	\begin{proof}
		We know that there exists such a map for any finite $X$.
		Let $Y \subseteq X$ and assume $\tau' \!: \mathscr{F}(Y) \times \mathscr{F}(Y) \to R$
		is grassmannian on $Y$.
		If there exists $z \in X \smallsetminus Y$ we can, by proceeding as in the
		proof of Lemma \ref{lem_finite_tau}, extend $\tau'$ to a grassmannian map 
		$\tau \!: \mathscr{F}(Y\cup\{z\}) \times \mathscr{F}(Y\cup\{z\}) \to R$
		on $Y\cup\{z\}$ such that $\tau|_{\mathscr{F}(Y) \times \mathscr{F}(Y)} = \tau'$.
		
		We will now use \emph{transfinite induction}, or the \emph{Hausdorff maximality theorem}\footnote{This 
		theorem should actually be regarded as an axiom of set theory since it is equivalent to the Axiom of Choice.},
		to prove that $\tau$ can be extended to all of $\mathscr{F}(X) \subseteq \mathscr{P}(X)$.
		Note that if $\tau$ is grassmannian on $Y \subseteq X$ then $\tau$
		is also a relation $\tau \subseteq \mathscr{P}(X) \times \mathscr{P}(X) \times R$.
		Let
		\begin{equation}
			\mathcal{H} := \Big\{ (Y,\tau) \in \mathscr{P}(X) \times \mathscr{P}\big( \mathscr{P}(X) \times \mathscr{P}(X) \times R \big) 
				: \textrm{$\tau$ is grassmannian on $Y$} \Big\}.
		\end{equation}
		Then $\mathcal{H}$ is partially ordered by
		\begin{equation}
			(Y,\tau) \leq (Y',\tau') \qquad \textrm{iff} \qquad 
				Y \subseteq Y' \quad \textrm{and} \quad \tau'|_{\mathscr{F}(Y) \times \mathscr{F}(Y)} = \tau.
		\end{equation}
		By the Hausdorff maximality theorem, there exists a maximal totally 
		ordered ``chain" $\mathcal{K} \subseteq \mathcal{H}$.
		Put $Y^* := \bigcup_{(Y,\tau) \in \mathcal{K}} Y$. We want to define a 
		grassmannian map $\tau^*$ on $Y^*$, for if we succeed in that, we find 
		$(Y^*,\tau^*) \in \mathcal{H} \cap \mathcal{K}$ and can conclude that $Y^* = X$
		by maximality and the former result.
		
		Take finite subsets $A$ and $B$ of $Y^*$. Each of the finite elements
		in $A \cup B$ lies in some $Y$ such that $(Y,\tau) \in \mathcal{K}$.
		Therefore, by the total ordering of $\mathcal{K}$, there exists one
		such $Y$ containing $A \cup B$.
		Put $\tau^* (A,B) := \tau(A,B)$, where $(Y,\tau)$ is this chosen element in $\mathcal{K}$.
		$\tau^*$ is well-defined since if $A \cup B \subseteq Y$ and $A \cup B \subseteq Y'$
		where $(Y,\tau), (Y',\tau') \in \mathcal{K}$ then $Y \subseteq Y'$ or $Y' \subseteq Y$
		and $\tau,\tau'$ agree on $(A,B)$.
		It is easy to verify that this $\tau^*$ is grassmannian on $Y^*$,
		since for each $A,B,C \in \mathscr{F}(X)$ there exists $(Y,\tau) \in \mathcal{K}$ 
		such that $A \cup B \cup C \subseteq Y$.
	\end{proof}
	
	We have shown that there exists a map $\tau \!: \mathscr{F}(X) \times \mathscr{F}(X) \to R$
	with the properties in Lemma \ref{lem_finite_tau}.
	We can then define the Clifford product on $\cl_\mathscr{F}(X)$ as usual by
	$AB := \tau(A,B) A \symdiff B$ for $A,B \in \mathscr{F}(X)$ and linear extension.
	Since only finite subsets are included, most of the previous
	constructions for finite-dimensional $\cl$ carry over to $\cl_\mathscr{F}$.
	For example, the decomposition into graded subspaces remains but now goes up towards infinity,
	\begin{equation} \label{infinite_grade_decomposition}
		\cl_\mathscr{F} = \bigoplus_{k=0}^{\infty} \cl_\mathscr{F}^k.
	\end{equation}
	Furthermore, Proposition \ref{prop_reverse} still holds, so the reverse and
	all other involutions behave as expected.
	
	The following theorem shows that it is possible to extend $\tau$ all the
	way to $\mathscr{P}(X)$ even in the infinite case. 
	We therefore have a Clifford product also on $\cl(X)$.
	
	\begin{thm} \label{thm_infinite_tau}
		For any set $X$ there exists a map $|\cdot|_2 \!: \mathscr{P}\big(\mathscr{P}(X)\big) \to \mathbb{Z}_2$
		such that
		\begin{displaymath}
		\setlength\arraycolsep{2pt}
		\begin{array}{rl}
			i)   & |\mathcal{A}|_2 \equiv |\mathcal{A}| \pmod{2} 
				\qquad \textrm{for finite $\mathcal{A} \subseteq \mathscr{P}(X)$}, \\[5pt]
			ii)  & |\mathcal{A} \cup \mathcal{B}|_2 = |\mathcal{A}|_2 + |\mathcal{B}|_2 \pmod{2}
				\qquad \textrm{if} \quad \mathcal{A} \cap \mathcal{B} = \varnothing.
		\end{array}
		\end{displaymath}
		Furthermore, for any commutative ring $R$ and signature $r\!: X \to R$
		such that $r(X)$ is contained in a finite and multiplicatively closed subset of $R$,
		there exists a map  $\tau \!: \mathscr{P}(X) \times \mathscr{P}(X) \to R$ such that
		properties \emph{(\emph{i})-(\emph{v})} in Lemma \ref{lem_finite_tau} hold, plus
		\begin{displaymath}
			vi) \ \tau(A,B) = (-1)^{\left|\binom{A}{2}\right|_2 + \left|\binom{B}{2}\right|_2 + \left|\binom{A \bigtriangleup B}{2}\right|_2}\ \tau(B,A)
				\quad \forall\ A,B \in \mathscr{P}(X).
		\end{displaymath}
	\end{thm}
	
	\noindent
	Here, $\binom{A}{n}$ denotes the set of all subsets of $A$ with $n$ elements.
	Note that for a finite set $A$, $\big|\binom{A}{n}\big| = \binom{|A|}{n}$ so that for example
	$\big|\binom{A}{1}\big| = |A|$ (in general, $\textrm{card}\ \binom{A}{1} = \textrm{card}\ A$)
	and $\big|\binom{A}{2}\big| = \frac{1}{2}|A|(|A|-1)$.
	This enables us to extend the basic involutions $\star$, $\dagger$ and $\cliffconj$
	to infinite sets as
	\begin{displaymath}
	\setlength\arraycolsep{2pt}
	\begin{array}{rcl}
		A^\star   &:=& (-1)^{\left|\binom{A}{1}\right|_2}\ A, \\[5pt]
		A^\dagger &:=& (-1)^{\left|\binom{A}{2}\right|_2}\ A,
	\end{array}
	\end{displaymath}
	and because $\big|\binom{A \bigtriangleup B}{1}\big|_2 = \big|\binom{A}{1}\big|_2 + \big|\binom{B}{1}\big|_2 \pmod{2}$
	still holds, we find that they satisfy the fundamental properties
	\eqref{grade_inv_property}-\eqref{cliffconj_property} for all elements of $\cl(X)$.
	The extra requirement (\emph{vi}) on $\tau$ was necessary here since we 
	cannot use Proposition \ref{prop_reverse} for infinite sets.
	Moreover, we can no longer write the decomposition \eqref{infinite_grade_decomposition} 
	since it goes beyond finite grades.
	We do have even and odd subspaces, though, defined by
	\begin{equation} \label{infinite_evenodd_subspaces}
		\cl^{\pm} := \{ x \in \cl : x^\star = \pm x \}.
	\end{equation}
	$\cl^+$ and $\cl_\mathscr{F}$ (with this $\tau$) are both subalgebras of $\cl$.
	
	It should be emphasized that $\tau$ needs not be zero on intersecting infinite sets
	(a rather trivial solution), but if e.g. $r \!: X \to \{\pm 1\}$ we can
	also demand that $\tau \!: \mathscr{P}(X) \times \mathscr{P}(X) \to \{\pm 1\}$.
	Theorem \ref{thm_infinite_tau} can be proved using \emph{nonstandard analysis} /
	\emph{internal set theory} and we will not consider this here.

	Let us now see how $\cl_\mathscr{F}$ and $\cl$ can be applied to the 
	setting of an infinite-dimensional vector space $\mathcal{V}$ 
	over a field $\mathbb{F}$ and with a quadratic form $q$.
	By the Hausdorff maximality theorem one can actually find a (necessarily infinite) orthogonal
	basis $E$ for this space in the sense that any vector in $\mathcal{V}$
	can be written as a finite linear combination of elements in $E$ and
	that $\beta_q(e,e') = 0$ for any pair of disjoint elements $e,e' \in E$.
	We then have
	\begin{equation}
		\mathcal{G}(\mathcal{V},q) \cong \cl_\mathscr{F}(E,\mathbb{F},q|_E),
	\end{equation}
	which is proved just like in the finite-dimensional case, using
	Proposition \ref{prop_universality}. The only difference is that one
	needs to check that the homomorphism $F\!: \mathcal{G} \to \cl_\mathscr{F}$
	is also injective.
	
	The $k$-blades of $\cl_\mathscr{F}$ represent $k$-dimensional subspaces
	of $\mathcal{V}$ even in the infinite case. Due to the intuitive
	and powerful handling of finite-dimensional geometry which was possible
	in a finite-dimensional $\mathcal{G}$, it would be extremely satisfying
	to be able to generalize the blade concept to e.g. closed subspaces
	of an infinite-dimensional Hilbert space. One could hope that the
	infinite basis subsets in $\cl(E)$ provide this generalization.
	Unfortunately, this is not so easy since $\cl(E)$ depends heavily
	on the choice of basis $E$.
	Let us sketch an intuitive picture of why this is so.
	
	With a countable basis $E = \{e_i\}_{i=1}^\infty$, an
	infinite basis blade in $\cl$ could be thought of as an infinite product
	$A = e_{i_1} e_{i_2} e_{i_3} \ldots = \prod_{k=1}^\infty e_{i_k}$.
	A change of basis to $E'$ would turn each $e \in E$ into a finite linear combination
	of elements in $E'$, e.g. $e_j = \sum_k \beta_{jk} e_k'$. However,
	this would require $A$ to be an \emph{infinite sum} of basis blades in $E'$, which is not allowed.
	Note that this is no problem in $\cl_\mathscr{F}$ since a
	basis blade $A = \prod_{k=1}^N e_{i_k}$ is a finite product
	and the change of basis therefore results in a finite sum.
	It may be possible to treat infinite sums in $\cl$ by taking the topology 
	of $\mathcal{V}$ into account, but at present this issue is not clear.
	
	
	Finally, we consider a nice application of the infinite-dimensional
	Clifford algebra $\cl_\mathscr{F}$. For a vector space $\mathcal{V}$,
	define the \emph{simplicial complex algebra}
	\begin{equation}
		\mathscr{C}(\mathcal{V}) := \cl_\mathscr{F}(\mathcal{V},R,1),
	\end{equation}
	where we forget about the vector space structure of $\mathcal{V}$ and treat
	individual points $\dot{v} \in \mathcal{V}$ as orthogonal basis 1-vectors 
	in $\cl_\mathscr{F}^1$ with $\dot{v}^2 = 1$.
	The dot indicates that we think of $v$ as a point rather than a vector.
	A basis $(k+1)$-blade in $\mathscr{C}(\mathcal{V})$ consists of a product
	$\dot{v}_0 \dot{v}_1 \ldots \dot{v}_k$ of individual points and represents
	a (possibly degenerate) oriented $k$-simplex in $\mathcal{V}$. 
	This simplex is given by the convex hull
	\begin{equation}
		\textrm{Conv}\{v_0,,v_1,\ldots,v_k\} := \left\{ \sum_{i=0}^k \alpha_i v_i \in \mathcal{V} : \alpha_i \geq 0, \sum_{i=0}^k \alpha_i = 1 \right\}.
	\end{equation}
	Hence, an arbitrary element in $\mathscr{C}(\mathcal{V})$ is a linear
	combination of simplices and can therefore represent a simplicial complex in $\mathcal{V}$.
	The restriction of $\mathscr{C}(\mathcal{V})$ to the $k$-simplices of
	a simplicial complex $K$ is usually called the \emph{$k$-chain group} $C_k(K)$.
	Here the generality of the ring $R$ comes in handy because one often
	works with $R = \mathbb{Z}$ in this context.
	
	The Clifford algebra structure of $\mathscr{C}(\mathcal{V})$ handles the
	orientation of the simplices, so that e.g. the line from the point $\dot{v}_0$
	to $\dot{v}_1$ is $\dot{v}_0\dot{v}_1 = - \dot{v}_1\dot{v}_0$.
	Furthermore, it allows us to define the \emph{boundary operator}
	\begin{displaymath}
	\setlength\arraycolsep{2pt}
	\begin{array}{l}
		\partial \!: \mathscr{C}(\mathcal{V}) \to \mathscr{C}(\mathcal{V}), \\[5pt]
		\displaystyle \quad \partial(x) := \sum_{\dot{v} \in \mathcal{V}} \dot{v} \liprod x.
	\end{array}
	\end{displaymath}
	Note that this is well-defined since only a finite number of points $\dot{v}$ 
	can be present in any fixed $x$. For a $k$-simplex, we have
	\begin{equation}
		\partial(\dot{v}_0 \dot{v}_1 \ldots \dot{v}_k) 
			= \sum_{i=0}^k \dot{v}_i \liprod (\dot{v}_0 \dot{v}_1 \ldots \dot{v}_k)
			= \sum_{i=0}^k (-1)^i \dot{v}_0 \dot{v}_1 \ldots \check{\dot{v}}_i \ldots \dot{v}_k.
	\end{equation}
	This shows that $\partial$ really is the traditional boundary operator on simplices.
	Proposition \ref{prop_trix} now makes the proof of $\partial^2 = 0$ a triviality,
	\begin{equation}
		\partial^2(x) 
			= \sum_{\dot{u} \in \mathcal{V}} \dot{u} \liprod \left( \sum_{\dot{v} \in \mathcal{V}} \dot{v} \liprod x \right)
			= \sum_{\dot{u}, \dot{v} \in \mathcal{V}} \dot{u} \liprod (\dot{v} \liprod x)
			= \sum_{\dot{u}, \dot{v} \in \mathcal{V}} (\dot{u} \wedge \dot{v}) \liprod x = 0.
	\end{equation}
	
	We can also assign a \emph{geometric measure} $\sigma$ to simplices,
	by mapping a $k$-simplex to a corresponding $k$-blade in $\mathcal{G}(\mathcal{V})$
	representing the directed volume of the simplex.
	Define $\sigma \!: \mathscr{C}(\mathcal{V}) \to \mathcal{G}(\mathcal{V})$ by
	\begin{displaymath}
	\setlength\arraycolsep{2pt}
	\begin{array}{rcl}
		\sigma(1) &:=& 0, \\[5pt]
		\sigma(\dot{v}) &:=& 1, \\[5pt]
		\sigma(\dot{v}_0 \dot{v}_1 \ldots \dot{v}_k) &:=& \frac{1}{k!} (v_1-v_0) \wedge (v_2-v_0) \wedge \cdots \wedge (v_k-v_0),
	\end{array}
	\end{displaymath}
	and extending linearly. One can verify that this is well-defined and that the
	geometric measure of a boundary is zero, i.e. $\sigma \circ \partial = 0$.
	One can take this construction even further and arrive at
	``discrete" equivalents of differentials, integrals and Stokes' theorem.
	See \cite{svensson} or \cite{naeve_svensson} for more on this.
	
	This completes our excursion to infinite-dimensional Clifford algebras.
	In the following sections we will always assume that $X$ is finite
	and $\mathcal{V}$ finite-dimensional.

	\newpage

\section{Isomorphisms} \label{sec_isomorphisms}

	In this section we establish an extensive set of relations between real and complex
	geometric algebras of varying signature. This eventually leads to an identification
	of these algebras as matrix algebras over $\mathbb{R}$, $\mathbb{C}$, or the
	quaternions $\mathbb{H}$. The complete listing of such identifications is usually 
	called the \emph{classification} of geometric algebras.
	
	We have seen that the even subspace $\mathcal{G}^{+}$ of $\mathcal{G}$ constitutes
	a subalgebra. The following proposition tells us that this subalgebra actually
	is the geometric algebra of a space of one dimension lower.
	
	\begin{prop} \label{prop_iso_even}
		We have the algebra isomorphisms
		\begin{displaymath}
		\setlength\arraycolsep{2pt}
		\begin{array}{c}
			\mathcal{G}^{+}(\mathbb{R}^{s,t}) \cong \mathcal{G}(\mathbb{R}^{s,t-1}), \\[5pt]
			\mathcal{G}^{+}(\mathbb{R}^{s,t}) \cong \mathcal{G}(\mathbb{R}^{t,s-1}),
		\end{array}
		\end{displaymath}
		for all $s,t$ for which the expressions make sense.
	\end{prop}
	\begin{proof}
		Take an orthonormal basis $\{e_1,\ldots,e_s,\epsilon_1,\ldots,\epsilon_t\}$ of $\mathbb{R}^{s,t}$ 
		such that $e_i^2 = 1$, $\epsilon_i^2 = -1$, and
		a corresponding basis $\{\underline{e}_1,\ldots,\underline{e}_s,\underline{\epsilon}_1,\ldots,\underline{\epsilon}_{t-1}\}$ of $\mathbb{R}^{s,t-1}$.
		Define $f\!: \mathbb{R}^{s,t-1} \to \mathcal{G}^{+}(\mathbb{R}^{s,t})$ by mapping
		\begin{displaymath}
		\setlength\arraycolsep{2pt}
		\begin{array}{ccl}
			\underline{e}_i &\mapsto& e_i \epsilon_t, \quad i=1,\ldots,s, \\[3pt]
			\underline{\epsilon}_i &\mapsto& \epsilon_i \epsilon_t, \quad i=1,\ldots,t-1,
		\end{array}
		\end{displaymath}
		and extending linearly. We then have
		\begin{displaymath}
		\setlength\arraycolsep{2pt}
		\begin{array}{rcl}
			f(\underline{e}_i) f(\underline{e}_j) &=& -f(\underline{e}_j) f(\underline{e}_i), \\[3pt]
			f(\underline{\epsilon}_i) f(\underline{\epsilon}_j) &=& -f(\underline{\epsilon}_j) f(\underline{\epsilon}_i)
		\end{array}
		\end{displaymath}
		for $i \neq j$, and
		\begin{displaymath}
		\setlength\arraycolsep{2pt}
		\begin{array}{c}
			f(\underline{e}_i) f(\underline{\epsilon}_j) = -f(\underline{\epsilon}_j) f(\underline{e}_i), \\[3pt]
			f(\underline{e}_i)^2 = 1, \quad f(\underline{\epsilon}_i)^2 = -1 
		\end{array}
		\end{displaymath}
		for all reasonable $i,j$.
		By Proposition \ref{prop_universality} (universality) we can extend $f$ to a homomorphism $F\!: \mathcal{G}(\mathbb{R}^{s,t-1}) \to \mathcal{G}^{+}(\mathbb{R}^{s,t})$.
		Since $\dim \mathcal{G}(\mathbb{R}^{s,t-1}) = 2^{s+t-1} = 2^{s+t}/2 = \dim \mathcal{G}^{+}(\mathbb{R}^{s,t})$ and
		$F$ is easily seen to be surjective, we have that $F$ is an isomorphism.
		
		For the second statement, we take a corresponding basis 
		$\{\underline{e}_1,\ldots,\underline{e}_t,\underline{\epsilon}_1,\ldots \\ \ldots,\underline{\epsilon}_{s-1}\}$ of $\mathbb{R}^{t,s-1}$
		and define $f\!: \mathbb{R}^{t,s-1} \to \mathcal{G}^{+}(\mathbb{R}^{s,t})$ by
		\begin{displaymath}
		\setlength\arraycolsep{2pt}
		\begin{array}{ccl}
			\underline{e}_i &\mapsto& \epsilon_i e_s, \quad i=1,\ldots,t, \\[3pt]
			\underline{\epsilon}_i &\mapsto& e_i e_s, \quad i=1,\ldots,s-1.
		\end{array}
		\end{displaymath}
		Proceeding as above, we obtain the isomorphism.
	\end{proof}
	
	\begin{cor}
		It follows immediately that
		\begin{displaymath}
		\setlength\arraycolsep{2pt}
		\begin{array}{c}
			\mathcal{G}(\mathbb{R}^{s,t}) \cong \mathcal{G}(\mathbb{R}^{t+1,s-1}), \\[5pt]
			\mathcal{G}^{+}(\mathbb{R}^{s,t}) \cong \mathcal{G}^{+}(\mathbb{R}^{t,s}).
		\end{array}
		\end{displaymath}
	\end{cor}
	
	\noindent
	In the above and further on we use the notation 
	$\mathcal{G}(\mathbb{F}^{0,0}) := \cl(\varnothing,\mathbb{F},\varnothing) = \mathbb{F}$
	for completeness.
	
	The property of geometric algebras that leads us to their eventual classification
	as matrix algebras is that they can be split up into tensor products of
	geometric algebras of lower dimension.
	
	\begin{prop} \label{prop_iso_tensor}
		We have the algebra isomorphisms
		\begin{displaymath}
		\setlength\arraycolsep{2pt}
		\begin{array}{c}
			\mathcal{G}(\mathbb{R}^{n+2,0}) \cong \mathcal{G}(\mathbb{R}^{0,n}) \otimes \mathcal{G}(\mathbb{R}^{2,0}), \\[5pt]
			\mathcal{G}(\mathbb{R}^{0,n+2}) \cong \mathcal{G}(\mathbb{R}^{n,0}) \otimes \mathcal{G}(\mathbb{R}^{0,2}), \\[5pt]
			\mathcal{G}(\mathbb{R}^{s+1,t+1}) \cong \mathcal{G}(\mathbb{R}^{s,t}) \otimes \mathcal{G}(\mathbb{R}^{1,1}),
		\end{array}
		\end{displaymath}
		for all $n$, $s$ and $t$ for which the expressions make sense.
	\end{prop}
	\begin{proof}
		For the first expression, take orthonormal bases 
		$\{e_i\}$ of $\mathbb{R}^{n+2}$, $\{\underline{\epsilon}_i\}$ of $\mathbb{R}^{0,n}$ 
		and $\{\overline{e}_i\}$ of $\mathbb{R}^2$.
		Define a mapping $f\!: \mathbb{R}^{n+2} \to \mathcal{G}(\mathbb{R}^{0,n}) \otimes \mathcal{G}(\mathbb{R}^2)$ by
		\begin{displaymath}
		\setlength\arraycolsep{1pt}
		\begin{array}{ccrcll}
			e_j &\ \mapsto\ & \underline{\epsilon}_j &\otimes& \overline{e}_1 \overline{e}_2, &\quad j=1,\ldots,n, \\[3pt]
			e_j &\ \mapsto\ & 1 &\otimes& \overline{e}_{j-n}, &\quad j=n+1,n+2,
		\end{array}
		\end{displaymath}
		and extend to an algebra homomorphism $F$ using the universal property.
		Since $F$ maps onto a set of generators for $\mathcal{G}(\mathbb{R}^{0,n}) \otimes \mathcal{G}(\mathbb{R}^2)$
		it is clearly surjective. Furthermore, $\dim \mathcal{G}(\mathbb{R}^{n+2}) = 2^{n+2} = \dim \mathcal{G}(\mathbb{R}^{0,n}) \otimes \mathcal{G}(\mathbb{R}^2)$,
		so $F$ is an isomorphism.
		
		The second expression is proved similarly. For the third expression, take orthonormal bases 
		$\{e_1,\ldots,e_{s+1},\epsilon_1,\ldots,\epsilon_{t+1}\}$ of $\mathbb{R}^{s+1,t+1}$,
		$\{\underline{e}_1,\ldots,\underline{e}_s,\underline{\epsilon}_1,\ldots,\underline{\epsilon}_t\}$ of $\mathbb{R}^{s,t}$
		and $\{\overline{e},\overline{\epsilon}\}$ of $\mathbb{R}^{1,1}$, where
		$e_i^2 = 1, \epsilon_i^2 = -1$ etc. Define $f\!: \mathbb{R}^{s+1,t+1} \to \mathcal{G}(\mathbb{R}^{s,t}) \otimes \mathcal{G}(\mathbb{R}^{1,1})$ by
		\begin{displaymath}
		\setlength\arraycolsep{1pt}
		\begin{array}{ccrcll}
			e_j &\ \mapsto\ & \underline{e}_j &\otimes& \overline{e} \overline{\epsilon}, &\quad j=1,\ldots,s, \\[3pt]
			\epsilon_j &\ \mapsto\ & \underline{\epsilon}_j &\otimes& \overline{e} \overline{\epsilon}, &\quad j=1,\ldots,t, \\[3pt]
			e_{s+1} &\ \mapsto\ & 1 &\otimes& \overline{e}, \\[3pt]
			\epsilon_{t+1} &\ \mapsto\ & 1 &\otimes& \overline{\epsilon}.
		\end{array}
		\end{displaymath}
		Proceeding as above, we can extend $f$ to an algebra isomorphism.
	\end{proof}
	
	We can also relate certain real geometric algebras to complex equivalents.
	
	\begin{prop} \label{prop_iso_real_complex}
		If $s+t$ is odd and $I^2 = -1$ then
		\begin{displaymath}
			\mathcal{G}(\mathbb{R}^{s,t}) 
				\cong \mathcal{G}^{+}(\mathbb{R}^{s,t}) \otimes \mathbb{C} 
				\cong \mathcal{G}(\mathbb{C}^{s+t-1}).
		\end{displaymath}
	\end{prop}
	\begin{proof}
		Since $s+t$ is odd, the pseudoscalar $I$ commutes with all other elements.
		This, together with the property $I^2 = -1$, makes it a good candidate for a scalar imaginary.
		Define $F\!: \mathcal{G}^{+}(\mathbb{R}^{s,t}) \otimes \mathbb{C} \to \mathcal{G}(\mathbb{R}^{s,t})$ 
		by linear extension of
		\begin{displaymath}
		\setlength\arraycolsep{1pt}
		\begin{array}{ccccll}
			E &\otimes& 1 &\ \mapsto\ & E & \quad \in \mathcal{G}^+, \\[3pt]
			E &\otimes& i &\ \mapsto\ & EI & \quad \in \mathcal{G}^-,
		\end{array}
		\end{displaymath}
		for even basis blades $E$. $F$ is easily seen to be an injective algebra homomorphism. 
		Using that the dimensions of these algebras are equal, we have an isomorphism.
		
		For the second isomorphism, note that Proposition \ref{prop_iso_even} gives us
		$\mathcal{G}^{+}(\mathbb{R}^{s,t}) \otimes \mathbb{C} \cong \mathcal{G}(\mathbb{R}^{s,t-1}) \otimes \mathbb{C}$.
		Finally, the order of complexification is unimportant since all nondegenerate
		complex quadratic forms are equivalent.
	\end{proof}
	
	\begin{cor}
		It follows immediately that, for these conditions,
		\begin{displaymath}
			\mathcal{G}(\mathbb{R}^{s,t}) \cong \mathcal{G}(\mathbb{R}^{p,q-1}) \otimes \mathbb{C}
		\end{displaymath}
		for any $p \geq 0$, $q \geq 1$ such that $p+q=s+t$.
	\end{cor}
	
	One important consequence of the tensor algebra isomorphisms in Proposition \ref{prop_iso_tensor}
	is that geometric algebras experience a kind of periodicity over 8 real 
	dimensions in the underlying vector space.
	
	\begin{prop} \label{prop_iso_periodicity}
		For all $n \geq 0$, there are periodicity isomorphisms
		\begin{displaymath}
		\setlength\arraycolsep{2pt}
		\begin{array}{c}
			\mathcal{G}(\mathbb{R}^{n+8,0}) \cong \mathcal{G}(\mathbb{R}^{n,0}) \otimes \mathcal{G}(\mathbb{R}^{8,0}), \\[5pt]
			\mathcal{G}(\mathbb{R}^{0,n+8}) \cong \mathcal{G}(\mathbb{R}^{0,n}) \otimes \mathcal{G}(\mathbb{R}^{0,8}), \\[5pt]
			\mathcal{G}(\mathbb{C}^{n+2}) \cong \mathcal{G}(\mathbb{C}^{n}) \otimes_\mathbb{C} \mathcal{G}(\mathbb{C}^{2}).
		\end{array}
		\end{displaymath}
	\end{prop}
	\begin{proof}
		Using Proposition \ref{prop_iso_tensor} 
		repeatedly, we obtain
		\begin{displaymath}
		\setlength\arraycolsep{2pt}
		\begin{array}{rcl}
			\mathcal{G}(\mathbb{R}^{n+8,0})
				&\cong& \mathcal{G}(\mathbb{R}^{0,n+6}) \otimes \mathcal{G}(\mathbb{R}^{2,0}) \\[5pt]
				&\cong& \mathcal{G}(\mathbb{R}^{n,0}) \otimes \mathcal{G}(\mathbb{R}^{0,2}) \otimes \mathcal{G}(\mathbb{R}^{2,0}) \otimes \mathcal{G}(\mathbb{R}^{0,2}) \otimes \mathcal{G}(\mathbb{R}^{2,0}) \\[5pt]
				&\cong& \mathcal{G}(\mathbb{R}^{n,0}) \otimes \mathcal{G}(\mathbb{R}^{8,0}),
		\end{array}
		\end{displaymath}
		and analogously for the second statement.
		
		For the last statement we take orthonormal bases $\{e_i\}$ of $\mathbb{C}^{n+2}$, $\{\underline{e}_i\}$ of $\mathbb{C}^n$ and $\{\overline{e}_i\}$ of $\mathbb{C}^2$.
		Define a mapping $f\!: \mathbb{C}^{n+2} \to \mathcal{G}(\mathbb{C}^n) \otimes_\mathbb{C} \mathcal{G}(\mathbb{C}^2)$ by
		\begin{displaymath}
		\setlength\arraycolsep{1pt}
		\begin{array}{ccrcll}
			e_j &\ \mapsto\ & i\ \underline{e}_j &\otimes_\mathbb{C}& \overline{e}_1 \overline{e}_2, &\quad j=1,\ldots,n, \\[3pt]
			e_j &\ \mapsto\ & 1 &\otimes_\mathbb{C}& \overline{e}_{j-n}, &\quad j=n+1,n+2,
		\end{array}
		\end{displaymath}
		and extend to an algebra isomorphism as usual.
	\end{proof}
	
	\begin{thm} \label{thm_classification_real}
		We obtain the classification of real geometric algebras as matrix algebras,
		given by Table \ref{table_classification_real} together with the periodicity
		\begin{displaymath}
			\mathcal{G}(\mathbb{R}^{s+8,t}) \cong \mathcal{G}(\mathbb{R}^{s,t+8}) \cong \mathcal{G}(\mathbb{R}^{s,t}) \otimes \mathbb{R}^{16 \times 16}.
		\end{displaymath}
	\end{thm}
	
	\begin{sidewaystable}\centering
		\begin{displaymath}
		\setlength\arraycolsep{2pt}
		\begin{array}{|c||c|c|c|c|c|c|c|c|c|}
			\hline &&&&&&&&& \\[-10pt]
			8	& \mathbb{R}[16]						& \mathbb{R}[16] \oplus \mathbb{R}[16]	& \mathbb{R}[32]						& \mathbb{C}[32]						& \mathbb{H}[32]						& \mathbb{H}[32] \oplus \mathbb{H}[32]	& \mathbb{H}[64]						& \mathbb{C}[128]						& \mathbb{R}[256] \\
			\hline &&&&&&&&& \\[-10pt]
			7	& \mathbb{R}[8] \oplus \mathbb{R}[8]	& \mathbb{R}[16]						& \mathbb{C}[16]						& \mathbb{H}[16]						& \mathbb{H}[16] \oplus \mathbb{H}[16]	& \mathbb{H}[32]						& \mathbb{C}[64]						& \mathbb{R}[128]						& \mathbb{R}[128] \oplus \mathbb{R}[128] \\
			\hline &&&&&&&&& \\[-10pt]
			6	& \mathbb{R}[8]							& \mathbb{C}[8]							& \mathbb{H}[8]							& \mathbb{H}[8] \oplus \mathbb{H}[8]	& \mathbb{H}[16]						& \mathbb{C}[32]						& \mathbb{R}[64]						& \mathbb{R}[64] \oplus \mathbb{R}[64]	& \mathbb{R}[128] \\
			\hline &&&&&&&&& \\[-10pt]
			5	& \mathbb{C}[4]							& \mathbb{H}[4]							& \mathbb{H}[4] \oplus \mathbb{H}[4]	& \mathbb{H}[8]							& \mathbb{C}[16]						& \mathbb{R}[32]						& \mathbb{R}[32] \oplus \mathbb{R}[32]	& \mathbb{R}[64]						& \mathbb{C}[64] \\
			\hline &&&&&&&&& \\[-10pt]
			4	& \mathbb{H}[2]							& \mathbb{H}[2] \oplus \mathbb{H}[2]	& \mathbb{H}[4]							& \mathbb{C}[8]							& \mathbb{R}[16]						& \mathbb{R}[16] \oplus \mathbb{R}[16]	& \mathbb{R}[32]						& \mathbb{C}[32]						& \mathbb{H}[32] \\
			\hline &&&&&&&&& \\[-10pt]
			3	& \mathbb{H} \oplus \mathbb{H}			& \mathbb{H}[2]							& \mathbb{C}[4]							& \mathbb{R}[8]							& \mathbb{R}[8] \oplus \mathbb{R}[8]	& \mathbb{R}[16]						& \mathbb{C}[16]						& \mathbb{H}[16]						& \mathbb{H}[16] \oplus \mathbb{H}[16] \\
			\hline &&&&&&&&& \\[-10pt]
			2	& \mathbb{H}							& \mathbb{C}[2]							& \mathbb{R}[4]							& \mathbb{R}[4] \oplus \mathbb{R}[4]	& \mathbb{R}[8]							& \mathbb{C}[8]							& \mathbb{H}[8]							& \mathbb{H}[8] \oplus \mathbb{H}[8]	& \mathbb{H}[16] \\
			\hline &&&&&&&&& \\[-10pt]
			1	& \mathbb{C}							& \mathbb{R}[2]							& \mathbb{R}[2] \oplus \mathbb{R}[2]	& \mathbb{R}[4]							& \mathbb{C}[4]							& \mathbb{H}[4]							& \mathbb{H}[4] \oplus \mathbb{H}[4]	& \mathbb{H}[8]							& \mathbb{C}[16] \\
			\hline &&&&&&&&& \\[-10pt]
			0	& \mathbb{R}							& \mathbb{R} \oplus \mathbb{R}			& \mathbb{R}[2]							& \mathbb{C}[2]							& \mathbb{H}[2]							& \mathbb{H}[2] \oplus \mathbb{H}[2]	& \mathbb{H}[4]							& \mathbb{C}[8]							& \mathbb{R}[16] \\
			\hline
			\hline
				& 0 & 1 & 2 & 3 & 4 & 5 & 6 & 7 & 8 \\
			\hline
		\end{array}
		\end{displaymath}
		\caption{The algebra $\mathcal{G}(\mathbb{R}^{s,t})$ in the box (s,t), where $\mathbb{F}[N] = \mathbb{F}^{N \times N}$. \label{table_classification_real}}
	\end{sidewaystable}
	
	\begin{proof}
		We have the following easily verified isomorphisms:
		\begin{displaymath}
		\setlength\arraycolsep{1pt}
		\begin{array}{l}
			\mathcal{G}(\mathbb{R}^{1,0}) \cong \mathbb{R} \oplus \mathbb{R}, \\[3pt]
			\mathcal{G}(\mathbb{R}^{0,1}) \cong \mathbb{C}, \\[3pt]
			\mathcal{G}(\mathbb{R}^{2,0}) \cong \mathbb{R}^{2 \times 2}, \\[3pt]
			\mathcal{G}(\mathbb{R}^{0,2}) \cong \mathbb{H}.
		\end{array}
		\end{displaymath}
		Some of these will be explained in detail in Section \ref{sec_lower_dim}.
		We can now work out the cases $(n,0)$ and $(0,n)$ for $n=0,1,\ldots,7$
		in a criss-cross fashion using Proposition \ref{prop_iso_tensor} and the
		tensor algebra isomorphisms
		\begin{displaymath}
		\setlength\arraycolsep{1pt}
		\begin{array}{l}
			\mathbb{C} \otimes_\mathbb{R} \mathbb{C} \cong \mathbb{C} \oplus \mathbb{C}, \\[3pt]
			\mathbb{C} \otimes_\mathbb{R} \mathbb{H} \cong \mathbb{C}^{2 \times 2}, \\[3pt]
			\mathbb{H} \otimes_\mathbb{R} \mathbb{H} \cong \mathbb{R}^{4 \times 4}.
		\end{array}
		\end{displaymath}
		For proofs of these, see e.g. \cite{lawson_michelsohn}.
		With $\mathcal{G}(\mathbb{R}^{1,1}) \cong \mathcal{G}(\mathbb{R}^{2,0})$ 
		and Proposition \ref{prop_iso_tensor} we can then work our way
		through the whole table diagonally. The periodicity follows from Proposition \ref{prop_iso_periodicity}
		and $\mathcal{G}(\mathbb{R}^{8,0}) \cong \mathbb{H} \otimes \mathbb{R}^{2 \times 2} \otimes \mathbb{H} \otimes \mathbb{R}^{2 \times 2} \cong \mathbb{R}^{16 \times 16}$.
	\end{proof}
	
	\noindent
	Because all nondegenerate complex quadratic forms on $\mathbb{C}^n$ are equivalent, 
	the complex version of this theorem turns out to be much simpler.
	
	\begin{thm} \label{thm_classification_complex}
		We obtain the classification of complex geometric algebras as matrix algebras,
		given by
		\begin{displaymath}
		\setlength\arraycolsep{2pt}
		\begin{array}{rcl}
			\mathcal{G}(\mathbb{C}^{0}) &\cong& \mathbb{C}, \\[5pt]
			\mathcal{G}(\mathbb{C}^{1}) &\cong& \mathbb{C} \oplus \mathbb{C},
		\end{array}
		\end{displaymath}
		together with the periodicity
		\begin{displaymath}
			\mathcal{G}(\mathbb{C}^{n+2}) \cong \mathcal{G}(\mathbb{C}^{n}) \otimes_\mathbb{C} \mathbb{C}^{2 \times 2}.
		\end{displaymath}
	\end{thm}

	\noindent
	In other words,
		\begin{displaymath}
		\setlength\arraycolsep{2pt}
		\begin{array}{rcl}
			\mathcal{G}(\mathbb{C}^{2k}) &\cong& \mathbb{C}^{2^k \times 2^k}, \\[5pt]
			\mathcal{G}(\mathbb{C}^{2k+1}) &\cong& \mathbb{C}^{2^k \times 2^k} \oplus \mathbb{C}^{2^k \times 2^k},
		\end{array}
		\end{displaymath}
	for $k=0,1,2,\ldots$

	\begin{proof}
		The isomorphism $\mathcal{G}(\mathbb{C}^n) \cong \mathcal{G}(\mathbb{R}^n) \otimes \mathbb{C}$ gives us
		\begin{displaymath}
		\setlength\arraycolsep{1pt}
		\begin{array}{l}
			\mathcal{G}(\mathbb{C}^0) \cong \mathbb{C} \\[3pt]
			\mathcal{G}(\mathbb{C}^1) \cong (\mathbb{R} \oplus \mathbb{R}) \otimes \mathbb{C} \cong \mathbb{C} \oplus \mathbb{C} \\[3pt]
			\mathcal{G}(\mathbb{C}^2) \cong \mathbb{R}^{2 \times 2} \otimes \mathbb{C} \cong \mathbb{C}^{2 \times 2}.
		\end{array}
		\end{displaymath}
		Then use Proposition \ref{prop_iso_periodicity} for periodicity.
	\end{proof}
	
%

	The periodicity of geometric algebras actually has a number of
	far-reaching consequences. One example is \emph{Bott periodicity}, which simply
	put gives a periodicity in the homotopy groups $\pi_k$ of the unitary, orthogonal and symplectic groups. 
	See \cite{lawson_michelsohn} for proofs using K-theory or \cite{nakahara} for examples.
	
	\newpage

\section{Groups} \label{sec_groups}

	One of the foremost reasons that geometric algebras appear naturally
	in so many areas of mathematics and physics is the fact that they contain
	a number of important groups. These are groups under the geometric product
	and thus lie embedded within the group of invertible elements in $\mathcal{G}$.
	In this section we will discuss the properties of various embedded groups and their relation
	to other familiar transformation groups such as the orthogonal and Lorentz groups.
	We will also introduce a generalized concept of spinor and see how such objects
	are related to the embedded groups.


	
	\begin{defn} \label{def_groups}
		We identify the following groups embedded in $\mathcal{G}$:
		\begin{displaymath}
		\setlength\arraycolsep{2pt}
		\begin{array}{lcll}
			\mathcal{G}^\times &:=& \{ x \in \mathcal{G} : \exists y \in \mathcal{G} : xy = yx = 1 \}
				& \quad\textrm{\emph{the group of invertible elements}} \\[5pt]
			\tilde{\Gamma} &:=& \{ x \in \mathcal{G}^\times : x^\star \mathcal{V} x^{-1} \subseteq \mathcal{V} \}
				& \quad\textrm{\emph{the Lipschitz group}} \\[5pt]
			\Gamma &:=& \{ v_1 v_2 \ldots v_k \in \mathcal{G} : v_i \in \mathcal{V}^\times \}
				& \quad\textrm{\emph{the versor group}} \\[5pt]
			\textrm{Pin} &:=& \{ x \in \Gamma : x x^\dagger = \pm 1 \}
				& \quad\textrm{\emph{the group of unit versors}} \\[5pt]
			\textrm{Spin} &:=& \textrm{Pin} \cap \mathcal{G}^+
				& \quad\textrm{\emph{the group of even unit versors}} \\[5pt]
			\textrm{Spin}^+ &:=& \{ x \in \textrm{Spin} : x x^\dagger = 1 \}
				& \quad\textrm{\emph{the rotor group}}
		\end{array}
		\end{displaymath}
		where $\mathcal{V}^\times := \{ v \in \mathcal{V} : v^2 \neq 0 \}$ is the set of invertible vectors.
	\end{defn}
	
	\noindent
	The versor group $\Gamma$ is the smallest group which contains $\mathcal{V}^\times$. Its elements
	are finite products of invertible vectors called \emph{versors}.
	As is hinted in Definition \ref{def_groups}, many important groups are subgroups of this group.
	One of the central, and highly non-trivial, results of this section 
	is that the versor and Lipschitz groups actually are equal.
	Therefore, $\Gamma$ is also called the Lipschitz group in honor of its creator.
	Sometimes it is also given the name \emph{Clifford group}, but we will, in accordance
	with other conventions, use that name to denote the finite group generated
	by an orthonormal basis.
	
	The Pin and Spin groups are both generated by unit vectors, and in the case
	of Spin, only an even number of such vector factors can be present.
	The elements of Spin$^+$ are called \emph{rotors}. As we will see, these
	groups are intimately connected to orthogonal groups and rotations.
	
	Throughout this section we will always assume that our scalars are real numbers
	unless otherwise stated. This is reasonable both from a geometric
	viewpoint and from the fact that e.g. many complex groups can be represented
	by groups embedded in real geometric algebras.
	Furthermore, we assume that $\mathcal{G}$ is nondegenerate so that we are
	working with a vector space of type $\mathbb{R}^{s,t}$. The corresponding
	groups associated to this space will be denoted Spin$(s,t)$ etc.
	
\subsection{Group actions on $\mathcal{G}$}
	
	In order to understand how groups embedded in a geometric algebra are related to
	more familiar groups of linear transformations, it is necessary 
	to study how groups in $\mathcal{G}$ can act on the 
	vector space $\mathcal{G}$ itself and on the embedded underlying vector space $\mathcal{V}$.
	The following are natural candidates for such actions.
	
	\begin{defn} \label{def_actions}
		Using the geometric product, we have the following natural actions:
		\begin{displaymath}
		\setlength\arraycolsep{2pt}
		\begin{array}{llcll}
			L\!: & \mathcal{G} &\to& \textrm{End}\ \mathcal{G} 
				& \quad\textrm{\emph{left action}} \\[2pt]
				& x &\mapsto& L_x\!: y \mapsto xy \\[5pt]
			R\!: & \mathcal{G} &\to& \textrm{End}\ \mathcal{G} 
				& \quad\textrm{\emph{right action}} \\[2pt]
				& x &\mapsto& R_x\!: y \mapsto yx \\[5pt]
			\textrm{Ad}\!: & \mathcal{G}^\times &\to& \textrm{End}\ \mathcal{G} 
				& \quad\textrm{\emph{adjoint action}} \\[2pt]
				& x &\mapsto& \textrm{Ad}_x\!: y \mapsto xyx^{-1} \\[5pt]
			\widetilde{\textrm{Ad}}\!: & \mathcal{G}^\times &\to& \textrm{End}\ \mathcal{G} 
				& \quad\textrm{\emph{twisted adjoint action}} \\[2pt]
				& x &\mapsto& \widetilde{\textrm{Ad}}_x\!: y \mapsto x^\star yx^{-1}
		\end{array}
		\end{displaymath}
		where $\textrm{End}\ \mathcal{G}$ are the (vector space) 
		endomorphisms of $\mathcal{G}$.
	\end{defn}
	
	\noindent
	Note that $L$ and $R$ are algebra homomorphisms while Ad and $\widetilde{\textrm{Ad}}$ 
	are group homomorphisms.
	These actions give rise to canonical representations of the groups embedded in $\mathcal{G}$.
	The twisted adjoint action takes the graded structure of $\mathcal{G}$ into account and will be
	seen to play a more important role than the normal adjoint action in geometric algebra.
	Using the expansion \eqref{blade_expansion} one can verify that $\textrm{Ad}_x$ is always
	an outermorphism, while in general $\widetilde{\textrm{Ad}}_x$ is not.
	Note, however, that these actions agree on the subgroup of even elements 
	$\mathcal{G}^\times \cap \mathcal{G}^+$.

	\begin{rem}	
		We note that, because the algebra $\mathcal{G}$ is assumed to
		be finite-dimen\-sion\-al, left inverses are always right inverses and \emph{vice versa}.
		This can be seen as follows. First note that the left and right actions are injective.
		Namely, assume that $L_x = 0$. Then $L_x(y) = 0\ \forall y$ and in particular $L_x(1) = x = 0$.
		Suppose now that $xy = 1$ for some $x,y \in \mathcal{G}$.
		But then $L_x L_y = id$, so that $L_y$ is a right inverse to $L_x$. 
		Now, using the dimension theorem 
		\begin{displaymath}
			\dim \ker L_y + \dim \textrm{im}\ L_y = \dim \mathcal{G}
		\end{displaymath}
		with $\ker L_y = 0$, we can conclude that $L_y$ is also a left inverse to $L_x$.
		Hence, $L_y L_x = id \Rightarrow L_{yx-1} = 0$, and $yx = 1$.
	\end{rem}
	
%
	

	
	Let us study the properties of the twisted adjoint action.
	For $v \in \mathcal{V}^\times$ we obtain
	\begin{equation}
		\widetilde{\textrm{Ad}}_v (v) = v^\star v v^{-1} = -v,
	\end{equation}
	and if $w \in \mathcal{V}$ is orthogonal to $v$,
	\begin{equation}
		\widetilde{\textrm{Ad}}_v (w) = v^\star w v^{-1} = -vwv^{-1} = w vv^{-1} = w.
	\end{equation}
	Hence, $\widetilde{\textrm{Ad}}_v$ is a reflection in the hyperplane
	orthogonal to $v$.
	For a general versor $x = u_1 u_2 \ldots u_k \in \Gamma$ we have
	\begin{equation} \label{versor_action}
	\setlength\arraycolsep{2pt}
	\begin{array}{rcl}
		\widetilde{\textrm{Ad}}_x (v) 
			&=& (u_1 \ldots u_k)^\star v (u_1 \ldots u_k)^{-1} 
			= u_1^\star \ldots u_k^\star v u_k^{-1} \ldots u_1^{-1} \\[5pt]
			&=& \widetilde{\textrm{Ad}}_{u_1} \circ \ldots \circ \widetilde{\textrm{Ad}}_{u_k} (v),
	\end{array}
	\end{equation}
	i.e. the twisted adjoint representation (restricted to act only on $\mathcal{V}$
	which is clearly invariant) gives a homomorphism 
	from the versor group into the group of \emph{orthogonal transformations},
	\begin{displaymath}
		\textrm{O}(\mathcal{V},q) := \{f\!: \mathcal{V} \to \mathcal{V} : f\ \textrm{linear bijection s.t.}\ q \circ f = q \}.
	\end{displaymath}
	We have the following fundamental theorem regarding the orthogonal group.
	
	\begin{thm}[Cartan-Dieudonné] \label{thm_cartan_dieudonne}
		Every orthogonal map on a non-degen\-erate space $(\mathcal{V},q)$ is a product of reflections.
		The number of reflections required is at most equal to the dimension of $\mathcal{V}$.
	\end{thm}
	
	\noindent For a constructive proof which works well for arbitrary signatures, see \cite{svensson}.
	
	\begin{cor}
		$\widetilde{\textrm{\emph{Ad}}}\!: \Gamma \to \textrm{\emph{O}}(\mathcal{V},q)$ is surjective.
	\end{cor}
	\begin{proof}
		We know that any $R \in \textrm{O}(\mathcal{V},q)$ can be written
		$R = \widetilde{\textrm{Ad}}_{v_1} \circ \ldots \circ \widetilde{\textrm{Ad}}_{v_k}$
		for some invertible vectors $v_1,\ldots,v_k$, $k \leq n$.
		But then $R = \widetilde{\textrm{Ad}}_{v_1 \ldots v_k}$, where $v_1 v_2 \ldots v_k \in \Gamma$.
	\end{proof}

\subsection{The Lipschitz group}
	
	We saw above that the twisted adjoint representation maps the versor
	group onto the group of orthogonal transformations of $\mathcal{V}$.
	The largest group in $\mathcal{G}$ for which $\widetilde{\textrm{Ad}}$
	forms a representation on $\mathcal{V}$, i.e. leaves $\mathcal{V}$ invariant,
	is per definition the Lipschitz group $\tilde{\Gamma}$.
	We saw from \eqref{versor_action} that $\Gamma \subseteq \tilde{\Gamma}$.
	
	We will now introduce an important function on $\mathcal{G}$, traditionally
	called the \emph{norm function},
	\begin{equation}
	\setlength\arraycolsep{2pt}
	\begin{array}{c}
		N\!: \mathcal{G} \to \mathcal{G}, \\[5pt]
		N(x) := x^\cliffconj x.
	\end{array}
	\end{equation}
	The name is a bit misleading since $N$ is not even
	guaranteed to take values in $\mathbb{R}$. For some special cases of algebras,
	however, it does act as a natural norm and we will see that it can be extended
	in many lower-dimensional algebras where it will act as a kind of determinant.
	Our first main result for this function is that it acts as a determinant on $\tilde{\Gamma}$.
	This will help us prove that $\Gamma = \tilde{\Gamma}$.
	
	\begin{lem} \label{lem_grade_commute}
		Assume that $\mathcal{G}$ is nondegenerate.
		If $x \in \mathcal{G}$ and $x^\star v = vx$ for all $v \in \mathcal{V}$ then
		$x$ must be a scalar.
	\end{lem}
	\begin{proof}
		Using Proposition \ref{prop_trix2} we have that 
		$v \liprod x = 0$ for all $v \in \mathcal{V}$.
		This means that, for a $k$-blade, 
		$(v_1 \wedge \cdots \wedge v_{k-1} \wedge v_k) * x = (v_1 \wedge \cdots \wedge v_{k-1}) * (v_k \liprod x) = 0$
		whenever $k \geq 1$. The nondegeneracy of the scalar product implies that $x$ must have grade 0.
	\end{proof}
	
	\begin{thm} \label{thm_norm}
		The norm function is a group homomorphism $N\!: \tilde{\Gamma} \to \mathbb{R}^\times$.
	\end{thm}
	\begin{proof}
		First note that if $xx^{-1}=1$ then also $x^\star (x^{-1})^\star = 1$ and $(x^{-1})^\dagger x^\dagger = 1$, 
		hence $(x^\star)^{-1}=(x^{-1})^\star$ and $(x^\dagger)^{-1}=(x^{-1})^\dagger$.
		
		Now take $x \in \tilde{\Gamma}$. Then $x^\star v x^{-1} \in \mathcal{V}$ for all $v \in \mathcal{V}$
		and therefore 
		\begin{equation}
			x^\star v x^{-1} = (x^\star v x^{-1})^\dagger = (x^{-1})^\dagger v x^\cliffconj.
		\end{equation}
		This means that $x^\dagger x^\star v = v x^\cliffconj x$, or $N(x)^\star v = v N(x)$.
		By Lemma \ref{lem_grade_commute} we find that $N(x) \in \mathbb{R}$.
		The homomorphism property now follows easily, since for $x,y \in \tilde{\Gamma}$,
		\begin{equation}
			N(xy) = (xy)^\cliffconj xy = y^\cliffconj x^\cliffconj x y = y^\cliffconj N(x) y = N(x)N(y).
		\end{equation}
		Finally, because $1 = N(1) = N(xx^{-1}) = N(x)N(x^{-1})$, we must have that $N(x)$ is nonzero.
	\end{proof}
	
	\begin{lem} \label{lem_lipschitz_surjective}
		The homomorphism $\widetilde{\textrm{\emph{Ad}}}\!: \tilde{\Gamma} \to \textrm{\emph{O}}(\mathcal{V},q)$
		has kernel $\mathbb{R}^\times$.
	\end{lem}
	\begin{proof}
		We first prove that $\widetilde{\textrm{Ad}}_x$ is orthogonal for $x \in \tilde{\Gamma}$.
		Note that, for $v \in \mathcal{V}$,
		\begin{equation}
		\setlength\arraycolsep{2pt}
		\begin{array}{rcl}
			N(\widetilde{\textrm{Ad}}_x(v)) &=& N(x^\star v x^{-1}) 
				= (x^\star v x^{-1})^\cliffconj x^\star v x^{-1} \\[5pt]
				&=& (x^{-1})^\cliffconj v^\cliffconj x^\dagger x^\star v x^{-1}
				= (x^{-1})^\cliffconj v^\cliffconj N(x)^\star v x^{-1} \\[5pt]
				&=& (x^{-1})^\cliffconj v^\cliffconj v x^{-1} N(x)^\star
				= N(v) N(x^{-1}) N(x) = N(v).
		\end{array}
		\end{equation}
		Then, since $N(v) = v^\cliffconj v = -v^2$, we have that $\widetilde{\textrm{Ad}}_x(v)^2 = v^2$.
		
		Now, if $\widetilde{\textrm{Ad}}_x = id$ then $x^\star v = v x$ for all $v \in \mathcal{V}$
		and by Lemma \ref{lem_grade_commute}
		we must have $x \in \mathbb{R} \cap \tilde{\Gamma} = \mathbb{R}^\times$.
	\end{proof}
	
	\noindent We finally obtain 
	
	\begin{thm} \label{thm_lipschitz}
		We have that $\Gamma = \tilde{\Gamma}$.
	\end{thm}
	\begin{proof}
		We saw earlier that $\Gamma \subseteq \tilde{\Gamma}$.
		Take $x \in \tilde{\Gamma}$. By the above lemma we have $\widetilde{\textrm{Ad}}_x \in \textrm{O}(\mathcal{V},q)$.
		Using the corollary to Theorem \ref{thm_cartan_dieudonne} we then find that
		$\widetilde{\textrm{Ad}}_x = \widetilde{\textrm{Ad}}_y$ for some $y \in \Gamma$.
		Then $\widetilde{\textrm{Ad}}_{xy^{-1}} = id$ and $xy^{-1} = \lambda \in \mathbb{R}^\times$.
		Hence, $x = \lambda y \in \Gamma$.
	\end{proof}
	
\subsection{Properties of Pin and Spin groups}
	
	From the discussion above followed that $\widetilde{\textrm{Ad}}$ gives a
	surjective homomorphism from the versor, or Lipschitz, group $\Gamma$ to the orthogonal
	group. The kernel of this homomorphism is the set of invertible scalars.
	Because the Pin and Spin groups consist of normalized versors ($N(x) = \pm 1$) we
	find the following
	
	\begin{thm} \label{thm_double_cover}
		The homomorphisms
		\begin{displaymath}
		\setlength\arraycolsep{2pt}
		\begin{array}{llcl}
			\widetilde{\textrm{\emph{Ad}}}\!: & \textrm{\emph{Pin}} (s,t) &\to& \textrm{\emph{O}} (s,t) \\[5pt]
			\widetilde{\textrm{\emph{Ad}}}\!: & \textrm{\emph{Spin}}(s,t) &\to& \textrm{\emph{SO}}(s,t) \\[5pt]
			\widetilde{\textrm{\emph{Ad}}}\!: & \textrm{\emph{Spin}}^+(s,t) &\to& \textrm{\emph{SO}}^+(s,t)
		\end{array}
		\end{displaymath}
		are surjective with kernel $\{\pm 1\}$.
	\end{thm}
	
	\noindent
	The homomorphism onto the \emph{special orthogonal group},
	\begin{displaymath}
		\textrm{SO}(\mathcal{V},q) := \{f \in \textrm{O}(\mathcal{V},q) : \det f = 1 \}
	\end{displaymath}
	follows since it
	is generated by an \emph{even} number of reflections. $\textrm{SO}^+$ denotes
	the connected component of SO containing the identity. This will soon be explained.
	
	In other words, the Pin and Spin groups are two-sheeted coverings of the orthogonal groups.
	Furthermore, we have the following relations between these groups.
	
	Take a unit versor $\psi = u_1 u_2 \ldots u_k \in \textrm{Pin}(s,t)$.
	If $\psi$ is odd we can always multiply by a unit vector $e$ so that
	$\psi = \pm \psi ee$ and $\pm \psi e \in \textrm{Spin}(s,t)$.
	Furthermore, when the signature is euclidean we have $\psi \psi^\dagger = 1$ for all unit versors.
	The same holds for \emph{even} unit versors in anti-euclidean spaces since the signs cancel out.
	Hence, $\textrm{Spin} = \textrm{Spin}^+$ unless there is mixed signature. But in that case we can
	find two orthogonal unit vectors $e_+,e_-$ such that $e_+^2 = 1$ and $e_-^2 = -1$.
	Since $e_+e_-(e_+e_-)^\dagger = -1$ we then have that $\psi = \psi(e_+e_-)^2$, where
	$\psi e_+e_- (\psi e_+e_-)^\dagger = 1$ if $\psi \psi^\dagger = -1$.
	
	Summing up, we have that, for $s,t \geq 1$ and any pair of orthogonal vectors $e_+,e_-$ such that
	$e_+^2 = 1$, $e_-^2 = -1$,
	\begin{displaymath}
	\setlength\arraycolsep{2pt}
	\begin{array}{rcl}
		\textrm{Pin}(s,t) &=& \textrm{Spin}^+(s,t) \cdot \{ 1, e_+, e_-, e_+ e_- \}, \\[5pt]
		\textrm{Spin}(s,t) &=& \textrm{Spin}^+(s,t) \cdot \{ 1, e_+ e_- \}.
	\end{array}
	\end{displaymath}
	Otherwise,
	\begin{displaymath}
		\textrm{Pin}(s,t) = \textrm{Spin}^{(+)}(s,t) \cdot \{ 1, e \}
	\end{displaymath}
	for any $e \in \mathcal{V}$ such that $e^2 = \pm 1$.
	From the isomorphism $\mathcal{G}^+(\mathbb{R}^{s,t}) \cong \mathcal{G}^+(\mathbb{R}^{t,s})$ we also have
	the signature symmetry $\textrm{Spin}^{(+)}(s,t) \cong \textrm{Spin}^{(+)}(t,s)$.
	In all cases,
	\begin{displaymath}
		\Gamma(s,t) = \mathbb{R}^\times \cdot \textrm{Pin}(s,t).
	\end{displaymath}
	
	From these considerations it is sufficient to study the properties of the
	rotor groups in order to understand the Pin, Spin and orthogonal groups.
	Fortunately, it turns out that the rotor groups have very convenient topological
	features.
	
	\begin{thm} \label{thm_spin_connected}
		The groups $\textrm{\emph{Spin}}^+(s,t)$ are pathwise connected for $s \geq 2$ or $t \geq 2$.
	\end{thm}
	\begin{proof}
		Pick a rotor $R \in \textrm{Spin}^+(s,t)$, where $s$ or $t$ is greater than one.
		Then $R = v_1 v_2 \ldots v_{2k}$ with an even number of $v_i \in \mathcal{V}$
		such that $v_i^2=1$ and an even number such that $v_i^2=-1$.
		Note that for any two invertible vectors $a,b$ we have
		$ab = aba^{-1}a = b'a$, where $b'^2 = b^2$.
		Hence, we can rearrange the vectors so that those with positive square come first, i.e.
		\begin{equation}
			R = a_1 b_1 \ldots a_p b_p a_1' b_1' \ldots a_q' b_q' = R_1 \ldots R_p R_1' \ldots R_q',
		\end{equation}
		where $a_i^2 = b_i^2 = 1$ and $R_i = a_i b_i = a_i * b_i + a_i \wedge b_i$
		are so called \emph{simple rotors} which are connected to either 1 or -1.
		This holds because $1 = R_i R_i^\dagger = (a_i * b_i)^2 - (a_i \wedge b_i)^2$,
		so we can, as is easily verified, write 
		$R_i = \pm e^{\phi_i a_i \wedge b_i}$ for some $\phi_i \in \mathbb{R}$
		(exponentials of bivectors will be treated shortly).
		Depending on the signature of the plane associated to $a_i \wedge b_i$, 
		i.e. on the sign of $(a_i \wedge b_i)^2 \in \mathbb{R}$,
		the set $e^{\mathbb{R} a_i \wedge b_i} \subseteq \textrm{Spin}^+$
		forms either a circle, a line or a hyperbola. In any case, it goes through the unit element.
		Finally, since $s>1$ or $t>1$ we can connect -1 to 1 with for example
		the circle $e^{\mathbb{R}e_1e_2}$, where $e_1, e_2$ are two orthonormal basis elements with
		the same signature.
	\end{proof}
	
	\noindent
	Continuity of $\widetilde{\textrm{Ad}}$ now implies that the set of rotations
	represented by rotors, i.e. $\textrm{SO}^+$, forms a continuous subgroup containing the identity.
	For euclidean and lorentzian signatures, we have an even simpler situation.
	
	\begin{thm} \label{thm_spin_simply_conn}
		The groups $\textrm{\emph{Spin}}^+(s,t)$ are simply connected 
		for $(s,t) = (n,0)$, $(0,n)$, $(1,n)$ or $(n,1)$, where $n \geq 3$.
		Hence, these are the universal covering groups of $\textrm{\emph{SO}}^+(s,t)$.
	\end{thm}
	
	\noindent
	This follows because $\pi_1\big(\textrm{SO}^+(s,t)\big) = \mathbb{Z}_2$ for these
	signatures. See e.g. \cite{lawson_michelsohn} for details.
	
	This sums up the the situation nicely for higher-dimensional euclidean
	and lorentzian spaces: The Pin group, which is a double-cover
	of the orthogonal group, consists of two or four simply connected components.
	These components are copies of the rotor group.
	In physics-terminology these components correspond to time and parity reflections.
	
	It is also interesting to relate the rotor group to its Lie algebra,
	which is actually the bivector space $\mathcal{G}^2$ with the usual commutator $[\cdot,\cdot]$.
	This follows because there is a Lie algebra isomorphism between $\mathcal{G}^2$
	and the algebra of antisymmetric transformations of $\mathcal{V}$, given by
	\begin{equation}
		f = [B,\cdot] \quad \leftrightarrow \quad B = \frac{1}{2} \sum_{i,j} e^i*f(e^j)\ e_i \wedge e_j,
	\end{equation}
	where $\{e_i\}$ is some general basis of $\mathcal{V}$.
	One verifies, by expanding in the geometric product, that
	\begin{equation}
		\frac{1}{2} [e_i \wedge e_j, e_k \wedge e_l] = e_j*e_k e_i \wedge e_l - e_j*e_l e_i \wedge e_k + e_i*e_l e_j \wedge e_k - e_i*e_k e_j \wedge e_l.
	\end{equation}
	Actually, this bivector Lie algebra is more general than it might first
	seem. Doran and Lasenby (see \cite{doran_et_al} or \cite{doran_lasenby}) have shown
	that the Lie algebra $\mathfrak{gl}$ of the general linear group can be represented as
	a bivector algebra. From the fact that any finite-dimensional Lie algebra
	has a faithful finite-dimensional representation
	(Ado's Theorem for characteristic zero, Iwasawa's Theorem for nonzero characteristic,
	see e.g. \cite{jacobson})
	we have that \emph{any} finite-dimensional real or complex Lie algebra can be represented 
	as a bivector algebra.
	
	We define the exponential of a multivector $x \in \mathcal{G}$ as the usual power series
	\begin{equation} \label{exp_definition}
		e^x := \sum_{k=0}^\infty \frac{x^k}{k!}.
	\end{equation}
	Since $\mathcal{G}$ is finite-dimensional we have the following 
	for any choice of norm on $\mathcal{G}$.
	
	\begin{prop} \label{prop_exp_property}
		The sum in \eqref{exp_definition} converges for all $x \in \mathcal{G}$ and
		\begin{displaymath}
			e^x e^{-x} = e^{-x} e^x = 1.
		\end{displaymath}
	\end{prop}
	
	\noindent
	The following now holds for any signature.
	
	\begin{thm} \label{thm_exp_in_spin}
		For any bivector $B \in \mathcal{G}^2$ we have that $\pm e^B \in \textrm{\emph{Spin}}^+$.
	\end{thm}
	\begin{proof}
		It is obvious that $\pm e^B$ is an even multivector and that $e^B (e^B)^\dagger = e^B e^{-B} = 1$.
		Hence, it is sufficient to prove that $\pm e^B \in \Gamma$, or by Theorem \ref{thm_lipschitz},
		that $e^B \mathcal{V} e^{-B} \subseteq \mathcal{V}$. This can be done by considering
		derivatives of the function $f(t) := e^{tB} v e^{-tB}$ for $v \in \mathcal{V}$.
		See e.g. \cite{svensson} or \cite{riesz} for details.
	\end{proof}
	
	\noindent
	The converse is true for (anti-) euclidean and lorentzian spaces.
	
	\begin{thm} \label{thm_spin_is_exp}
		For $(s,t) = (n,0)$, $(0,n)$, $(1,n)$ or $(n,1)$, we have 
		\begin{displaymath}
			\textrm{\emph{Spin}}^+(s,t) = \pm e^{\mathcal{G}^2(\mathbb{R}^{s,t})},
		\end{displaymath}
		i.e. any rotor can be written as (minus) the exponential of a bivector.
		The minus sign is only required in the cases 
		$(0,0)$, $(1,0)$, $(0,1)$, $(1,1)$, $(1,2)$, $(2,1)$, $(1,3)$ and $(3,1)$.
	\end{thm}
	
	\noindent
	The proof can be found in \cite{riesz}. Essentially, it relies on the fact
	that any isometry of an euclidean or lorentzian space can be generated
	by a single infinitesimal transformation. This holds for these spaces \emph{only},
	so that for example $\textrm{Spin}^+(2,2) \neq \pm e^{\mathcal{G}^2(\mathbb{R}^{2,2})}$,
	where for instance $\pm e_1e_2e_3e_4 e^{\beta(e_1e_2 + 2e_1e_4 + e_3e_4)}$, $\beta>0$, 
	cannot be reduced to a single exponential.
	
	
\subsection{Spinors}
	
	
	Spinors are objects which originally appeared in physics in the early
	days of quantum mechanics, but have by now made their way into other fields as well,
	such as differential geometry and topology.
	They are traditionally represented as elements of a complex vector space
	since it was natural to add and subtract them and scale them by complex amplitudes
	in their original physical applications.
	However, what really characterizes them as spinors is the fact that they 
	can be acted upon by rotations, together with their rather special transformation
	properties under such rotations. While a rotation needs just one revolution
	to get back to the identity in a vector representation, it takes two
	revolutions to come back to the identity in a spinor representation.
	This is exactly the behaviour experienced by rotors, since they transform vectors double-sidedly
	with the (twisted) adjoint action as
	\begin{equation} \label{rotor_double_action}
		v \mapsto \psi v \psi^\dagger \mapsto \phi \psi v \psi^\dagger \phi^\dagger = (\phi\psi)v(\phi\psi)^\dagger,
	\end{equation}
	where we applied consecutive rotations represented by rotors $\psi$ and $\phi$.
	Note that the rotors themselves transform according to
	\begin{equation}
		1 \mapsto \psi \mapsto \phi\psi,
	\end{equation}
	that is single-sidedly with the left action.
	
	This hints that rotors could represent some form of spinors in geometric algebra.
	However, the rotors form a group and not a vector space, so we need
	to consider a possible enclosing vector space with similar properties.
	Some authors (see e.g. Hestenes \cite{hestenes_sobczyk}) have 
	considered so called \emph{operator spinors}, which are
	general even multivectors that leave $\mathcal{V}$ invariant 
	under the double-sided action \eqref{rotor_double_action},
	i.e. elements of
	\begin{equation}
		\Sigma := \{ \Psi \in \mathcal{G}^+ : \Psi \mathcal{V} \Psi^\dagger \subseteq \mathcal{V} \}.
	\end{equation}
	Since $\Psi v \Psi^\dagger$ is both odd and self-reversing, it is an element
	of $\bigoplus_{k=0}^{\infty} \mathcal{G}^{4k+1}$. Therefore $\Sigma$
	is only guaranteed to be a vector space for $\dim \mathcal{V} \leq 4$, 
	where it coincides with $\mathcal{G}^+$.
	
	For a general spinor space embedded in $\mathcal{G}$, we seek a subspace
	$\Sigma \subseteq \mathcal{G}$ that is invariant under left action by rotors,
	i.e. such that
	\begin{equation}
		\Psi \in \Sigma, \quad \psi \in \textrm{Spin}^+(s,t) \quad \Rightarrow \quad \psi\Psi \in \Sigma.
	\end{equation}
	One obvious and most general choice of such a spinor space is the whole space $\mathcal{G}$. 
	However, we will soon see from examples that spinors in lower dimensions
	are best represented by another natural suggestion, namely the set of even multivectors. 
	Hence, we follow Francis and Kosowsky \cite{francis_kosowsky}
	and define the space of spinors $\Sigma$ for arbitrary dimensions as the even subalgebra $\mathcal{G}^+$.
	
	Note that action by \emph{rotors}, i.e. $\textrm{Spin}^+$ instead of Spin, 
	ensures that $\Psi \Psi^\dagger$ remains invariant under right action
	and $\Psi^\dagger \Psi$ remains invariant under left action on $\Psi$.
	In lower dimensions these are invariant under both actions, so they are good
	candidates for invariant or \emph{observable} quantities in physics.
	Also note that if $\dim \mathcal{V} \leq 4$ then 
	the set of unit spinors and the set of rotors coincide, i.e.
	\begin{equation} \label{lowerdim_spinor_rotor}
		\textrm{Spin}^+(s,t) = \{\psi \in \mathcal{G}^+ : \psi \psi^\dagger = 1\} \quad \textrm{for}\ s+t \leq 4.
	\end{equation}
	This follows because then $\psi$ is invertible and $\psi^\star v \psi^{-1} = \psi v \psi^\dagger$
	is both odd and self-reversing, hence a vector.
	Thus, $\psi$ lies in $\tilde{\Gamma}$ and is therefore an even unit versor.
	
	A popular alternative to the above is to consider spinors as elements in 
	left-sided ideals of (mostly complex) geometric algebras. This is the view 
	which is closest related to the original complex vector space picture. 
	We will see an example of how these views are related in the case of the Dirac algebra.
	A motivation for this definition of spinor is that it admits a straightforward
	transition to basis-independent spinors, so called \emph{covariant spinors}.
	This is required for a treatment of spinor fields on curved
	manifolds, i.e. in the gravitational setting.
	However, covariant spinors lack the clearer geometrical picture provided by 
	operator and even subalgebra spinors. Furthermore, by reconsidering the definition and 
	interpretation of these geometric spinors, it is possible to deal with 
	basis-independence also in this case.
	These and other properties of spinors related to physics will be
	discussed in Section \ref{sec_spinors}.

	\newpage

\section{A study of lower-dimensional algebras} \label{sec_lower_dim}

	We have studied the structure of geometric algebras in general
	and saw that they are related to many other familiar 
	algebras and groups. We will now go through a number of lower-dimensional
	examples in detail to see just how structure-rich these algebras are.
	Although we go no higher than to a five-dimensional base vector space,
	we manage to find a variety of structures related to physics.
	
	We choose to focus around the groups and spinors of these example algebras.
	It is highly recommended that the reader aquires a more complete understanding
	of at least the plane, space and spacetime algebras from other sources.
	Good introductions aimed at physicists can be found in \cite{doran_lasenby}
	and \cite{gull_lasenby_doran}. A more mathematical treatment is given in \cite{lounesto}.

\subsection{$\mathcal{G}(\mathbb{R}^{1})$}
	
	Since $\mathcal{G}(\mathbb{R}^{0,0}) = \mathbb{R}$ is just the field of real numbers, which should be familiar,
	we start instead with $\mathcal{G}(\mathbb{R}^{1})$, the geometric algebra of the real line. 
	Let $\mathbb{R}^1$ be spanned by one basis element $e$ such that $e^2 = 1$. Then
	\begin{equation}
		\mathcal{G}(\mathbb{R}^{1}) = \textrm{Span}_\mathbb{R} \{1,e\}.
	\end{equation}
	This is a commutative algebra with pseudoscalar $e$. One easily finds the invertible
	elements $\mathcal{G}^\times$ by considering the norm function,
	which with a one-dimensional $\mathcal{V}$ is given by
	\begin{equation}
		N_1(x) := x^\cliffconj x = x^\star x.
	\end{equation}
	For an arbitrary element $x = \alpha + \beta e$ then
	\begin{equation}
		N_1(\alpha + \beta e) = (\alpha - \beta e)(\alpha + \beta e) = \alpha^2 - \beta^2 \in \mathbb{R}.
	\end{equation}
	When $N_1(x) \neq 0$ we find that $x$ has an inverse $x^{-1} = \frac{1}{N_1(x)} x^\star = \frac{1}{\alpha^2-\beta^2}(\alpha - \beta e)$.
	Hence,
	\begin{equation}
		\mathcal{G}^\times(\mathbb{R}^{1}) = \{ x \in \mathcal{G} : N_1(x) \neq 0 \} 
			= \{ \alpha + \beta e \in \mathcal{G} : \alpha^2 \neq \beta^2 \}.
	\end{equation}
	Note also that $N_1(xy) = x^\star y^\star x y = N_1(x)N_1(y)$ for all $x,y \in \mathcal{G}$
	since the algebra is commutative.
	
	The other groups are rather trivial in this space. Because $\mathcal{G}^+ = \mathbb{R}$,
	we have
	\begin{displaymath}
	\setlength\arraycolsep{2pt}
	\begin{array}{rcl}
		\textrm{Spin}^{(+)}(1,0) &=& \{ 1, -1 \}, \\[5pt]
		\textrm{Pin}(1,0) &=& \{ 1, -1, e, -e \}, \\[5pt]
		\Gamma(1,0) &=& \mathbb{R}^\times\ \bigsqcup\ \mathbb{R}^\times e,
	\end{array}
	\end{displaymath}
	where we write $\bigsqcup$ to emphasize a \emph{disjoint} union.
	The spinors in this algebra are just real scalars.
	
\subsection{$\mathcal{G}(\mathbb{R}^{0,1}) \stackrel{\sim}{=} \mathbb{C}$ - The complex numbers}
	
	As one might have noticed from previous discussions on isomorphisms,
	the complex numbers are in fact a real geometric algebra.
	Let $i$ span a one-dimensional anti-euclidean space and be normalized to $i^2 = -1$.
	Then
	\begin{equation}
		\mathcal{G}(\mathbb{R}^{0,1}) = \textrm{Span}_\mathbb{R} \{1,i\} \cong \mathbb{C}.
	\end{equation}
	This is also a commutative algebra, but, unlike the previous example, this is a \emph{field}
	since every nonzero element is invertible. The norm function is an
	actual norm (squared) in this case,
	\begin{equation}
		N_1(\alpha + \beta i) = (\alpha - \beta i)(\alpha + \beta i) = \alpha^2 + \beta^2 \in \mathbb{R}^+,
	\end{equation}
	namely the modulus of the complex number. Note that the grade involution
	represents the complex conjugate and $x^{-1} = \frac{1}{N_1(x)} x^\star$ as above.
	The relevant groups are
	\begin{displaymath}
	\setlength\arraycolsep{2pt}
	\begin{array}{rcl}
		\mathcal{G}^\times(\mathbb{R}^{0,1}) &=& \mathcal{G} \!\smallsetminus\! \{0\}, \\[5pt]
		\textrm{Spin}^{(+)}(0,1) &=& \{ 1, -1 \}, \\[5pt]
		\textrm{Pin}(0,1) &=& \{ 1, -1, i, -i \}, \\[5pt]
		\Gamma(0,1) &=& \mathbb{R}^\times\ \bigsqcup\ \mathbb{R}^\times i.
	\end{array}
	\end{displaymath}
	The spinor space is still given by $\mathbb{R}$.
	
\subsection{$\mathcal{G}(\mathbb{R}^{0,0,1})$}

	We include this as our only example of a degenerate algebra,
	just to see what such a situation might look like.
	Let $n$ span a one-dimensional space with quadratic form $q=0$.
	Then
	\begin{equation}
		\mathcal{G}(\mathbb{R}^{0,0,1}) = \textrm{Span}_\mathbb{R} \{1,n\}
	\end{equation}
	and $n^2=0$. The norm function depends only on the scalar part,
	\begin{equation}
		N_1(\alpha + \beta n) = (\alpha - \beta n)(\alpha + \beta n) = \alpha^2 \in \mathbb{R}^+.
	\end{equation}
	An element is invertible if and only if the scalar part is nonzero.
	Since no vectors are invertible, we are left with only the empty product
	in the versor group. This gives
	\begin{displaymath}
	\setlength\arraycolsep{2pt}
	\begin{array}{rcl}
		\mathcal{G}^\times(\mathbb{R}^{0,0,1}) &=& \{ \alpha + \beta n \in \mathcal{G} : \alpha \neq 0 \}, \\[5pt]
		\Gamma &=& \{1\}.
	\end{array}
	\end{displaymath}
	Note, however, that for $\alpha \neq 0$
	\begin{equation}
		(\alpha + \beta n)^\star n (\alpha + \beta n)^{-1} = (\alpha - \beta n) n {\textstyle \frac{1}{\alpha^2}} (\alpha - \beta n) = n,
	\end{equation}
	so the Lipschitz group is
	\begin{equation}
		\tilde{\Gamma} = \mathcal{G}^\times \neq \Gamma.
	\end{equation}
	This shows that the assumption on nondegeneracy was necessary in the
	discussion about the Lipschitz group in Section \ref{sec_groups}.
	
\subsection{$\mathcal{G}(\mathbb{R}^{2})$ - The plane algebra}
	
	Our previous examples were rather trivial, but we now come to our first 
	really interesting case, namely the geometric algebra of the euclidean plane.
	Let $\{e_1, e_2\}$ be an orthonormal basis of $\mathbb{R}^2$ and consider
	\begin{equation}
		\mathcal{G}(\mathbb{R}^2) = \textrm{Span}_\mathbb{R} \{1,\ e_1,e_2,\ e_1e_2\}.
	\end{equation}
	An important feature of this algebra is that the pseudoscalar $I := e_1e_2$ squares to $-1$.
	This makes the even subalgebra isomorphic to the complex numbers, 
	in correspondence with the relation $\mathcal{G}^+(\mathbb{R}^{2}) \cong \mathcal{G}(\mathbb{R}^{0,1})$.
	
	Let us find the invertible elements of the plane algebra.
	For two-dimensional algebras we use the original norm function
	\begin{equation}
		N_2(x) := x^\cliffconj x
	\end{equation}
	since it satisfies $N_2(x)^\cliffconj = N_2(x)$ for all $x \in \mathcal{G}$.
	The sign relations for involutions in Table \ref{table_involutions} then require
	this to be a scalar, so we have a map
	\begin{equation}
		N_2 \!: \mathcal{G} \to \mathcal{G}^0 = \mathbb{R}.
	\end{equation}
	For an arbitrary element $x = \alpha + a_1e_1 + a_2e_2 + \beta I \in \mathcal{G}$ we have
	\begin{equation}
	\setlength\arraycolsep{2pt}
	\begin{array}{rcl}
		N_2(x) &=& (\alpha - a_1e_1 - a_2e_2 - \beta I)(\alpha + a_1e_1 + a_2e_2 + \beta I) \\[5pt]
			&=& \alpha^2 - a_1^2 - a_2^2 + \beta^2.
	\end{array}
	\end{equation}
	Furthermore, $N_2(x^\cliffconj) = N_2(x)$ and $N_2(xy) = N_2(x)N_2(y)$ for all $x,y$.
	Proceeding as in the one-dimensional case, we find that $x$ has
	an inverse $x^{-1} = \frac{1}{N_2(x)} x^\cliffconj$ if and only if $N_2(x) \neq 0$, i.e.
	\begin{equation}
	\setlength\arraycolsep{2pt}
	\begin{array}{rcl}
		\mathcal{G}^\times(\mathbb{R}^{2}) &=& \{ x \in \mathcal{G} : N_2(x) \neq 0 \} \\[5pt]
			&=& \{ \alpha + a_1e_1 + a_2e_2 + \beta I \in \mathcal{G} : \alpha^2 - a_1^2 - a_2^2 + \beta^2 \neq 0 \}.
	\end{array}
	\end{equation}
	
	For $x = \alpha + \beta I$ in the even subspace we have $x^\cliffconj = x^\dagger = \alpha - \beta I$, 
	so the Clifford conjugate acts as
	complex conjugate in this case. Again, the norm function (here $N_2$) acts as modulus squared.
	We find that the rotor group, i.e. the group of even unit versors, corresponds to
	the group of unit complex numbers,
	\begin{equation}
		\textrm{Spin}^{(+)}(2,0) = e^{\mathbb{R}I} \cong \textrm{U}(1).
	\end{equation}
	Note that, because $e_1$ and $I$ anticommute, 
	\begin{equation}
		e^{\varphi I} e_1 e^{-\varphi I} = e_1 e^{-2\varphi I} = e_1(\cos 2\varphi - I \sin 2\varphi)
			= e_1 \cos 2\varphi - e_2 \sin 2\varphi,
	\end{equation}
	so a rotor $\pm e^{-\varphi I/2}$ represents a counter-clockwise\footnote{Assuming, 
	of course, that $e_1$ points at 3 o'clock and $e_2$ at 12 o'clock.}
	rotation in the plane by an angle $\varphi$.
	The Pin group is found by picking for example $e_1$;
	\begin{displaymath}
	\setlength\arraycolsep{2pt}
	\begin{array}{rcl}
		\textrm{Pin}(2,0) &=& e^{\mathbb{R}I}\ \bigsqcup\ e^{\mathbb{R}I} e_1, \\[5pt]
		\Gamma(2,0) &=& \mathbb{R}^\times e^{\mathbb{R}I}\ \bigsqcup\ \mathbb{R}^\times e^{\mathbb{R}I} e_1.
	\end{array}
	\end{displaymath}
	
	As we saw above, the spinors of $\mathbb{R}^2$ are nothing but complex numbers.
	We can write any spinor or complex number $\Psi \in \mathcal{G}^+$ in
	the polar form $\Psi = \rho e^{\varphi I}$, which is just a rescaled rotor.
	The spinor action
	\begin{equation}
		a \mapsto \Psi a \Psi^\dagger = \rho^2 e^{\varphi I} a e^{-\varphi I}
	\end{equation}
	then gives a geometric interpretation of the spinor $\Psi$ as an
	operation to rotate by an angle $-2\varphi$ and scale by $\rho^2$.
	
\subsection{$\mathcal{G}(\mathbb{R}^{0,2}) \stackrel{\sim}{=} \mathbb{H}$ - The quaternions}

	The geometric algebra of the anti-euclidean plane is isomorphic 
	to Hamilton's \emph{quaternion algebra} $\mathbb{H}$.
	This follows by taking an orthonormal basis $\{i,j\}$ of $\mathbb{R}^{0,2}$ and considering
	\begin{equation}
		\mathcal{G}(\mathbb{R}^{0,2}) = \textrm{Span}_\mathbb{R} \{1,\ i,j,\ k\},
	\end{equation}
	where $k := ij$ is the pseudoscalar. We then have the classic identities defining quaternions,
	\begin{equation}
		i^2 = j^2 = k^2 = ijk = -1.
	\end{equation}
	We write an arbitrary quaternion as $x = \alpha + a_1i + a_2j + \beta k$.
	The Clifford conjugate acts as the \emph{quaternion conjugate}, 
	$x^\cliffconj = \alpha - a_1i - a_2j - \beta k$.
	
	The norm function $N_2$ has the same properties as in the euclidean algebra, 
	but in this case it once again represents the square of an actual norm,
	namely the \emph{quaternion norm},
	\begin{equation}
		N_2(x) = \alpha^2 + a_1^2 + a_2^2 + \beta^2.
	\end{equation}
	Just as in the complex case then, all nonzero elements are invertible,
	\begin{equation}
		\mathcal{G}^\times(\mathbb{R}^{0,2}) = \{ x \in \mathcal{G} : N_2(x) \neq 0 \} = \mathcal{G} \smallsetminus \{0\}.
	\end{equation}
	The even subalgebra is also in this case isomorphic to the complex
	numbers, so the spinors and groups $\Gamma$, Pin and Spin are 
	no different than in the euclidean case.

\subsection{$\mathcal{G}(\mathbb{R}^{1,1})$}

	This is our simplest example of a lorentzian algebra. An orthonormal
	basis $\{e_+,e_-\}$ of $\mathbb{R}^{1,1}$ consists of a timelike vector, $e_+^2 = 1$,
	and a spacelike vector, $e_-^2 = -1$. 
	In general, a vector (or blade) $v$ is called \emph{timelike} if $v^2>0$, 
	\emph{spacelike} if $v^2<0$, and \emph{lightlike} or \emph{null} if $v^2=0$.
	This terminology is taken from relativistic physics. 
	The two-dimensional lorentzian algebra is given by
	\begin{equation}
		\mathcal{G}(\mathbb{R}^{1,1}) = \textrm{Span}_\mathbb{R} \{1,\ e_+,e_-,\ e_+e_-\}.
	\end{equation}
	The group of invertible elements is as usual given by
	\begin{equation}
	\setlength\arraycolsep{2pt}
	\begin{array}{rcl}
		\mathcal{G}^\times(\mathbb{R}^{2}) &=& \{ x \in \mathcal{G} : N_2(x) \neq 0 \} \\[5pt]
			&=& \{ \alpha + a_+e_+ + a_-e_- + \beta I \in \mathcal{G} : \alpha^2 - a_+^2 + a_-^2 - \beta^2 \neq 0 \}.
	\end{array}
	\end{equation}
	The pseudoscalar $I := e_+e_-$ squares to the identity in this case and
	the even subalgebra is therefore $\mathcal{G}^+(\mathbb{R}^{1,1}) \cong \mathcal{G}(\mathbb{R}^1)$.
	This has as an important consequence that the rotor group is fundamentally
	different from the euclidean case,
	\begin{equation}
		\textrm{Spin}^+(1,1) = \{ \psi = \alpha + \beta I \in \mathcal{G}^+ : \psi^\dagger \psi = \alpha^2 - \beta^2 = 1 \} 
			= \pm e^{\mathbb{R}I}.
	\end{equation}
	This is a pair of disjoint hyperbolas passing through the points 1 and -1, respectively.
	The Spin group consists of four such hyperbolas and the Pin group of eight,
	\begin{displaymath}
	\setlength\arraycolsep{2pt}
	\begin{array}{rcl}
		\textrm{Spin}(1,1) &=& \pm e^{\mathbb{R}I}\ \bigsqcup\ \pm e^{\mathbb{R}I}I, \\[5pt]
		\textrm{Pin}(1,1) &=& \pm e^{\mathbb{R}I}\ \bigsqcup\ \pm e^{\mathbb{R}I}e_+\ \bigsqcup\ \pm e^{\mathbb{R}I}e_-\ \bigsqcup\ \pm e^{\mathbb{R}I}I, \\[5pt]
		\Gamma(1,1) &=& \mathbb{R}^\times e^{\mathbb{R}I}\ \bigsqcup\ \mathbb{R}^\times e^{\mathbb{R}I}e_+\ \bigsqcup\ \mathbb{R}^\times e^{\mathbb{R}I}e_-\ \bigsqcup\ \mathbb{R}^\times e^{\mathbb{R}I}I.
	\end{array}
	\end{displaymath}
	The rotations that are represented by rotors of this kind are called \emph{Lorentz boosts}.
	We will return to the Lorentz group in the 4-dimensional spacetime, but for now note the 
	hyperbolic nature of these rotations,
	\begin{equation}
	\setlength\arraycolsep{2pt}
	\begin{array}{rcl}
		e^{\alpha I} e_+ e^{-\alpha I} &=& e_+ e^{-2\alpha I} = e_+(\cosh 2\alpha - I \sinh 2\alpha) \\[5pt]
			&=& e_+ \cosh 2\alpha - e_- \sinh 2\alpha.
	\end{array}
	\end{equation}
	Hence, a rotor $\pm e^{\alpha I/2}$ transforms (or \emph{boosts}) timelike vectors 
	by a hyperbolic angle $\alpha$ away from the positive spacelike direction.
	
	The spinor space $\mathcal{G}^+$ consists partly of scaled rotors, $\rho e^{\alpha I}$,
	but there are also two subspaces $\rho(1 \pm I)$ of \emph{null spinors} which
	cannot be represented in this way. Note that such a spinor $\Psi$ acts on vectors as
	\begin{equation}
	\setlength\arraycolsep{2pt}
	\begin{array}{rcl}
		\Psi e_+ \Psi^\dagger &=& \rho^2 (1 \pm I)e_+(1 \mp I) = 2\rho^2(e_+ \mp e_-), \\[5pt]
		\Psi e_- \Psi^\dagger &=& \mp 2\rho^2(e_+ \mp e_-),
	\end{array}
	\end{equation}
	so it maps the whole space into one of the two null-spaces.
	The action of a non-null spinor has a nice interpretation as a boost plus scaling.

\subsection{$\mathcal{G}(\mathbb{R}^{3}) \stackrel{\sim}{=} \mathcal{G}(\mathbb{C}^{2})$ - The space algebra / Pauli algebra}
	
	Since the 3-dimensional euclidean space is the space that is most
	familiar to us humans, one could expect its geometric algebra, the
	\emph{space algebra}, to be familiar as well. Unfortunately, this is
	generally not the case. Most of its features, however, are commonly
	known but under different names and in separate contexts. For example,
	using the isomorphism $\mathcal{G}(\mathbb{R}^{3}) \cong \mathcal{G}(\mathbb{C}^{2}) \cong \mathbb{C}^{2 \times 2}$
	from Proposition \ref{prop_iso_real_complex}, we find that this algebra
	also appears in quantum mechanics in the form of the complex \emph{Pauli algebra}.
	
	We take an orthonormal basis $\{e_1,e_2,e_3\}$ in $\mathbb{R}^3$ and obtain
	\begin{equation} \label{space_algebra}
		\mathcal{G}(\mathbb{R}^3) = \textrm{Span}_\mathbb{R} \{1,\ e_1,e_2,e_3,\ e_1I,e_2I,e_3I,\ I\},
	\end{equation}
	where $I := e_1e_2e_3$ is the pseudoscalar. We write in this way to emphasize 
	the duality between the vector and bivector spaces in this case.
	This duality can be used to define cross products and rotation axes etc.
	However, since the use of such concepts is limited to three dimensions \emph{only},
	it is better to work with their natural counterparts within geometric algebra.
	
	To begin with, we would like to find the invertible elements of the space algebra.
	An arbitrary element $x \in \mathcal{G}$ can be written as
	\begin{equation} \label{arbitrary_space_element}
		x = \alpha + a + bI + \beta I,
	\end{equation}
	where $\alpha,\beta \in \mathbb{R}$ and $a,b \in \mathbb{R}^3$.
	Note that, since the algebra is odd, the pseudoscalar commutes with everything
	and furthermore $I^2 = -1$.
	The norm function $N_2$ does not take values in $\mathbb{R}$ in this algebra,
	but due to the properties of the Clifford conjugate we have 
	$N_2(x) = N_2(x)^\cliffconj \in \mathcal{G}^0 \oplus \mathcal{G}^3$.
	This subspace is, from our observation, isomorphic to 
	$\mathbb{C}$ and its corresponding complex
	conjugate is given by $[x]_3$. Using this, we can construct a real-valued map 
	$N_3 \!: \mathcal{G} \to \mathbb{R}^+$ by taking the complex modulus,
	\begin{equation} \label{N3_definition}
		N_3(x) := [N_2(x)]_3 N_2(x) = [x^\cliffconj x] x^\cliffconj x.
	\end{equation}
	Plugging in \eqref{arbitrary_space_element} we obtain
	\begin{equation} \label{N3_expression_space}
		N_3(x) = (\alpha^2 - a^2 + b^2 - \beta^2)^2 + 4(\alpha\beta - a*b)^2.
	\end{equation}
	Although $N_3$ takes values in $\mathbb{R}^+$, it is not a real norm\footnote{This 
	will also be seen to be required from dimensional considerations and the 
	remark to Hurwitz' Theorem in Section \ref{sec_reps}}
	on $\mathcal{G}$ since there are nonzero elements with $N_3(x)=0$.
	It does however have the multiplicative property
	\begin{equation}
	\setlength\arraycolsep{2pt}
	\begin{array}{crcl}
					& N_2(xy) &=& (xy)^\cliffconj xy = y^\cliffconj N_2(x) y = N_2(x)N_2(y) \\[5pt]
		\Rightarrow & N_3(xy) &=& [N_2(xy)]_3 N_2(xy) = [N_2(x)]_3 [N_2(y)]_3 N_2(x) N_2(y) \\[5pt]
					&	&=& N_3(x)N_3(y),
	\end{array}
	\end{equation}
	for all $x,y$, since $N_2(x)$ commutes with all of $\mathcal{G}$. 
	We also observe from \eqref{N3_expression_space} that $N_3(x^\cliffconj) = N_3(x)$.
	The expression \eqref{N3_definition} singles out the invertible elements
	as those elements \eqref{arbitrary_space_element} for which $N_3(x) \neq 0$, i.e.
	\begin{equation}
	\setlength\arraycolsep{2pt}
	\begin{array}{rcl}
		\mathcal{G}^\times(\mathbb{R}^{3}) &=& \{ x \in \mathcal{G} : N_3(x) \neq 0 \} \\[5pt]
			&=& \{ x \in \mathcal{G} : (\alpha^2 - a^2 + b^2 - \beta^2)^2 + 4(\alpha\beta - a*b)^2 \neq 0 \}.
	\end{array}
	\end{equation}
	
	The even subalgebra of $\mathcal{G}(\mathbb{R}^3)$ is the quaternion algebra, as
	follows from the isomorphism $\mathcal{G}^+(\mathbb{R}^{3,0}) \cong \mathcal{G}(\mathbb{R}^{0,2}) \cong \mathbb{H}$.
	The rotor group is then, according to the observation \eqref{lowerdim_spinor_rotor},
	the group of unit quaternions 
	(note that the reverse here acts as the quaternion conjugate),
	\begin{equation} \label{space_rotor_group}
		\textrm{Spin}^{(+)}(3,0) = \{ \alpha + bI \in \mathcal{G}^+ : \alpha^2 + b^2 = 1 \} 
			= e^{\mathcal{G}^2(\mathbb{R}^3)} \cong \textrm{SU}(2),
	\end{equation}
	where the exponentiation of the bivector algebra followed from Theorem \ref{thm_spin_is_exp}.
	The last isomorphism shows the relation to the Pauli algebra and is
	perhaps the most famous representation of a Spin group.
	An arbitrary rotor $\psi$ can according to \eqref{space_rotor_group} be written
	in the polar form $\psi = e^{\varphi\hat{n}I}$, where $\hat{n}$ is a unit vector,
	and represents a rotation by an angle $-2\varphi$ in the plane $\hat{n}I = -\hat{n}\dual$
	(i.e. $2\varphi$ counter-clockwise around the axis $\hat{n}$).
	
	The Pin group consists of two copies of the rotor group,
	\begin{equation}
	\setlength\arraycolsep{2pt}
	\begin{array}{rcl}
		\textrm{Pin}(3,0) &=& e^{\mathcal{G}^2(\mathbb{R}^3)}\ \bigsqcup\ e^{\mathcal{G}^2(\mathbb{R}^3)} \hat{n}, \\[5pt]
		\Gamma(3,0) &=& \mathbb{R}^\times e^{\mathcal{G}^2(\mathbb{R}^3)}\ \bigsqcup\ \mathbb{R}^\times e^{\mathcal{G}^2(\mathbb{R}^3)} \hat{n},
	\end{array}
	\end{equation}
	for some unit vector $\hat{n}$. The Pin group can be visualized as two
	unit 3-spheres $S^3$ lying in the even and odd subspaces, respectively.
	The odd one includes a reflection and represents the non-orientation-preserving part of O$(3)$.
	
	As we saw above, the spinors in the space algebra are quaternions.
	An arbitrary spinor can be written $\Psi = \rho e^{\varphi\hat{n}I/2}$
	and acts on vectors by rotating in the plane $\hat{n}\dual$ with the angle $\varphi$
	and scaling with $\rho^2$. We will continue our discussion on these
	spinors in Section \ref{sec_spinors}.
	
\subsection{$\mathcal{G}(\mathbb{R}^{1,3})$ - The spacetime algebra}
	
	We take as a four-dimensional example the \emph{spacetime algebra} (STA),
	which is the geometric algebra of Minkowski spacetime, $\mathbb{R}^{1,3}$.
	This is the stage for special relativistic physics and what is fascinating
	with the STA is that it embeds a lot of important physical objects in a natural way.
	
	By convention, we denote an orthonormal basis of the Minkowski space
	by $\{\gamma_0,\gamma_1,\gamma_2,\gamma_3\}$, where $\gamma_0$ is timelike
	and the other $\gamma_i$ are spacelike.
	This choice of notation is motivated by the Dirac representation of the STA
	in terms of so called gamma matrices which will be explained in more detail later.
	The STA expressed in this basis is
	\begin{displaymath} \label{spacetime_algebra}
	\setlength\arraycolsep{2pt}
	\begin{array}{l}
		\mathcal{G}(\mathbb{R}^{1,3}) = \\[3pt]
		\quad \textrm{Span}_\mathbb{R} \{1,\ 
			\gamma_0,\gamma_1,\gamma_2,\gamma_3,\ 
			\boldsymbol{e}_1,\boldsymbol{e}_2,\boldsymbol{e}_3, 
			\boldsymbol{e}_1I,\boldsymbol{e}_2I,\boldsymbol{e}_3I,\ 
			\gamma_0I,\gamma_1I,\gamma_2I,\gamma_3I,\ I\},
	\end{array}
	\end{displaymath}
	where the pseudoscalar is $I := \gamma_0\gamma_1\gamma_2\gamma_3$ 
	and we set $\boldsymbol{e}_i := \gamma_i\gamma_0$, $i=1,2,3$. The form of the STA basis
	chosen above emphasizes the duality which exists between the graded subspaces.
	It also hints that the even subalgebra of the STA is the space algebra \eqref{space_algebra}.
	This is true from the isomorphism $\mathcal{G}^+(\mathbb{R}^{1,3}) \cong \mathcal{G}(\mathbb{R}^{3,0})$,
	but we can also verify this explicitly by noting that 
	$\boldsymbol{e}_i^2 = 1$ and $\boldsymbol{e}_1\boldsymbol{e}_2\boldsymbol{e}_3 = I$.
	Hence, the timelike (positive square) blades $\{\boldsymbol{e}_i\}$ form a basis of
	a 3-dimensional euclidean space called the \emph{relative space} to $\gamma_0$.
	For any timelike vector $a$ we can find a similar relative space spanned by the bivectors
	$\{b \wedge a\}$ for $b \in \mathbb{R}^{1,3}$. These spaces all generate the
	\emph{relative space algebra} $\mathcal{G}^+$.
	This is a very powerful	concept which helps us visualize and work 
	efficiently in Minkowski spacetime and the STA.
	
	Using boldface to denote relative space elements, 
	an arbitrary multivector $x \in \mathcal{G}$ can be written
	\begin{equation} \label{arbitrary_STA_element}
		x = \alpha + a + \boldsymbol{a} + \boldsymbol{b}I + bI + \beta I,
	\end{equation}
	where $\alpha,\beta \in \mathbb{R}$, $a,b \in \mathbb{R}^{1,3}$ and $\boldsymbol{a},\boldsymbol{b}$ 
	in relative space $\mathbb{R}^3$.
	As usual, we would like to find the invertible elements.
	Looking at the norm function $N_2 \!: \mathcal{G} \to \mathcal{G}^0 \oplus \mathcal{G}^3 \oplus \mathcal{G}^4$,
	it is not obvious that we can extend this to a real-valued function on $\mathcal{G}$.
	Fortunately, we have for $X = \alpha + bI + \beta I \in \mathcal{G}^0 \oplus \mathcal{G}^3 \oplus \mathcal{G}^4$
	that
	\begin{equation} \label{G_0_3_4_norm}
		[X]_{3,4}X = (\alpha - bI - \beta I)(\alpha + bI + \beta I) = \alpha^2 - b^2 + \beta^2 \in \mathbb{R}.
	\end{equation}
	Hence, we can define a map $N_4 \!: \mathcal{G} \to \mathbb{R}$ by
	\begin{equation} \label{N4_definition}
		N_4(x) := [N_2(x)]_{3,4} N_2(x) = [x^\cliffconj x] x^\cliffconj x.
	\end{equation}
	Plugging in \eqref{arbitrary_STA_element} into $N_2$, we obtain after a tedious calculation
	\begin{equation} \label{N2_STA_expression}
	\setlength\arraycolsep{2pt}
	\begin{array}{rcl}
		N_2(x) &=& \alpha^2 - a^2 - \boldsymbol{a}^2 + \boldsymbol{b}^2 + b^2 - \beta^2 \\[5pt]
		&&	+\ 2(\alpha b - \beta a - a \liprod \boldsymbol{b} + b \liprod \boldsymbol{a} - a \liprod \boldsymbol{a}\dual - b \liprod \boldsymbol{b}\dual)I \\[5pt]
		&&	+\ 2(\alpha \beta - a*b - \boldsymbol{a}*\boldsymbol{b})I
	\end{array}
	\end{equation}
	and, by \eqref{G_0_3_4_norm},
	\begin{equation} \label{N4_expression}
	\setlength\arraycolsep{2pt}
	\begin{array}{rcl}
		N_4(x) &=& (\alpha^2 - a^2 - \boldsymbol{a}^2 + \boldsymbol{b}^2 + b^2 - \beta^2)^2 \\[5pt]
		&&	-\ 4(\alpha b - \beta a - a \liprod \boldsymbol{b} + b \liprod \boldsymbol{a} - a \liprod \boldsymbol{a}\dual - b \liprod \boldsymbol{b}\dual)^2 \\[5pt]
		&&	+\ 4(\alpha \beta - a*b - \boldsymbol{a}*\boldsymbol{b})^2.
	\end{array}
	\end{equation}
	We will prove some rather non-trivial statements about this norm function
	where we need that $[xy]x=x[yx]$ for all $x,y \in \mathcal{G}$.
	This is a quite general property of this involution.
	
	\begin{lem} \label{lem_nonscalar_conj}
		In any Clifford algebra $\cl(X,R,r)$ (even when $X$ is infinite), we have
		\begin{displaymath}
			[xy]x = x[yx] \quad \forall x,y \in \cl.
		\end{displaymath}
	\end{lem}
	\begin{proof}
		Using linearity, we can set $y=A \in \mathscr{P}(X)$ and expand $x$ in 
		coordinates $x_B \in R$ as $x = \sum_{B \in \mathscr{P}(X)} x_B B$.
		We obtain
		\begin{displaymath}
		\setlength\arraycolsep{2pt}
		\begin{array}{rcl}
			x[Ax]
				&=& \sum_{B,C} x_B x_C\ B[AC] \\[5pt]
				&=& \sum_{B,C} x_B x_C\ \big((A \symdiff C = \varnothing) - (A \symdiff C \neq \varnothing)\big)\ BAC \\[5pt]
				&=& \sum_{B,C} x_B x_C\ \big((A=C) - (A \neq C)\big)\ BAC \\[5pt]
				&=& \sum_B x_B x_A\ BAA - \sum_{C \neq A} \sum_B x_B x_C\ BAC
		\end{array}
		\end{displaymath}
		and
		\begin{displaymath}
		\setlength\arraycolsep{2pt}
		\begin{array}{rcl}
			[xA]x
				&=& \sum_{B,C} x_B x_C\ [BA]C \\[5pt]
				&=& \sum_{B,C} x_B x_C\ \big((B=A) - (B \neq A)\big)\ BAC \\[5pt]
				&=& \sum_C x_A x_C\ AAC - \sum_{B \neq A} \sum_C x_B x_C\ BAC \\[5pt]
				&=& x_A^2\ AAA + \sum_{C \neq A} x_A x_C\ AAC - \sum_{B \neq A} x_B x_A \underbrace{BAA}_{AAB} \\[5pt]
				&& \quad - \sum_{B \neq A} \sum_{C \neq A} x_B x_C\ BAC \\[5pt]
				&=& x_A^2\ AAA - \sum_{B \neq A} \sum_{C \neq A} x_B x_C\ BAC \\[5pt]
				&=& x[Ax].
		\end{array}
		\end{displaymath}
	\end{proof}
	
	\noindent
	We now have the following
	
	\begin{thm} \label{thm_norm_conj_4}
		$N_4(x^\cliffconj) = N_4(x)$ for all $x \in \mathcal{G}(\mathbb{R}^{1,3})$.
	\end{thm}
	
	\begin{rem}
		Note that this is not at all obvious from the expression \eqref{N4_expression}.
	\end{rem}
	
	\begin{proof}
		Using Lemma \ref{lem_nonscalar_conj} we have that
		\begin{equation}
			N_4(x^\cliffconj) = [xx^\cliffconj]xx^\cliffconj = x[x^\cliffconj x]x^\cliffconj.
		\end{equation}
		Since $N_4$ takes values in $\mathbb{R}$, this must be a scalar, so that
		\begin{equation}
			N_4(x^\cliffconj) = \langle x[x^\cliffconj x]x^\cliffconj \rangle_0 
				= \langle [x^\cliffconj x]x^\cliffconj x \rangle_0 = \langle N_4(x) \rangle_0 = N_4(x),
		\end{equation}
		where we used the symmetry of the scalar product.
	\end{proof}
	
	\begin{lem} \label{lem_conj_034}
		For all $X,Y \in \mathcal{G}^0 \oplus \mathcal{G}^3 \oplus \mathcal{G}^4$ we have
		\begin{displaymath}
			[XY] = [Y][X].
		\end{displaymath}
	\end{lem}
	\begin{proof}
		Take arbitrary elements $X = \alpha + bI + \beta I$ and $Y = \alpha' + b'I + \beta' I$. Then
		\begin{displaymath}
		\setlength\arraycolsep{2pt}
		\begin{array}{rcl}
			[XY]
				&=& [(\alpha + bI + \beta I)(\alpha' + b'I + \beta' I)]\\[5pt]
				&=& \alpha \alpha' - \alpha b'I - \alpha \beta' I - bI\alpha' + b*b' - b \wedge b' + b \beta' - \beta I \alpha' - \beta b' - \beta \beta'
		\end{array}
		\end{displaymath}
		and
		\begin{displaymath}
		\setlength\arraycolsep{2pt}
		\begin{array}{rcl}
			[Y][X]
				&=& (\alpha' - b'I - \beta' I)(\alpha - bI - \beta I)\\[5pt]
				&=& \alpha' \alpha - \alpha' bI - \alpha' \beta I - b'I\alpha + b'*b + b' \wedge b - b' \beta - \beta' I \alpha + \beta' b - \beta' \beta.
		\end{array}
		\end{displaymath}
		Comparing these expressions we find that they are equal.
	\end{proof}
	
	\noindent
	We can now prove that $N_4$ really acts as a determinant on the STA.
	
	\begin{thm} \label{thm_norm_prod_4}
		The norm function $N_4$ satisfies the product property
		\begin{displaymath}
			N_4(xy) = N_4(x) N_4(y) \quad \forall x,y \in \mathcal{G}(\mathbb{R}^{1,3}).
		\end{displaymath}
	\end{thm}
	
	\begin{proof}
		Using that $N_4(xy)$ is a scalar and that $N_2$ takes values in 
		$\mathcal{G}^0 \oplus \mathcal{G}^3 \oplus \mathcal{G}^4$, we obtain
		\begin{displaymath}
		\setlength\arraycolsep{2pt}
		\begin{array}{rcl}
			N_4(xy) 
				&=& [(xy)^\cliffconj xy](xy)^\cliffconj xy 
				= [y^\cliffconj x^\cliffconj x y] y^\cliffconj x^\cliffconj x y \\[5pt]
				&=& \langle [y^\cliffconj x^\cliffconj x y] y^\cliffconj x^\cliffconj x y \rangle_0
				= \langle x^\cliffconj x y [y^\cliffconj x^\cliffconj x y] y^\cliffconj \rangle_0 \\[5pt]
				&=& \langle x^\cliffconj x [y y^\cliffconj x^\cliffconj x] y y^\cliffconj \rangle_0
				= \langle N_2(x) [N_2(y^\cliffconj) N_2(x)] N_2(y^\cliffconj) \rangle_0 \\[5pt]
				&=& \langle N_2(x) [N_2(x)] [N_2(y^\cliffconj)] N_2(y^\cliffconj) \rangle_0
				= \langle N_4(x) N_4(y^\cliffconj) \rangle_0 \\[5pt]
				&=& N_4(x) N_4(y^\cliffconj),
		\end{array}
		\end{displaymath}
		where we applied Lemma \ref{lem_nonscalar_conj} and then Lemma \ref{lem_conj_034}.
		Theorem \ref{thm_norm_conj_4} now gives the claimed identity.
	\end{proof}

	From \eqref{N4_definition} we find that the group of invertible elements is given by
	\begin{equation}
		\mathcal{G}^\times(\mathbb{R}^{1,3}) = \{ x \in \mathcal{G} : N_4(x) \neq 0 \}
	\end{equation}
	and the inverse of $x \in \mathcal{G}^\times$ is
	\begin{equation}
		x^{-1} = \frac{1}{N_4(x)} [x^\cliffconj x] x^\cliffconj.
	\end{equation}
	Note that the above theorems regarding $N_4$ only rely on the commutation properties
	of the different graded subspaces and not on the actual signature
	and field of the vector space. Therefore, these hold for all $\cl(X,R,r)$
	such that $|X|=4$, and 
	\begin{equation}
		\cl^\times = \{ x \in \cl : [x^\cliffconj x] x^\cliffconj x \in R\ \textrm{is invertible} \}.
	\end{equation}
	
	Let us now turn our attention to the rotor group of the STA.
	The reverse equals the Clifford conjugate on the even subalgebra (it also corresponds
	to the Clifford conjugate defined on the relative space),
	so we find from \eqref{N2_STA_expression} that the rotor group is
	\begin{equation} \label{STA_rotor_group}
	\setlength\arraycolsep{2pt}
	\begin{array}{l}
		\textrm{Spin}^{+}(1,3) 
				= \{ x \in \mathcal{G}^+ : N_2(x) = x^\cliffconj x = 1  \} \\[5pt]
		\quad	= \{ \alpha + \boldsymbol{a} + \boldsymbol{b}I + \beta I \in \mathcal{G}^+ : \alpha^2 - \boldsymbol{a}^2 + \boldsymbol{b}^2 - \beta^2 = 1,\ \alpha\beta = \boldsymbol{a}*\boldsymbol{b} \} \\[5pt]
		\quad	= \pm e^{\mathcal{G}^2(\mathbb{R}^{1,3})} \cong \textrm{SL}(2,\mathbb{C}).
	\end{array}
	\end{equation}
	The last isomorphism is related to the Dirac representation of the STA, while
	the exponentiation identity was obtained from Theorem \ref{thm_spin_is_exp} and gives a better picture
	of what the rotor group looks like. Namely, any rotor $\psi$ can be written $\psi = \pm e^{\boldsymbol{a}+\boldsymbol{b}I}$
	for some relative vectors $\boldsymbol{a},\boldsymbol{b}$.
	A pure $\pm e^{\boldsymbol{b}I}$ corresponds to a rotation in the spacelike plane
	$\boldsymbol{b}\dual$ with angle $2|\boldsymbol{b}|$ 
	(which is a corresponding rotation also in relative space),
	while $\pm e^{\boldsymbol{a}}$ corresponds to a ``rotation" in the timelike plane $\boldsymbol{a}$,
	i.e. a boost in the relative space direction $\boldsymbol{a}$ with 
	velocity $\textrm{arctanh}(2|\boldsymbol{a}|)$ times the speed of light.
	
	Picking for example $\gamma_0$ and $\gamma_1$, we obtain the Spin and Pin groups,
	\begin{displaymath}
	\setlength\arraycolsep{2pt}
	\begin{array}{rcl}
		\textrm{Spin}(1,3) &=& \pm e^{\mathcal{G}^2(\mathbb{R}^{1,3})}\ \bigsqcup\ \pm e^{\mathcal{G}^2(\mathbb{R}^{1,3})} \gamma_0\gamma_1, \\[5pt]
		\textrm{Pin}(1,3) &=& \pm e^{\mathcal{G}^2(\mathbb{R}^{1,3})}\ \bigsqcup\ \pm e^{\mathcal{G}^2(\mathbb{R}^{1,3})} \gamma_0\ \bigsqcup\ \pm e^{\mathcal{G}^2(\mathbb{R}^{1,3})} \gamma_1\ \bigsqcup\ \pm e^{\mathcal{G}^2(\mathbb{R}^{1,3})} \gamma_0\gamma_1, \\[5pt]
		\Gamma(1,3) &=& \mathbb{R}^\times e^{\mathcal{G}^2(\mathbb{R}^{1,3})} \bigsqcup \mathbb{R}^\times e^{\mathcal{G}^2(\mathbb{R}^{1,3})} \gamma_0 \bigsqcup \mathbb{R}^\times e^{\mathcal{G}^2(\mathbb{R}^{1,3})} \gamma_1 \bigsqcup \mathbb{R}^\times e^{\mathcal{G}^2(\mathbb{R}^{1,3})} \gamma_0\gamma_1.
	\end{array}
	\end{displaymath}
	The Pin group forms a double-cover of the so called \emph{Lorentz group} O(1,3).
	Since the rotor group is connected, we find that O(1,3) has four connected components.
	The Spin group covers the subgroup of \emph{proper} Lorentz transformations preserving orientation,
	while the rotor group covers the connected \emph{proper orthochronous} Lorentz group
	which also preserves the direction of time.
	
	The spinor space of the STA is the relative space algebra.
	We will discuss these spinors in more detail later, but for now note that
	an invertible spinor 
	\begin{equation}
		\Psi \Psi^\dagger = \rho e^{I\varphi} \in \mathcal{G}^0 \oplus \mathcal{G}^4
		\quad \Rightarrow \quad
		\Psi = \rho^{1/2} e^{I\varphi/2} \psi
	\end{equation}
	is the product of a rotor $\psi$, a \emph{duality rotor} $e^{I\varphi/2}$ 
	and a scale factor $\rho^{1/2}$.
	
	
\subsection{$\mathcal{G}(\mathbb{R}^{4,1}) \stackrel{\sim}{=} \mathcal{G}(\mathbb{C}^{4})$ - The Dirac algebra}
	
	The Dirac algebra is the representation of the STA which is most
	commonly used in physics. This is due to historic reasons, since the geometric
	nature of this algebra from its relation to the spacetime algebra was not uncovered
	until the 1960s. 
	The relation between these algebras is observed by noting that the
	pseudoscalar in $\mathcal{G}(\mathbb{R}^{4,1})$ commutes with all
	elements and squares to minus the identity. By Proposition \ref{prop_iso_real_complex}
	we have that the Dirac algebra is the complexification of the STA,
	\begin{equation}
		\mathcal{G}(\mathbb{R}^{4,1}) \cong \mathcal{G}(\mathbb{R}^{1,3}) \otimes \mathbb{C} \cong \mathcal{G}(\mathbb{C}^4) \cong \mathbb{C}^{4 \times 4}.
	\end{equation}
	We construct this isomorphism explicitly by taking bases
	$\{\gamma_0,\gamma_1,\gamma_2,\gamma_3\}$ of $\mathbb{R}^{1,3}$ as usual
	and $\{e_0,\ldots,e_4\}$ of $\mathbb{R}^{4,1}$ such that $e_0^2 = -1$ and the other $e_j^2 = 1$.
	We write $\mathcal{G}_5 := \mathcal{G}(\mathbb{R}^{4,1})$ and 
	$\mathcal{G}_4^\mathbb{C} := \mathcal{G}(\mathbb{R}^{1,3}) \otimes \mathbb{C}$, 
	and use the convention that Greek indices run from 0 to 3.
	The isomorphism $F \!: \mathcal{G}_5 \to \mathcal{G}_4^\mathbb{C}$
	is given by the following correspondence of basis elements:
	\begin{displaymath}
	\begin{array}{rccccccc}
		\mathcal{G}_4^\mathbb{C}:
		& 1 \otimes 1
		& \gamma_\mu \otimes 1
		& \gamma_\mu \wedge \gamma_\nu \otimes 1
		& \gamma_\mu \wedge \gamma_\nu \wedge \gamma_\lambda \otimes 1
		& I_4 \otimes 1
		& 1 \otimes i \\[5pt]
		\mathcal{G}_5:
		& 1
		& e_\mu e_4
		& -e_\mu \wedge e_\nu
		& -e_\mu \wedge e_\nu \wedge e_\lambda e_4
		& e_0 e_1 e_2 e_3
		& I_5 \\[5pt]
		x^\cliffconj\ \textrm{in}\ \mathcal{G}_4^\mathbb{C}:
		& + & - & - & + & + & + \\[5pt]
		[x]\ \textrm{in}\ \mathcal{G}_4^\mathbb{C}:
		& + & - & - & - & - & + \\[5pt]
		\overline{x}\ \textrm{in}\ \mathcal{G}_4^\mathbb{C}:
		& + & + & + & + & + & - \\[5pt]
	\end{array}
	\end{displaymath}
	The respective pseudoscalars are $I_4 := \gamma_0\gamma_1\gamma_2\gamma_3$
	and $I_5 := e_0e_1e_2e_3e_4$.
	We have also noted the correspondence between involutions in
	the different algebras. Clifford conjugate in $\mathcal{G}_4^\mathbb{C}$
	corresponds to reversion in $\mathcal{G}_5$, 
	the $[\ ]$-involution becomes the $[\ ]_{1,2,3,4}$-involution,
	while complex conjugation in $\mathcal{G}_4^\mathbb{C}$ corresponds 
	to grade involution in $\mathcal{G}_5$. In other words,
	\begin{equation}
		F(x^\dagger) = F(x)^\cliffconj,
		\quad F([x]_{1,2,3,4}) = [F(x)],
		\quad F(x^\star) = \overline{F(x)}.
	\end{equation}
	
	We can use the correspondence above to find a norm function on $\mathcal{G}_5$.
	Since $N_4 \!: \mathcal{G}(\mathbb{R}^{1,3}) \to \mathbb{R}$ was independent of the choice of field,
	we have that the complexification of $N_4$ satisfies
	\begin{displaymath}
	\setlength\arraycolsep{2pt}
	\begin{array}{rcl}
		N_4^\mathbb{C} \!: \mathcal{G}(\mathbb{C}^{4}) &\to& \mathbb{C}, \\[5pt]
			x &\mapsto& [x^\cliffconj x] x^\cliffconj x.
	\end{array}
	\end{displaymath}
	Taking the modulus of this complex number, we arrive at a real-valued map
	$N_5 \!: \mathcal{G}(\mathbb{R}^{4,1}) \to \mathbb{R}$ with
	\begin{displaymath}
	\setlength\arraycolsep{2pt}
	\begin{array}{rrl}
		N_5(x) &:=& \overline{N^{\mathbb{C}}_4\big(F(x)\big)} N^{\mathbb{C}}_4\big(F(x)\big) \\[5pt]
			&=& \overline{ [F(x)^\cliffconj F(x)] F(x)^\cliffconj F(x) }\ [F(x)^\cliffconj F(x)] F(x)^\cliffconj F(x) \\[5pt]
			&=& \big[ [x^\dagger x]_{1,2,3,4} x^\dagger x \big]_5 [x^\dagger x]_{1,2,3,4} x^\dagger x \\[5pt]
			&=& \big[ [x^\dagger x]_{1,4} x^\dagger x \big] [x^\dagger x]_{1,4} x^\dagger x.
	\end{array}
	\end{displaymath}
	In the final steps we noted that $x^\dagger x \in \mathcal{G}^0 \oplus \mathcal{G}^1 \oplus \mathcal{G}^4 \oplus \mathcal{G}^5$
	and that $\mathbb{C} \subseteq \mathcal{G}_4^\mathbb{C}$ corresponds
	to $\mathcal{G}^0 \oplus \mathcal{G}^5 \subseteq \mathcal{G}_5$.
	Furthermore, since $N^{\mathbb{C}}_4(xy) = N^{\mathbb{C}}_4(x) N^{\mathbb{C}}_4(y)$, we have
	\begin{displaymath}
	\setlength\arraycolsep{2pt}
	\begin{array}{rrl}
		N_5(xy) 
			&=& \overline{N^{\mathbb{C}}_4\big(F(x)F(y)\big)}\ N^{\mathbb{C}}_4\big(F(x)F(y)\big) \\[5pt]
			&=& \overline{N^{\mathbb{C}}_4\big(F(x)\big)}\ \overline{N^{\mathbb{C}}_4\big(F(y)\big)}\ N^{\mathbb{C}}_4\big(F(x)\big)\ N^{\mathbb{C}}_4\big(F(y)\big) \\[5pt]
			&=& N_5(x) N_5(y)
	\end{array}
	\end{displaymath}
	for all $x,y \in \mathcal{G}$.
	The invertible elements of the Dirac algebra are then as usual
	\begin{equation}
		\mathcal{G}^\times(\mathbb{R}^{4,1}) = \{ x \in \mathcal{G} : N_5(x) \neq 0 \}
	\end{equation}
	and the inverse of $x \in \mathcal{G}^\times$ is
	\begin{equation}
		x^{-1} = \frac{1}{N_5(x)} \big[ [x^\dagger x]_{1,4} x^\dagger x \big] [x^\dagger x]_{1,4} x^\dagger.
	\end{equation}
	The above strategy could also have been used to obtain the expected 
	result for $N_3$ on $\mathcal{G}(\mathbb{R}^{3,0}) \cong \mathcal{G}(\mathbb{C}^2)$
	(with a corresponding isomorphism $F$):
	\begin{equation}
		N_3(x) := \overline{N^{\mathbb{C}}_2\big(F(x)\big)} N^{\mathbb{C}}_2\big(F(x)\big) = [x^\cliffconj x] x^\cliffconj x.
	\end{equation}
	
	We briefly describe how spinors are dealt with in this representation.
	This will not be the the same as the even subspace spinors which we
	usually consider. For the selected basis of $\mathcal{G}_4^\mathbb{C}$ we form
	the idempotent element
	\begin{equation}
		\textstyle
		f := \frac{1}{2}(1 + \gamma_0)\frac{1}{2}(1 + i\gamma_1\gamma_2).
	\end{equation}
	Spinors are now defined as elements of the ideal $\mathcal{G}_4^\mathbb{C} f$
	and one can show that every such element can be written as
	\begin{equation}
		\textstyle
		\Psi = \sum_{i=1}^4 \psi_i f_i, \quad \psi_i \in \mathbb{C},
	\end{equation}
	where
	\begin{equation}
		f_1 := f, \quad f_2 := -\gamma_1\gamma_3f, \quad f_3 := \gamma_3\gamma_0f, \quad f_4 := \gamma_1\gamma_0f.
	\end{equation}
	With the standard representation of $\gamma_\mu$ as generators of $\mathbb{C}^{4 \times 4}$,
	\begin{equation} \label{gamma_matrices}
		\gamma_0 = \left[ \begin{array}{cc} I & 0 \\ 0 & -I \end{array} \right], \quad
		\gamma_i = \left[ \begin{array}{cc} 0 & -\sigma_i \\ \sigma_i & 0 \end{array} \right],\ i=1,2,3,
	\end{equation}
	where $\sigma_i$ are the \emph{Pauli matrices} which generate a representation
	of the Pauli algebra, one finds that
	\begin{equation}
		\Psi = \left[ \begin{array}{cccc}
			\psi_1 &0&0&0 \\ 
			\psi_2 &0&0&0 \\ 
			\psi_3 &0&0&0 \\ 
			\psi_4 &0&0&0 \\ 
			\end{array} \right].
	\end{equation}
	Hence, these spinors can be thought of as complex column vectors.
	Upon a transformation of the basis $\{\gamma_\mu\}$, the components $\psi_i$ will
	transform according to the representation $D^{(1/2,0)} \oplus D^{(0,1/2)}$ of $\textrm{SL}(2,\mathbb{C})$.
	We will come back to discuss these spinors in Section \ref{sec_spinors}.
	See \cite{rodrigues_et_al} for more details on this correspondence between
	spinors and ideals.
	
\subsection{Summary of norm functions}
	
	
	The norm functions 
	\begin{displaymath}
	\setlength\arraycolsep{2pt}
	\begin{array}{rcl}
		N_0(x) &:=& x, \\[5pt]
		N_1(x) &=& x^\cliffconj x, \\[5pt]
		N_2(x) &=& x^\cliffconj x, \\[5pt]
		N_3(x) &=& [x^\cliffconj x] x^\cliffconj x, \\[5pt]
		N_4(x) &=& [x^\cliffconj x] x^\cliffconj x, \\[5pt]
		N_5(x) &=& \big[ [x^\dagger x]_{1,4} x^\dagger x \big] [x^\dagger x]_{1,4} x^\dagger x
	\end{array}
	\end{displaymath}
	constructed above (where we added $N_0$ for completeness)
	all have the product property
	\begin{equation} \label{norm_function_property}
		N_k(xy) = N_k(x) N_k(y)
	\end{equation}
	for all $x,y \in \mathcal{G}(\mathcal{V})$ when $\mathcal{V}$ is $k$-dimensional.
	Because these functions only involve products and involutions,
	and the proofs of the above identities only rely on commutation properties in the
	respective algebras, they even hold for any Clifford algebra $\cl(X,R,r)$
	with $|X|=0,1,\ldots,5$, respectively.
	
	For matrix algebras, a similar product property is satisfied by the determinant.
	On the other hand, we have the following theorem for matrices.
	
	\begin{thm} \label{thm_determinant}
		Assume that $d\!: \mathbb{R}^{n \times n} \to \mathbb{R}$ is continuous
		and satisfies
		\begin{equation} \label{determinant_property}
			d(AB) = d(A)d(B)
		\end{equation}
		for all $A,B \in \mathbb{R}^{n \times n}$.
		Then $d$ must be either $0$, $1$, $|\det|^\alpha$ or $(\textrm{\emph{sign}} \circ \det) |\det|^\alpha$
		for some $\alpha > 0$.
	\end{thm}
	
	\noindent
	In other words, we must have that $d = d_1 \circ \det$, where $d_1\!: \mathbb{R} \to \mathbb{R}$
	is continuous and $d_1(\lambda \mu) = d_1(\lambda) d_1(\mu)$.
	This $d_1$ is uniquely determined e.g. by whether $d$ takes negative values, together with
	the value of $d(\lambda I)$ for any $\lambda > 1$.
	This means that the determinant is the $unique$ real-valued function on real matrices with
	the product property \eqref{determinant_property}.
	The proof of this theorem can be found in the appendix.
	
	Looking at Table \ref{table_classification_real}, we see that 
	$\mathcal{G}(\mathbb{R}^{k,k}) \cong \mathbb{R}^{2^k \times 2^k}$ for $k=0,1,2,\ldots$
	From the above theorem we then know that there are \emph{unique}\footnote{Actually,
	the functions are either $\det$ or $|\det|$. $N_2$ and $N_4$ constructed previously
	are smooth, however, so they must be equal to $\det$.}
	continuous functions $N_{2k}\!: \mathcal{G}(\mathbb{R}^{k,k}) \to \mathbb{R}$
	such that $N_{2k}(xy) = N_{2k}(x) N_{2k}(y)$ and $N_{2k}(\lambda) = \lambda^{2^k}$.
	These are given by the determinant on the corresponding matrix algebra.
	What we do not know, however, is if every one of these can be expressed in the same simple
	form as $N_0$, $N_2$ and $N_4$, i.e. as a composition of products and grade-based involutions.
	Due to the complexity of higher-dimensional algebras, it is not obvious 
	whether a continuation of the strategy employed so far can be successful or not.
	It is even difficult\footnote{The first couple of $N_k$ can be verified directly
	using a geometric algebra package in Maple, but already for $N_4$ this
	becomes impossible to do straight-away on a standard desktop computer.}
	to test out suggestions of norm functions on a computer,
	since the number of operations involved grows as $2^{2k \cdot 2^k}$.
	We therefore leave this question as a suggestion for further investigation.
	
	Because of the product property \eqref{norm_function_property},
	the norm functions also lead to interesting factorization
	identities on rings. An example is $N_2$ for quaternions,
	\begin{equation} \label{lagrange_identity}
	\setlength\arraycolsep{2pt}
	\begin{array}{l}
		(x_1^2 + x_2^2 + x_3^2 + x_4^2)(y_1^2 + y_2^2 + y_3^2 + y_4^2) \\[5pt]
		\quad	= (x_1y_1 - x_2y_2 - x_3y_3 - x_4y_4)^2
				+ (x_1y_2 + x_2y_1 + x_3y_4 - x_4y_3)^2 \\[5pt]
		\qquad +\ (x_1y_3 - x_2y_4 + x_3y_1 + x_4y_2)^2
				+ (x_1y_4 + x_2y_3 - x_3y_2 + x_4y_1)^2.
	\end{array}
	\end{equation}
	This is called the \emph{Lagrange identity}.
	These types of identities can be used to prove theorems in number
	theory. Using \eqref{lagrange_identity}, one can for example prove that every integer
	can be written as a sum of four squares of integers. Or, in other
	words, every integer is the norm (squared) of an integral quaternion.
	See e.g. \cite{herstein} for the proof.
	
	Another possible application of norm functions could be in public key
	cryptography and one-way trapdoor functions. We have not investigated
	this idea further, however.

	\newpage

\section{Representation theory} \label{sec_reps}

	In this section we will use the classification of geometric algebras as
	matrix algebras, which was developed in Section \ref{sec_isomorphisms},
	to work out the representation theory of these algebras.
	Since one can find representations of geometric algebras in many
	areas of mathematics and physics, this leads to a number of interesting
	applications.
	We will consider two main examples in detail, namely normed division algebras
	and vector fields on higher-dimensional spheres.
	
	\begin{defn} \label{def_representation}
		For $\mathbb{K} = \mathbb{R}$, $\mathbb{C}$ or $\mathbb{H}$, we define a 
		\emph{$\mathbb{K}$-representation} of $\mathcal{G}(\mathcal{V},q)$
		as an $\mathbb{R}$-algebra homomorphism
		\begin{displaymath}
			\rho\!: \mathcal{G}(\mathcal{V},q) \to \textrm{End}_\mathbb{K}(W),
		\end{displaymath}
		where $W$ is a finite-dimensional vector space over $\mathbb{K}$.
		$W$ is called a \emph{$\mathcal{G}(\mathcal{V},q)$-module} over $\mathbb{K}$.
	\end{defn}
	
	\noindent
	Note that a vector space over $\mathbb{C}$ or $\mathbb{H}$ can be considered as a real vector space together
	with operators $J$ or $I,J,K$ in $\textrm{End}_\mathbb{R}(W)$ that anticommute and 
	square to minus the identity.
	In the definition above we assume that these operators commute with $\rho(x)$
	for all $x \in \mathcal{G}$, so that $\rho$ can be said to respect the $\mathbb{K}$-structure
	of the space $W$. When talking about the dimension of the module $W$
	we will always refer to its dimension as a real vector space.
	
	The standard strategy when studying representation theory is to look
	for irreducible representations.
	
	\begin{defn} \label{def_reducible_rep}
		A representation $\rho$ is called \emph{reducible} if $W$ can be written
		as a sum of proper (not equal to $0$ or $W$) invariant subspaces, i.e.
		\begin{displaymath}
			W = W_1 \oplus W_2 \quad \textrm{and} \quad \rho(x)(W_j) \subseteq W_j \quad \forall\ x \in \mathcal{G}.
		\end{displaymath}
		In this case we can write $\rho = \rho_1 \oplus \rho_2$, where $\rho_j(x) := \rho(x)|_{W_j}$.
		A representation is called \emph{irreducible} if it is not reducible.
	\end{defn}
	
	\noindent
	The traditional definition of an irreducible representation is that it
	does not have any proper invariant subspaces. However, because $\mathcal{G}$ is
	generated by a finite group (the Clifford group) one can verify that these 
	two definitions are equivalent in this case.
	
	\begin{prop} \label{prop_irreps}
		Every $\mathbb{K}$-representation $\rho$ of a geometric algebra $\mathcal{G}(\mathcal{V},q)$ 
		can be split up into a direct sum $\rho = \rho_1 \oplus \ldots \oplus \rho_m$
		of irreducible representations.
	\end{prop}
	\begin{proof}
		This follows directly from the definitions and the fact that $W$ is finite-dimensional.
	\end{proof}
	
	\begin{defn} \label{def_equivalent_reps}
		Two $\mathbb{K}$-representations $\rho_j\!: \mathcal{G}(\mathcal{V},q) \to \textrm{End}_\mathbb{K}(W_j)$,
		$j=1,2$, are said to be \emph{equivalent} if there exists a $\mathbb{K}$-linear 
		isomorphism $F\!: W_1 \to W_2$ such that
		\begin{displaymath}
			F \circ \rho_1(x) \circ F^{-1} = \rho_2(x) \quad \forall\ x \in \mathcal{G}.
		\end{displaymath}
	\end{defn}
	
	\begin{thm} \label{thm_matrix_reps}
		Up to equivalence, the only irreducible representations of the matrix
		algebras $\mathbb{K}^{n \times n}$ and $\mathbb{K}^{n \times n} \oplus \mathbb{K}^{n \times n}$ are
		\begin{displaymath}
			\rho\!: \mathbb{K}^{n \times n} \to \textrm{\emph{End}}_\mathbb{K}(\mathbb{K}^n)
		\end{displaymath}
		and
		\begin{displaymath}
			\rho_{1,2}\!: \mathbb{K}^{n \times n} \oplus \mathbb{K}^{n \times n} \to \textrm{\emph{End}}_\mathbb{K}(\mathbb{K}^n)
		\end{displaymath}
		respectively, where $\rho$ is the defining representation and
		\begin{displaymath}
		\setlength\arraycolsep{2pt}
		\begin{array}{c}
			\rho_1(x,y):=\rho(x), \\
			\rho_2(x,y):=\rho(y).
		\end{array}
		\end{displaymath}
	\end{thm}
	\begin{proof}
		This follows from the classical fact that the algebras 
		$\mathbb{K}^{n \times n}$ are simple and that simple algebras have only
		one irreducible representation up to equivalence. See e.g. \cite{lang} for details.
	\end{proof}
	
	\begin{thm} \label{thm_representations}
		From the above, together with the classification of real geometric algebras,
		follows the table of representations in Table \ref{table_representations},
		where $\nu_{s,t}$ is the number of inequivalent irreducible representations
		and $d_{s,t}$ is the dimension of an irreducible representation for $\mathcal{G}(\mathbb{R}^{s,t})$.
		The cases for $n>8$ are obtained using the periodicity
		\begin{equation} \label{rep_periodicity}
		\setlength\arraycolsep{2pt}
		\begin{array}{lcl}
			\nu_{m+8k} &=& \nu_m, \\[5pt]
			d_{m+8k} &=& 16^k d_m.
		\end{array}
		\end{equation}
	\end{thm}
	
	\begin{table}[ht]
		\begin{displaymath}
		\begin{array}{c|l|c|c|l|c|c}
				n & \mathcal{G}(\mathbb{R}^{n,0}) & \nu_{n,0} & d_{n,0} & \mathcal{G}(\mathbb{R}^{0,n}) & \nu_{0,n} & d_{0,n} \\
			\hline
			&&&&&&\\[-1.8ex]
			0 & \mathbb{R}												& 1 & 1 	& \mathbb{R}												& 1 & 1 \\
			1 & \mathbb{R} \oplus \mathbb{R}							& 2 & 1 	& \mathbb{C}												& 1 & 2 \\
			2 & \mathbb{R}^{2 \times 2}									& 1 & 2 	& \mathbb{H}												& 1 & 4 \\
			3 & \mathbb{C}^{2 \times 2}									& 1 & 4 	& \mathbb{H} \oplus \mathbb{H}								& 2 & 4 \\
			4 & \mathbb{H}^{2 \times 2}									& 1 & 8 	& \mathbb{H}^{2 \times 2}									& 1 & 8 \\
			5 & \mathbb{H}^{2 \times 2} \oplus \mathbb{H}^{2 \times 2}	& 2 & 8 	& \mathbb{C}^{4 \times 4}									& 1 & 8 \\
			6 & \mathbb{H}^{4 \times 4}									& 1 & 16 	& \mathbb{R}^{8 \times 8}									& 1 & 8 \\
			7 & \mathbb{C}^{8 \times 8}									& 1 & 16 	& \mathbb{R}^{8 \times 8} \oplus \mathbb{R}^{8 \times 8}	& 2 & 8 \\
			8 & \mathbb{R}^{16 \times 16}								& 1 & 16 	& \mathbb{R}^{16 \times 16}									& 1 & 16
		\end{array}
		\end{displaymath}
		\caption{Number and dimension of irreducible representations
		of euclidean and anti-euclidean geometric algebras. \label{table_representations}}
	\end{table}
	
	We will now consider the situation when the representation space $W$
	is endowed with an inner product.
	Note that if $W$ is a vector space over $\mathbb{K}$ with an inner product,
	we can always find a \emph{$\mathbb{K}$-invariant} inner product on $W$, i.e.
	such that the operators $J$ or $I,J,K$ are orthogonal.
	Namely, let $\langle\cdot,\cdot\rangle_\mathbb{R}$ be an inner product on $W$
	and put
	\begin{equation}
		\langle x,y \rangle_\mathbb{C} := \sum_{\Gamma \in \{id,J\}} \langle \Gamma x, \Gamma y \rangle_\mathbb{R}, \quad
		\langle x,y \rangle_\mathbb{H} := \sum_{\Gamma \in \{id,I,J,K\}} \langle \Gamma x, \Gamma y \rangle_\mathbb{R}.
	\end{equation}
	Then $\langle Jx,Jy \rangle_\mathbb{K} = \langle x,y \rangle_\mathbb{K}$ 
	and $\langle Jx,y \rangle_\mathbb{K} = -\langle x,Jy \rangle_\mathbb{K}$, etc.
	
	In the same way, when $\mathcal{V}$ is euclidean or anti-euclidean, we can 
	for a representation $\rho\!: \mathcal{G}(\mathcal{V}) \to \textrm{End}_\mathbb{K}(W)$
	find an inner product such that $\rho$ acts orthogonally with unit vectors, 
	i.e. such that $\langle \rho(e)x,\rho(e)y \rangle = \langle x,y \rangle$ for
	all $x,y \in W$ and $e \in \mathcal{V}$ with $e^2 = \pm 1$.
	We construct such an inner product by averaging a, possibly $\mathbb{K}$-invariant,
	inner product $\langle \cdot,\cdot \rangle_\mathbb{K}$ over the Clifford group.
	Take an orthonormal basis $E$ of $\mathcal{V}$ and put
	\begin{equation} \label{rho_inv_inner_prod}
		\langle x,y \rangle := \sum_{\Gamma \in \mathcal{B}_E} \langle \rho(\Gamma) x, \rho(\Gamma) y \rangle_\mathbb{K}.
	\end{equation}
	We then have that
	\begin{equation}
		\langle \rho(e_i)x,\rho(e_i)y \rangle = \langle x,y \rangle
	\end{equation}
	and
	\begin{equation}
	\setlength\arraycolsep{2pt}
	\begin{array}{lcl}
		\langle \rho(e_i)x,\rho(e_j)y \rangle 
			= \langle \rho(e_i)\rho(e_i)x,\rho(e_i)\rho(e_j)y \rangle
			= \pm \langle x,\rho(e_i)\rho(e_j)y \rangle \\[3pt]
			\quad = \mp \langle x,\rho(e_j)\rho(e_i)y \rangle
			= \mp \langle \rho(e_j)x,\rho(e_j)\rho(e_j)\rho(e_i)y \rangle \\[3pt]
			\quad = - \langle \rho(e_j)x,\rho(e_i)y \rangle
	\end{array}
	\end{equation}
	for $e_i \neq e_j$ in $E$. Thus, if $e = \sum_i a_i e_i$ and $\sum_i a_i^2 = 1$, we obtain
	\begin{equation}
		\langle \rho(e)x,\rho(e)y \rangle = \sum_{i,j} a_i a_j \langle \rho(e_i)x,\rho(e_i)y \rangle = \langle x,y \rangle.
	\end{equation}
	Hence, this inner product has the desired property.
	Also note that, for $v \in \mathcal{V} = \mathbb{R}^{n,0}$, we have
	\begin{equation} \label{rho_acts_symmetric}
		\langle \rho(v)x,y \rangle = \langle x,\rho(v)y \rangle,
	\end{equation}
	while for $\mathcal{V} = \mathbb{R}^{0,n}$,
	\begin{equation} \label{rho_acts_antisymmetric}
		\langle \rho(v)x,y \rangle = - \langle x,\rho(v)y \rangle,
	\end{equation}
	i.e. $\rho(v)$ is symmetric for euclidean spaces and antisymmetric for
	anti-euclidean spaces.
	
	We are now ready for some examples which illustrate how representations
	of geometric algebras can appear in various contexts and how their
	representation theory can be used to prove important theorems.
	
	
\subsection{Example I: Normed division algebras}
	
	Our first example concerns the possible dimensions of normed division algebras. 
	A \emph{normed division algebra} is an algebra $\mathcal{A}$ over $\mathbb{R}$ 
	(not necessarily associative) with a norm $|\cdot|$ such that 
	\begin{equation} \label{norm_div_algebra}
		|xy|=|x||y|
	\end{equation}
	for all $x,y \in \mathcal{A}$
	and such that every nonzero element is invertible. We will prove the following
	
	\begin{thm}[Hurwitz' Theorem] \label{thm_hurwitz}
		If $\mathcal{A}$ is a finite-dimensional normed division algebra 
		over $\mathbb{R}$, then its dimension is either 1, 2, 4 or 8.
	\end{thm}
	
	\begin{rem}
		This corresponds uniquely to $\mathbb{R}$, $\mathbb{C}$, $\mathbb{H}$, and the octonions $\mathbb{O}$, respectively.
		The proof of unicity requires some additional steps, see e.g. \cite{baez}.
	\end{rem}
	
	Let us first consider the restrictions that the requirement \eqref{norm_div_algebra} puts on the norm.
	Assume that $\mathcal{A}$ has dimension $n$.
	For every $a \in \mathcal{A}$ we have a linear transformation
	\begin{displaymath}
	\setlength\arraycolsep{2pt}
	\begin{array}{lccl}
		L_a\!:	& \mathcal{A}	& \to 		& \mathcal{A}, \\
				& x 			& \mapsto 	& ax
	\end{array}
	\end{displaymath}
	given by left multiplication by $a$. When $|a|=1$ we then have
	\begin{equation}
		|L_a x|=|ax|=|a||x|=|x|,
	\end{equation}
	i.e. $L_a$ preserves the norm. Hence, it maps the unit sphere $S:=\{x \in \mathcal{A}: |x|=1\}$
	in $\mathcal{A}$ into itself. Furthermore, since every element in $\mathcal{A}$ is invertible,
	we can for each pair $x,y \in S$ find an $a \in S$ such that $L_a x = ax = y$.
	Now, these facts imply a large amount of symmetry of $S$. In fact, we have the following
	
	\begin{lem} \label{lem_operator_symmetry}
		Assume that $V$ is a finite-dimensional normed vector space.
		Let $S_V$ denote the unit sphere in $V$. If, for every $x,y \in S_V$,
		there exists an operator $L \in \textrm{\emph{End}}(V)$ such that
		$L(S_V) \subseteq S_V$ and $L(x)=y$, then $V$ must be an inner product space.
	\end{lem}
	\begin{proof}
		We will need the following fact:
		Every compact subgroup $G$ of $\textrm{GL}(n)$ preserves some inner product on $\mathbb{R}^n$.
		This can be proven by picking a Haar-measure $\mu$ on $G$ and averaging any
		inner product $\langle \cdot,\cdot \rangle$ on $\mathbb{R}^n$ over $G$ using this measure,
		\begin{equation}
			\langle x,y \rangle_G := \int_G \langle gx,gy \rangle\ d\mu(g).
		\end{equation}
		
		Now, let $G$ be the group of linear transformations on $V \cong \mathbb{R}^n$ 
		which preserve its norm $|\cdot|$.
		$G$ is compact in the finite-dimensional operator norm topology, since $G = \bigcap_{x \in V} \{L \in \textrm{End}(\mathbb{R}^n) : |Lx|=|x|\}$
		is closed and bounded by 1. Furthermore, $L \in G$ is injective and therefore an isomorphism.
		The group structure is obvious. Hence, $G$ is a compact subgroup of $\textrm{GL}(n)$.
		
		From the above we know that there exists an inner product $\langle \cdot,\cdot \rangle$ on $\mathbb{R}^n$
		which is preserved by $G$. Let $|\cdot|_\circ$ denote the norm associated to this inner
		product, i.e. $|x|_\circ^2 = \langle x,x \rangle$.
		Take a point $x \in \mathbb{R}^n$ with $|x|=1$ and rescale the inner 
		product so that also $|x|_\circ = 1$.
		Let $S$ and $S_\circ$ denote the unit spheres associated to $|\cdot|$ and $|\cdot|_\circ$,
		respectively. By the conditions in the lemma, there is for every $y \in S$
		an $L \in G$ such that $L(x)=y$.
		But $G$ also preserves the norm $|\cdot|_\circ$, so $y$ must also lie in $S_\circ$.
		Hence, $S$ is a subset of $S_\circ$. However, being unit spheres associated to norms,
		$S$ and $S_\circ$ are both homeomorphic
		to the standard sphere $S^{n-1}$, so we must have that they are equal.
		Therefore, the norms must be equal.
	\end{proof}
	
	We now know that our normed division algebra $\mathcal{A}$ has some
	inner product $\langle \cdot,\cdot \rangle$ such that $\langle x,x \rangle = |x|^2$.
	We call an element $a \in \mathcal{A}$ \emph{imaginary} if $a$ is orthogonal to the
	unit element, i.e. if $\langle a,1_\mathcal{A} \rangle = 0$.
	Let $\textrm{Im}\ \mathcal{A}$ denote the $(n-1)$-dimensional subspace of imaginary elements.
	We will observe that $\textrm{Im}\ \mathcal{A}$ acts on $\mathcal{A}$ in a special way.
	
	Take a curve $\gamma\!: (-\epsilon,\epsilon) \to S$ on the unit sphere such that
	$\gamma(0)=1_\mathcal{A}$ and $\gamma'(0)=a \in \textrm{Im}\ \mathcal{A}$.
	(Note that $\textrm{Im}\ \mathcal{A}$ is the tangent space to $S$ at the unit element.)
	Then, because the product in $\mathcal{A}$ is continuous,
	\begin{equation}
	\setlength\arraycolsep{2pt}
	\begin{array}{rcl}
		\frac{d}{dt}\big|_{t=0} L_{\gamma(t)} x 
			&=& {\displaystyle \lim_{h \to 0}}\ \frac{1}{h}\big(L_{\gamma(h)} x - L_{\gamma(0)} x\big) \\[8pt]
			&=& {\displaystyle \lim_{h \to 0}}\ \frac{1}{h}\big(\gamma(h) - \gamma(0)\big)x
			= \gamma'(0)x = ax = L_a x
	\end{array}
	\end{equation}
	and
	\begin{equation}
	\setlength\arraycolsep{2pt}
	\begin{array}{rcl}
		0 &=& \frac{d}{dt}\big|_{t=0} \langle x,y \rangle 
			= \frac{d}{dt}\big|_{t=0} \langle L_{\gamma(t)}x,L_{\gamma(t)}y \rangle \\[8pt]
			&=& \langle \frac{d}{dt}\big|_{t=0} L_{\gamma(t)}x,L_{\gamma(0)}y \rangle + \langle L_{\gamma(0)}x,\frac{d}{dt}\big|_{t=0} L_{\gamma(t)}y \rangle \\[8pt]
			&=& \langle L_a x,y \rangle + \langle x,L_a y \rangle.
	\end{array}
	\end{equation}
	Hence, $L_a^* = -L_a$ for $a \in \textrm{Im}\ \mathcal{A}$. 
	If, in addition, $|a|=1$ we have that $L_a \in \textrm{O}(\mathcal{A},|\cdot|^2)$, so
	$L_a^2 = - L_a L_a^* = - id$. For an arbitrary imaginary element $a$ we obtain by rescaling
	\begin{equation} \label{normed_div_square}
		L_a^2 = -|a|^2.
	\end{equation}
	This motivates us to consider the geometric algebra $\mathcal{G}(\textrm{Im}\ \mathcal{A},q)$
	with quadratic form $q(a) := -|a|^2$.
	By \eqref{normed_div_square} and the universal property of geometric algebras (Proposition \ref{prop_universality}) 
	we find that $L$ extends to a representation of $\mathcal{G}(\textrm{Im}\ \mathcal{A},q)$ on $\mathcal{A}$,
	\begin{equation} \label{normed_div_rep}
		\hat{L}\!: \mathcal{G}(\textrm{Im}\ \mathcal{A},q) \to \textrm{End}(\mathcal{A}),
	\end{equation}
	i.e. a representation of $\mathcal{G}(\mathbb{R}^{0,n-1})$ on $\mathbb{R}^n$.
	The representation theory now demands that $n$ is 
	a multiple of $d_{0,n-1}$. By studying Table \ref{table_representations}
	and taking periodicity \eqref{rep_periodicity} into account
	we find that this is only possible for $n=1,2,4,8$.

\subsection{Example II: Vector fields on spheres}
	
	In our next example we consider the $N$-dimensional unit spheres $S^N$
	and use representations of geometric algebras to construct vector fields on them.
	The number of such vector fields that can be found gives us information
	about the topological features of these spheres.
	
	\begin{thm}[Radon-Hurwitz] \label{thm_vec_field}
		On $S^N$ there exist $n_N$ pointwise linearly independent vector fields, where, if we write $N$ uniquely as
		\begin{equation}
			N+1 = (2t+1)2^{4a+b}, \quad t,a \in \mathbb{N},\ b \in \{0,1,2,3\},
		\end{equation}
		then
		\begin{equation}
			n_N = 8a + 2^b - 1.
		\end{equation}
		For example,
		\begin{displaymath}
		\begin{array}{c|ccccccccccccccccc}
			N	& 0 & 1 & 2 & 3 & 4 & 5 & 6 & 7 & 8 & 9 &10 &11 &12 &13 &14 &15 &16 \\
			\hline
			n_N	& 0 & 1 & 0 & 3 & 0 & 1 & 0 & 7 & 0 & 1 & 0 & 3 & 0 & 1 & 0 & 8 & 0
		\end{array}
		\end{displaymath}
	\end{thm}

	\begin{cor}
		$S^1$, $S^3$ and $S^7$ are parallelizable.
	\end{cor}

	\begin{rem}
		The number of vector fields constructed in this way 
		is actually the maximum number of possible such fields on $S^N$.
		This is a much deeper result proven by Adams \cite{adams} using algebraic topology.
	\end{rem}
	
	Our main observation is that if $\mathbb{R}^{N+1}$ is a $\mathcal{G}(\mathbb{R}^{0,n})$-module
	then we can construct $n$ pointwise linearly independent vector fields on
	$S^N = \{x \in \mathbb{R}^{N+1} : \langle x,x \rangle=1\}$. 
	Namely, suppose we have a representation $\rho$ of $\mathcal{G}(\mathbb{R}^{0,n})$ on $\mathbb{R}^{N+1}$.
	Take an inner product $\langle\cdot,\cdot\rangle$ on $\mathbb{R}^{N+1}$ such that
	the action of $\rho$ is orthogonal and pick any basis $\{e_1,\ldots,e_n\}$ of $\mathbb{R}^{0,n}$.
	We can now define a collection of smooth vector fields $\{V_1,\ldots,V_n\}$ on $\mathbb{R}^{N+1}$ by
	\begin{equation}
		V_i(x) := \rho(e_i)x, \quad i=1,\ldots,n.
	\end{equation}
	According to the observation \eqref{rho_acts_antisymmetric} this action is
	antisymmetric, so that
	\begin{equation}
		\langle V_i(x),x \rangle = \langle \rho(e_i)x,x \rangle = 0.
	\end{equation}
	Hence, $V_i(x) \in T_x S^N$ for $x \in S^N$. By restricting to $S^N$ we therefore
	have $n$ tangent vector fields. It remains to show that these are pointwise
	linearly independent.
	Take $x \in S^N$ and consider the linear map
	\begin{equation}
	\setlength\arraycolsep{2pt}
	\begin{array}{rccl}
		i_x\!: 	& \mathbb{R}^{0,n} &\to& T_x S^N \\[3pt]
				& v &\mapsto& i_x(v) := \rho(v)x
	\end{array}
	\end{equation}
	Since the image of $i_x$ is $\textrm{Span}_\mathbb{R} \{V_i(x)\}$ it is sufficient
	to prove that $i_x$ is injective. But if $i_x(v)=\rho(v)x=0$ then also
	$v^2 x = \rho(v)^2 x = 0$, so we must have $v=0$.
	
	Now, for a fixed $N$ we want to find as many vector fields as possible, so we seek the
	highest $n$ such that $\mathbb{R}^{N+1}$ is a $\mathcal{G}(\mathbb{R}^{0,n})$-module.
	From the representation theory we know that this requires that $N+1$ is a multiple
	of $d_{0,n}$. Furthermore, since $d_{0,n}$ is a power of 2 we obtain the maximal
	such $n$ when $N+1=p 2^m$, where $p$ is odd and $d_{0,n} = 2^m$.
	Using Table \ref{table_representations} and the periodicity \eqref{rep_periodicity}
	we find that if we write $N+1 = p2^{4a+b}$, with $0 \leq b \leq 3$, then $n = 8a + 2^b - 1$.
	This proves the theorem.

	\newpage

\section{Spinors in physics} \label{sec_spinors}

	
	In this final section we discuss how the view of spinors as even
	multivectors can be used to reformulate physical theories in a way
	which clearly expresses the geometry of these theories, and
	therefore leads to conceptual simplifications.
	
	In the geometric picture provided by geometric algebra we consider
	spinor fields as (smooth) maps from the space or spacetime $\mathcal{V}$
	into the spinor space of $\mathcal{G}(\mathcal{V})$,
	\begin{equation}
		\Psi \!: \mathcal{V} \to \mathcal{G}^+(\mathcal{V}).
	\end{equation}
	The field could also take values in the spinor space of a subalgebra of $\mathcal{G}$.
	For example, a relativistic complex scalar field living
	on Minkowski spacetime could be considered as a spinor field taking values
	in a plane subalgebra $\mathcal{G}^+(\mathbb{R}^{0,2}) \subset \mathcal{G}(\mathbb{R}^{1,3})$.
	
	In the following we will use the summation convention that matching
	upper and lower Greek indices implies summation over 0,1,2,3.
	We will not write out physical constants such as $c,e,\hbar$.
	
\subsection{Pauli spinors}
	
	Pauli spinors describe the spin state of a non-relativistic fermionic
	particle such as the non-relativistic electron. 
	Since this is the non-relativistic limit of the Dirac theory discussed
	below, we will here just state the corresponding
	representation of Pauli spinors as even multivectors of the space algebra.
	We saw that such an element can be written as $\Psi = \rho^{1/2} e^{\varphi \hat{n}I/2}$, 
	i.e. a scaled rotor.
	For this spinor field, the physical state is expressed by the observable
	vector (field)
	\begin{equation} \label{nonrel_observable}
		s := \Psi e_3 \Psi^\dagger = \rho e^{\varphi \hat{n}I/2} e_3 e^{-\varphi \hat{n}I/2},
	\end{equation}
	which is interpreted as the expectation value of the particle's spin,
	scaled by the spatial probability amplitude $\rho$.
	The vector $e_3$ acts as a reference axis for the spin.
	The up and down spin basis states in the ordinary complex representation
	correspond to the rotors which leave $e_3$ invariant, respectively
	the rotors which rotate $e_3$ into $-e_3$. Observe the invariance of $s$
	under right-multiplication of $\Psi$ by $e^{\varphi e_3I}$.
	This corresponds to the complex phase invariance in the
	conventional formulation.
	
\subsection{Dirac-Hestenes spinors}
	
	Dirac spinors describe the state of a relativistic Dirac particle,
	such as an electron or neutrino.
	Conventionally, Dirac spinors are represented by four-component complex
	column vectors, $\psi = [\psi_1,\psi_2,\psi_3,\psi_4]^T \in \mathbb{C}^4$.
	For a spinor \emph{field} the components will be complex-valued functions on spacetime.
	Acting on these spinors are the complex $4 \times 4$-matrices 
	$\{\gamma_0,\gamma_1,\gamma_2,\gamma_3\}$ given in \eqref{gamma_matrices},
	which generate a matrix representation of the Dirac algebra.
	The \emph{Dirac adjoint} of a column spinor is a row matrix 
	\begin{equation}
		\psi^\dagger\gamma_0 = [\psi_1^*,\psi_2^*,-\psi_3^*,-\psi_4^*],
	\end{equation}
	where, in this context, complex conjugation is denoted by $^*$
	and hermitian conjugation by $^\dagger$.
	The physical state of a Dirac particle is determined by the following 16
	so called \emph{bilinear covariants}:
	\begin{equation} \label{bincov_dirac}
	\setlength\arraycolsep{2pt}
	\begin{array}{lcc}
		\alpha &:=& \psi^\dagger \gamma_0 \psi \\[3pt]
		J^\mu  &:=& \psi^\dagger \gamma_0 \gamma^\mu \psi \\[3pt]
		S^{\mu\nu} &:=& \psi^\dagger \gamma_0 i\gamma^\mu \gamma^\nu \psi \\[3pt]
		K^\mu  &:=& \psi^\dagger \gamma_0 iI^{-1} \gamma_\mu \psi \\[3pt]
		\beta  &:=& \psi^\dagger \gamma_0 I^{-1} \psi
	\end{array}
	\end{equation}
	Their integrals over space give expectation values of the physical observables.
	For example, $J^0$ integrated over a spacelike domain gives the
	probability\footnote{Or rather the probability multiplied with the charge of the particle.
	For a large number of particles this can be interpreted as a charge density.}
	of finding the particle in that domain, and
	$J^k$, $k=1,2,3$, give the current of probability. These are components
	of a spacetime current vector $J$.
	The quantities $S^{\mu\nu}$ describe the probability density of the
	particle's electromagnetic moment, while $K^\mu$ gives the direction
	of the particle's spin\footnote{In the formulation below, we obtain the
	relative space spin vector as $\boldsymbol{K} = \frac{\hbar}{2} K \wedge \gamma_0 / |K \wedge \gamma_0|$.}.
	
	In Hestenes' reformulation of the Dirac theory, we represent spinors by
	even multivectors $\Psi \in \mathcal{G}^+$ in the real spacetime algebra $\mathcal{G}(\mathbb{R}^{1,3})$.
	Note that both $\Psi$ and $\psi$ have eight real components, so this is no limitation.
	In this representation, the gamma matrices are considered as orthonormal basis 
	vectors of the Minkowski spacetime and the bilinear covariants are given by
	\begin{equation} \label{bincov_hestenes}
	\setlength\arraycolsep{2pt}
	\begin{array}{rcc}
		\alpha + \beta I &=& \Psi \Psi^\dagger \\[3pt]
		J &=& \Psi \gamma_0 \Psi^\dagger \\[3pt]
		S &=& \Psi \gamma_1\gamma_2 \Psi^\dagger \\[3pt]
		K &=& \Psi \gamma_3 \Psi^\dagger
	\end{array}
	\end{equation}
	where $J=J^\mu \gamma_\mu$, $K=K^\mu \gamma_\mu$ are spacetime
	vectors and $S = \frac{1}{2} S^{\mu\nu} \gamma_\mu \wedge \gamma_\nu$
	a bivector.
	This reformulation allows for a nice geometric interpretation of the Dirac theory.
	Here, spinors are objects that transform the reference basis $\{\gamma_\mu\}$
	into the observable quantities.
	
	Since a spinor only has eight components, the bilinear covariants cannot
	be independent. From \eqref{bincov_hestenes} we easily find a number of relations
	called the \emph{Fierz indentities},
	\begin{equation} \label{fierz_identities}
	\setlength\arraycolsep{2pt}
	\begin{array}{rcc}
		J^2 = - K^2 = \alpha^2 + \beta^2, \quad
		J * K = 0, \quad
		J \wedge K = -(\alpha I + \beta)S.
	\end{array}
	\end{equation}
	The Fierz identities also include a bunch of relations which in the case $\alpha^2+\beta^2 \neq 0$ can be
	derived directly from these three. In total, there are seven degrees of freedom,
	given for example by the spacetime current $J$, the relative space direction of the spin vector $K$
	(two angles) and the so called Yvon-Takabayasi angle $\chi := \arctan (\beta/\alpha)$.
	The eighth degree of freedom present in a spinor is the phase-invariance,
	which in the original Dirac theory corresponds to the overall complex phase of $\psi$,
	while in the Dirac-Hestenes picture corresponds to a rotational freedom in the $\gamma_1\gamma_2$-plane,
	or equivalently around the spin axis in relative space. This is explained by the invariance
	of \eqref{bincov_hestenes} under a transformation $\Psi \mapsto \Psi e^{\varphi \boldsymbol{e}_3 I}$.
	
	In the null case, i.e. when $\alpha = \beta = 0$, we have the additional identities
	\begin{equation} \label{fierz_identities_null}
		S^2 = 0, \quad JS = SJ = 0, \quad KS = SK = 0.
	\end{equation}
	The geometric interpretation is that $J \propto K$ are both null vectors
	and $S$ is a null bivector blade with $J$ and $K$ in its null subspace.
	Hence, the remaining five degrees of freedom are given by the direction
	of the plane represented by S, which must be tangent to the light-cone (two angles),
	plus the magnitudes of $S$, $J$ and $K$.

	The equation which describes the evolution of a Dirac spinor in
	spacetime is the Dirac equation, which in this representation is given by
	the \emph{Dirac-Hestenes equation},
	\begin{equation} \label{dirac_eq}
		\nabla\Psi\gamma_1\gamma_2 - A\Psi = m\Psi\gamma_0,
	\end{equation}
	where $\nabla := \gamma^\mu \frac{\partial}{\partial x^\mu}$
	(we use the spacetime coordinate expansion $x = x^\mu \gamma_\mu$)
	and $A$ is the electromagnetic potential vector field.
	Here, $\gamma_1\gamma_2$ again plays the role of the complex imaginary unit $i$.
	
	Another interesting property of the STA is that the electomagnetic
	field is most naturally represented as a bivector field in $\mathcal{G}(\mathbb{R}^{1,3})$.
	We can write any such bivector field as $F = \boldsymbol{E} + I\boldsymbol{B}$,
	where $\boldsymbol{E}$ and $\boldsymbol{B}$ are relative space vector fields.
	In the context of relativistic electrodynamics,
	these are naturally interpreted as the electric and magnetic fields, respectively.
	Maxwell's equations are compactly written as
	\begin{equation} \label{maxwells_eq}
		\nabla F = J
	\end{equation}
	in this formalism, where $J$ is the source current.
	The physical quantity describing the energy and momentum present
	in an electromagnetic field is the \emph{Maxwell stress-energy tensor}
	which in the STA formulation can be interpreted as a map 
	$T\!: \mathbb{R}^{1,3} \to \mathbb{R}^{1,3}$, given by
	\begin{equation} \label{stress_energy}
		\textstyle
		T(x) := -\frac{1}{2}FxF = \frac{1}{2}FxF^\dagger.
	\end{equation}
	For example, the energy of the field $F$ relative to the $\gamma_0$-direction
	is $\gamma_0 * T(\gamma_0) = \frac{1}{2}(\boldsymbol{E}^2 + \boldsymbol{B}^2$).
	
	Rodrigues and Vaz \cite{vaz_rodrigues_93}, \cite{vaz_rodrigues_95} have studied 
	an interesting correspondence between
	the Dirac and Maxwell equations. With the help of the following theorem,
	they have proved that the electromagnetic field can be obtained from a
	spinor field satisfying an equation similar to the Dirac equation.
	This theorem also serves to illustrate how efficiently computations can be
	performed in the STA framework.
	
	\begin{thm} \label{thm_bivector_spinor}
		Any bivector $F \in \mathcal{G}^{2}(\mathbb{R}^{1,3})$ can be written as
		\begin{displaymath}
			F = \Psi \gamma_0 \gamma_1 \Psi^\dagger,
		\end{displaymath}
		for some (nonzero) spinor $\Psi \in \mathcal{G}^{+}(\mathbb{R}^{1,3}).$
	\end{thm}
	\begin{proof}
		Take any bivector $F = \boldsymbol{E} + I\boldsymbol{B} \in \mathcal{G}^2$.
		Note that 
		\begin{equation}
			F^2 = (\boldsymbol{E}^2 - \boldsymbol{B}^2) + 2 \boldsymbol{E}*\boldsymbol{B} I = \rho e^{\phi I}
		\end{equation}
		for some $\rho \geq 0$ and $0 \leq \phi < 2\pi$.
		We consider the cases $F^2 \neq 0$ and $F^2 = 0$ separately.
		
		If $F^2 \neq 0$ then $\boldsymbol{E}^2 - \boldsymbol{B}^2$ and $\boldsymbol{E}*\boldsymbol{B}$ are not both zero
		and we can apply a \emph{duality rotation} of $F$ into
		\begin{equation}
			F' = \boldsymbol{E}' + I\boldsymbol{B}' := e^{-\phi I/4} F e^{-\phi I/4} \quad \Rightarrow \quad F'^2 = \rho, 
		\end{equation}
		i.e. such that $\boldsymbol{E}'^2 - \boldsymbol{B}'^2 > 0$ and $\boldsymbol{E}'*\boldsymbol{B}' = 0$.
		Hence, we can select an orthonormal basis $\{\boldsymbol{e}_i\}$ of the relative space,
		aligned so that $\boldsymbol{E}' = E' \boldsymbol{e}_1$ and $\boldsymbol{B}' = B' \boldsymbol{e}_2$,
		where $E' = |\boldsymbol{E}'|$ etc.
		Consider now a boost $\boldsymbol{a} = \alpha \boldsymbol{e}_3$ of 
		angle $\alpha$ in the direction orthogonal to both $\boldsymbol{E}'$ and $\boldsymbol{B}'$.
		Using that
		\begin{equation}
			e^{-\frac{\alpha}{2} \boldsymbol{e}_3} \boldsymbol{e}_1 e^{\frac{\alpha}{2} \boldsymbol{e}_3} 
				= \boldsymbol{e}_1 e^{\alpha \boldsymbol{e}_3} 
				= \boldsymbol{e}_1 (\cosh \alpha + \sinh \alpha\ \boldsymbol{e}_3)
		\end{equation}
		and likewise for $\boldsymbol{e}_2$, we obtain
		\begin{equation}
		\setlength\arraycolsep{2pt}
		\begin{array}{rrl}
			F'' &:=& e^{-\boldsymbol{a}/2} F' e^{\boldsymbol{a}/2} 
				= E' e^{-\frac{\alpha}{2} \boldsymbol{e}_3} \boldsymbol{e}_1 e^{\frac{\alpha}{2} \boldsymbol{e}_3} + IB' e^{-\frac{\alpha}{2} \boldsymbol{e}_3} \boldsymbol{e}_2 e^{\frac{\alpha}{2} \boldsymbol{e}_3} \\[5pt]
				&=& E' \boldsymbol{e}_1 (\cosh \alpha + \sinh \alpha\ \boldsymbol{e}_3) + IB' \boldsymbol{e}_2 (\cosh \alpha + \sinh \alpha\ \boldsymbol{e}_3) \\[5pt]
				&=& (E' \cosh \alpha - B' \sinh \alpha) \boldsymbol{e}_1 + I (B' \cosh \alpha - E' \sinh \alpha) \boldsymbol{e}_2 \\[5pt]
				&=& \cosh \alpha \big( (E' - B' \tanh \alpha) \boldsymbol{e}_1 + I (B' - E' \tanh \alpha) \boldsymbol{e}_2 \big),
		\end{array}
		\end{equation}
		where we also noted that $\boldsymbol{e}_1 \boldsymbol{e}_3 = -I \boldsymbol{e}_2$
		and $I \boldsymbol{e}_2 \boldsymbol{e}_3 = -\boldsymbol{e}_1$.
		Since $E'^2 - B'^2 > 0$ we can choose $\alpha := \textrm{arctanh}(\frac{B'}{E'})$ 
		and obtain $F'' = \sqrt{1 - (\frac{B'}{E'})^2} \boldsymbol{E}' = E'' \boldsymbol{e}_1$,
		where $E'' > 0$.
		Finally, some relative space rotor $e^{I\boldsymbol{b}/2}$ takes $\boldsymbol{e}_1$ to 
		our timelike target blade (relative space vector) $\gamma_0 \gamma_1$, i.e.
		\begin{equation}
			F'' = E'' e^{I\boldsymbol{b}/2} \gamma_0 \gamma_1 e^{-I\boldsymbol{b}/2}.
		\end{equation}
		Summing up, we have that $F = \Psi \gamma_0 \gamma_1 \Psi^\dagger$, where
		\begin{equation}
			\Psi = \sqrt{E''} e^{\phi I/4} e^{\boldsymbol{a}/2} e^{I\boldsymbol{b}/2} \in \mathcal{G}^+.
		\end{equation}
		
		When $F^2 = 0$ we have that both $\boldsymbol{E}^2 = \boldsymbol{B}^2$ 
		and $\boldsymbol{E}*\boldsymbol{B}=0$. Again, we select an orthonormal basis $\{\boldsymbol{e}_i\}$
		of the relative space so that $\boldsymbol{E} = E \boldsymbol{e}_1$ 
		and $\boldsymbol{B} = B \boldsymbol{e}_2 = E \boldsymbol{e}_2$.
		Note that
		\begin{equation}
		\setlength\arraycolsep{2pt}
		\begin{array}{l}
			(1 - I \boldsymbol{e}_1 \boldsymbol{e}_2) \boldsymbol{e}_1 (1 + I \boldsymbol{e}_1 \boldsymbol{e}_2)
				= \boldsymbol{e}_1 - I \boldsymbol{e}_1 \boldsymbol{e}_2 \boldsymbol{e}_1 + I \boldsymbol{e}_2 - I \boldsymbol{e}_1 \boldsymbol{e}_2 I \boldsymbol{e}_2 \\[5pt]
				\qquad = 2(\boldsymbol{e}_1 + I \boldsymbol{e}_2).
		\end{array}
		\end{equation}
		Thus, $\frac{1}{\sqrt{2}}(1 - I \boldsymbol{e}_1 \boldsymbol{e}_2) \boldsymbol{E} \frac{1}{\sqrt{2}}(1 + I \boldsymbol{e}_1 \boldsymbol{e}_2) = \boldsymbol{E} + I \boldsymbol{B}$.
		Using that $\boldsymbol{e}_1$ can be obtained from $\gamma_0 \gamma_1$ 
		with some relative space rotor $e^{I\boldsymbol{b}/2}$, we have that
		$F = \Psi \gamma_0 \gamma_1 \Psi^\dagger$, where
		\begin{equation}
			\Psi = ({\textstyle \frac{E}{2}})^{1/2} (1 - {\textstyle \frac{1}{E^2}}I\boldsymbol{E}\boldsymbol{B}) e^{I\boldsymbol{b}/2} \in \mathcal{G}^+.
		\end{equation}
		The case $F=0$ can be achieved not only using $\Psi=0$, but also with e.g. $\Psi=(1+\gamma_0\gamma_1)$.
	\end{proof}

	\noindent
	Note that we can switch $\gamma_0\gamma_1$ for any other non-null reference blade, e.g. $\gamma_1\gamma_2$.
	
	\begin{rem}
		In the setting of electrodynamics, where $F = \boldsymbol{E} + I\boldsymbol{B}$ is an
		electromagnetic field, we obtain as a consequence of this theorem and proof the following
		result due to Rainich, Misner and Wheeler.
		If we define an \emph{extremal field} as a field for which the magnetic (electric)
		field is zero and the electric (magnetic) field is parallel to one coordinate axis,
		the theorem of Rainich-Misner-Wheeler says that:
		``At any point of Minkowski spacetime any nonnull electromagnetic field can
		be reduced to an extremal field by a Lorentz transformation and a duality rotation."
	\end{rem}

	
	The reformulation of the Pauli and Dirac theory 
	observables above depended on the choice
	of fixed reference bases $\{e_1,e_2,e_3\}$ and $\{\gamma_0,\gamma_1,\gamma_2,\gamma_3\}$.
	When a different basis $\{e_1',e_2',e_3'\}$ is selected, we cannot apply the
	same spinor $\Psi$ in \eqref{nonrel_observable} since this
	would in general yield an $s' = \Psi e_3' \Psi^\dagger \neq s$.
	This is not a problem in flat space since we can set up a globally defined
	field of such reference frames without ambiguity. However, in the
	covariant setting of a curved manifold, i.e. when gravitation is involved,
	we cannot fix a certain field of reference frames, but must allow a variation
	in these and equations which transform covariantly under such variations.
	In fibre bundle theory this corresponds to picking different sections
	of an orthonormal frame bundle.
	We therefore seek a formulation of spinor that takes care of this
	required covariance. 
	Rodrigues, de Souza, Vaz and Lounesto \cite{rodrigues_et_al} have
	considered the following definition.
	
	\begin{defn} \label{def_DHS}
		A \emph{Dirac-Hestenes spinor} (DHS) is an equivalence class of triplets $(\Sigma,\psi,\Psi)$,
		where $\Sigma$ is an oriented orthonormal basis of $\mathbb{R}^{1,3}$, 
		$\psi$ is an element in $\textrm{Spin}^+(1,3)$, and $\Psi \in \mathcal{G}^+(\mathbb{R}^{1,3})$ is
		the representative of the spinor in the basis $\Sigma$.
		We define the equivalence relation by $(\Sigma,\psi,\Psi) \sim (\Sigma_0,\psi_0,\Psi_0)$
		if and only if $\Sigma = \widetilde{\textrm{Ad}}_{\psi \psi_0^{-1}} \Sigma_0$ 
		and $\Psi = \Psi_0 \psi_0 \psi^{-1}$.
		The basis $\Sigma_0$ should be thought of as a fixed \emph{reference basis} and 
		the choice of $\psi_0$ is arbitrary but fixed for this basis.
		We suppress this choice and write just $\Psi_\Sigma$ for the spinor $(\Sigma,\psi,\Psi)$.
	\end{defn}
	
	\noindent
	Note that when for example $J = \Psi_E^{\phantom{\dagger}} e_0 \Psi_E^\dagger$ for some basis $E = \{e_i\}$
	we now have the desired invariance property $J = \Psi_{E'}^{\phantom{\dagger}} e_0' \Psi_{E'}^\dagger$
	for some other basis $E' = \{e_i'\}$. Hence, $J$ is now a completely basis independent
	object which in the Dirac theory represents the physical and observable local current 
	produced by a Dirac particle.
	
	The definition above allows for the construction of a covariant Dirac-Hestenes
	spinor field. The possibility of defining such a field on a certain
	manifold depends on the existence of a so called spin structure on it.
	Geroch \cite{geroch} has shown that in the spacetime case, i.e. when the tangent space is $\mathbb{R}^{1,3}$,
	this is equivalent to the existence of a globally defined field of
	time-oriented orthonormal reference frames. 
	In other words, the principal $\textrm{SO}^+$-bundle
	of the manifold must be trivial.
	We direct the reader to \cite{rodrigues_et_al} for a continued discussion.
	
	We end by mentioning that other types of spinors can be represented in the STA as well.
	See e.g. \cite{francis_kosowsky} for a discussion on Lorentz,
	Majorana and Weyl spinors.
	
	


	\newpage

\section{Summary and discussion}

	We have seen that a vector space endowed with a quadratic form
	naturally embeds in an associated geometric algebra. This
	algebra depends on the signature and dimension of the underlying
	vector space, and expresses the geometry of the space through the 
	properties of its multivectors.
	By introducing a set of products, involutions and other operations,
	we got access to the rich structure of this algebra and could
	identify certain significant types of multivectors, such as blades, rotors, and spinors.
	Blades were found to represent subspaces of the underlying vector
	space and gave a geometric interpretation to multivectors and the
	various algebraic operations, while rotors connected the groups
	of structure-respecting transformations to corresponding
	groups embedded in the algebra. This enabled a powerful
	encoding of rotations using geometric products and 
	allowed us to identify candidates for spinors in arbitrary dimensions.
	
	The introduced concepts were put to practice when we worked
	out a number of lower-dimensional examples. These had obvious
	applications in mathematics and physics.
	Norm functions were found to act as determinants on the respective 
	algebras and helped us find the corresponding groups of invertible elements.
	We noted that the properties of such norm functions also lead to totally non-geometric
	applications in number theory.
	
	We also studied the relation between geometric algebras and
	matrix algebras, and used the well-known representation theory
	of such algebras to work out the corresponding representations
	of geometric algebras. The dimensional restrictions of such
	representations led to 
	proofs of classic theorems regarding normed division algebras
	and vector fields on spheres.

	Throughout our examples, we saw that complex structures appear naturally within real geometric
	algebras and that many formulations in physics which involve complex
	numbers can be identified as structures within real geometric
	algebras. Such identifications also resulted in various geometric interpretations
	of complex numbers.
	This suggests that, whenever complex numbers appear in an otherwise
	real or geometric context, one should ask oneself if not
	a real geometric interpretation can be given to them.
	
	
	The combinatorial construction of Clifford algebra which we introduced mainly
	served as a tool for understanding the structure of geometric or Clifford algebras
	and the behaviour of, and relations between, the different products.
	This construction also expresses the generality of Clifford algebras in that they can
	be defined and find applications in general algebraic contexts. 
	Furthermore, it gives new suggestions for how to proceed with the
	infinite-dimensional case.
	Combinatorial Clifford algebra has previously been applied to simplify proofs
	in graph theory \cite{svensson_naeve}.

	Finally, we considered examples in physics and in particular
	relativistic quantum mechanics, where the representation of spinors
	as even multivectors in the geometric algebra of spacetime led
	to conceptual simplifications.
	The resulting picture is a rather classical
	one, with particles as fields of operations which rotate and
	scale elements of a reference basis into the observable expectation
	values.
	Although this is a geometric and conceptually powerful view,
	it is unfortunately not that enlightening with respect to the quantum mechanical
	aspects of states and measurement. This requires an operator-eigenvalue formalism 
	which of course can be formulated in geometric algebra, but
	sort of breaks the geometric picture.
	The geometric view of spinors \emph{does} fit in the context of
	quantum field theory, however, since spinor fields
	there already assume a classical character. 
	It is not clear what conceptual simplifications that geometric
	algebras can bring to other quantum mechanical theories than
	the Pauli and Dirac ones, since
	most realistic particle theories are formulated in infinite-dimensional 
	spaces. 
	Doran and Lasenby \cite{doran_lasenby} have presented suggestions for a 
	multi-particle formulation in geometric algebra,
	but it still involves a fixed and finite number of particles.
	
	Motivated by the conceptual simplifications of Dirac theory brought
	by the spacetime algebra, one can argue about the geometric significance of all particles.
	The traditional classification of particles in terms of spin
	quantum numbers relies on the complex representation theory of the
	(inhomogeneous) Lorentz group. 
	There are complex (or rather complexified) representations of the
	STA-em\-bedded scalar, spinor (through the Dirac algebra), and
	vector fields, corresponding to spin 0, $\frac{1}{2}$, and 1, respectively.
	Coincidentally, the fundamental particles that have been experimentally verified
	all have spin quantum numbers $\frac{1}{2}$ or 1,
	corresponding to spinor fields and vector fields. Furthermore,
	the proposed Higgs particle is a scalar field with spin 0.
	Since all these types of fields are naturally represented within the STA
	it then seems natural to me that exactly these spins have turned up.
	

\addcontentsline{toc}{section}{Acknowledgements}
\section*{Acknowledgements}

	I would like to thank my supervisor Lars Svensson for many long and 
	interesting discussions, and for giving me a lot of freedom to 
	follow my own interests in the subject of geometric algebra.
	I am also grateful for his decision to give introductory lectures
	to first year students about geometric algebra and other
	mathematical topics which for some reason are considered
	controversial. This is what spawned my interest in mathematics
	in general and geometric algebra in particular.

	\newpage

\addcontentsline{toc}{section}{Appendix: Matrix theorems}
\section*{Appendix: Matrix theorems}
	\renewcommand{\thesection}{A}
	\renewcommand{\theequation}{\arabic{equation}}

	In order to avoid long digressions in the text, we have placed
	proofs to some, perhaps not so familiar, theorems in this appendix.
	
	In the following theorem we assume that $R$ is an arbitrary commutative ring and
	\begin{displaymath}
		A = \left[
		\begin{array}{ccc}
			a_{11} & \cdots & a_{1m} \\
			\vdots &		& \vdots \\
			a_{n1} & \cdots & a_{nm} \\
		\end{array}
		\right] \in R^{n \times m}, \qquad
		a_j = \left[
		\begin{array}{c}
			a_{1j} \\
			\vdots \\
			a_{nj} \\
		\end{array}
		\right],
	\end{displaymath}
	i.e. $a_j$ denotes the $j$:th column in $A$.
	If $I \subseteq \{1,\ldots,n\}$ and $J \subseteq \{1,\ldots,m\}$ we
	let $A_{I,J}$ denote the $|I| \times |J|$-matrix minor obtained from $A$
	by deleting the rows and columns not in $I$ and $J$.
	Further, let $k$ denote the rank of $A$, i.e. the highest integer $k$
	such that there exists $I,J$ as above with $|I|=|J|=k$ and $\det A_{I,J} \neq 0$.
	By renumbering the $a_{ij}$:s we can without loss of generality assume 
	that $I=J=\{1,2,\ldots,k\}$.
	
	\begin{thm}[Basis minor]
		If the rank of $A$ is $k$, and
		\begin{displaymath}
			d := \det \left[
			\begin{array}{ccc}
				a_{11} & \cdots & a_{1k} \\
				\vdots &		& \vdots \\
				a_{k1} & \cdots & a_{kk} \\
			\end{array}
			\right] \neq 0,
		\end{displaymath}
		then every $d \cdot a_j$ is a linear combination of $a_1,\ldots,a_k$.
	\end{thm}
	
	\begin{proof}
		Pick $i \in \{1,\ldots,n\}$ and $j \in \{1,\ldots,m\}$ and consider
		the $(k+1)\times(k+1)$-matrix
		\begin{displaymath}
			B_{i,j} := \left[
			\begin{array}{cccc}
				a_{11} & \cdots & a_{1k} & a_{1j} \\
				\vdots &		& \vdots & \vdots \\
				a_{k1} & \cdots & a_{kk} & a_{kj} \\
				a_{i1} & \cdots & a_{ik} & a_{ij} \\
			\end{array}
			\right].
		\end{displaymath}
		Then $\det B_{i,j} = 0$. Expanding $\det B_{i,j}$ along the bottom row
		for fixed $i$ we obtain
		\begin{equation}
			a_{i1} C_1 + \ldots + a_{ik} C_k + a_{ij} d = 0,
		\end{equation}
		where the $C_l$ are independent of the choice of $i$ (but of course dependent on $j$).
		Hence,
		\begin{equation}
			C_1 a_1 + \ldots + C_k a_k + d a_j = 0,
		\end{equation}
		and similarly for all $j$.
	\end{proof}
	
	The following shows that the factorization
	$\det(AB) = \det(A) \det(B)$
	is a unique property of the determinant.
	
	\begin{thm}[Uniqueness of determinant]
		Assume that $d\!: \mathbb{R}^{n \times n} \to \mathbb{R}$ is continuous
		and satisfies
		\begin{equation}
			d(AB) = d(A)d(B)
		\end{equation}
		for all $A,B \in \mathbb{R}^{n \times n}$.
		Then $d$ must be either $0$, $1$, $|\det|^\alpha$ or $(\textrm{\emph{sign}} \circ \det) |\det|^\alpha$
		for some $\alpha > 0$.
	\end{thm}
	
	\begin{proof}
		First, we have that
		\begin{equation}
		\setlength\arraycolsep{2pt}
		\begin{array}{rcl}
			d(0) = d(0^2) = d(0)^2, \\[5pt]
			d(I) = d(I^2) = d(I)^2,
		\end{array}
		\end{equation}
		so $d(0)$ and $d(I)$ must be either 0 or 1.
		Furthermore,
		\begin{equation}
		\setlength\arraycolsep{2pt}
		\begin{array}{rcl}
			d(0) = d(0A) = d(0)d(A), \\[5pt]
			d(A) = d(IA) = d(I)d(A),
		\end{array}
		\end{equation}
		for all $A \in \mathbb{R}^{n \times n}$,
		which implies that $d=1$ if $d(0)=1$ and $d=0$ if $d(I)=0$.
		We can therefore assume that $d(0)=0$ and $d(I)=1$.
		
		Now, an arbitrary matrix $A$ can be written as
		\begin{equation}
			A = E_1 E_2 \ldots E_k R,
		\end{equation}
		where $R$ is on reduced row-echelon form (as close to the identity
		matrix as possible) and $E_i$ are elementary row operations of the form
		\begin{equation}
		\setlength\arraycolsep{2pt}
		\begin{array}{lcl}
			R_{ij} &:=& (\textrm{swap rows $i$ and $j$}), \\[5pt]
			E_i(\lambda) &:=& (\textrm{scale row $i$ by $\lambda$}),\ \textrm{or} \\[5pt]
			E_{ij}(c) &:=& (\textrm{add $c$ times row $j$ to row $i$}).
		\end{array}
		\end{equation}
		Because $R_{ij}^2 = I$, we must have $d(R_{ij}) = \pm 1$. This gives, since 
		\begin{equation}
			E_i(\lambda) = R_{1i} E_1(\lambda) R_{1i}, 
		\end{equation}
		that $d\big(E_i(\lambda)\big) = d\big(E_1(\lambda)\big)$ and
		\begin{equation}
			d(\lambda I) = d\big(E_1(\lambda) \ldots E_n(\lambda)\big) = d\big(E_1(\lambda)\big) \ldots d\big(E_n(\lambda)\big) = d\big(E_1(\lambda)\big)^n.
		\end{equation}
		In particular, we have $d\big(E_1(0)\big) = 0$ and of course $d\big(E_1(1)\big) = d(I) = 1$.
		
		If $A$ is invertible, then $R=I$. Otherwise, $R$ must contain a row of zeros
		so that $R = E_i(0)R$ for some $i$. But then $d(R)=0$ and $d(A)=0$.
		When $A$ is invertible we have $I=AA^{-1}$ and $1=d(I)=d(A)d(A^{-1})$,
		i.e. $d(A) \neq 0$ and $d(A^{-1}) = d(A)^{-1}$. Hence,
		\begin{equation}
			A \in \textrm{GL}(n) \quad \Leftrightarrow \quad d(A) \neq 0.
		\end{equation}
		We thus have that $d$ is completely determined by its values on $R_{ij}$,
		$E_1(\lambda)$ and $E_{ij}(c)$. Note that we have not yet used the
		continuity of $d$, but it is time for that now.
		We can split $\mathbb{R}^{n \times n}$ into three connected components,
		namely $\textrm{GL}^-(n)$, $\det^{-1}(0)$ and $\textrm{GL}^+(n)$,
		where the determinant is less than, equal to, and greater than zero, respectively.
		Since $E_1(1), E_{ij}(c) \in \textrm{GL}^+(n)$ and 
		$E_1(-1), R_{ij} \in \textrm{GL}^-(n)$, we have by continuity of $d$ that
		\begin{equation} \label{sign_of_d}
		\setlength\arraycolsep{2pt}
		\begin{array}{lcl}
			d(R_{ij}) =+1 &\quad \Rightarrow \quad& \textrm{$d$ is $>0$, $=0$, resp. $>0$} \\[5pt]
			d(R_{ij}) =-1 &\quad \Rightarrow \quad& \textrm{$d$ is $<0$, $=0$, resp. $>0$} \\[5pt]
		\end{array}
		\end{equation}
		on these parts.
		Using that $d\big(E_1(-1)\big)^2 = d\big(E_1(-1)^2\big) = d(I) = 1$, we have
		$d\big(E_1(-1)\big) = \pm 1$ and $d\big(E_1(-\lambda)\big) = d\big(E_1(-1)\big) d\big(E_1(\lambda)\big) = \pm d\big(E_1(\lambda)\big)$
		where the sign depends on \eqref{sign_of_d}.
		On $\mathbb{R}^{++} := \{ \lambda \in \mathbb{R} : \lambda > 0 \}$ we have a continuous map
		$d \circ E_1\!: \mathbb{R}^{++} \to \mathbb{R}^{++}$ such that
		\begin{equation}
			d \circ E_1(\lambda \mu) = d \circ E_1(\lambda) \cdot d \circ E_1(\mu) \quad \forall\ \lambda,\mu \in \mathbb{R}^{++}.
		\end{equation}
		Forming $f := \ln \circ\ d \circ E_1 \circ \exp$, we then have a continuous map
		$f\!: \mathbb{R} \to \mathbb{R}$ such that
		\begin{equation}
			f(\lambda + \mu) = f(\lambda) + f(\mu).
		\end{equation}
		By extending linearity from $\mathbb{Z}$ to $\mathbb{Q}$ and $\mathbb{R}$ by continuity,
		we must have that $f(\lambda) = \alpha \lambda$ for some $\alpha \in \mathbb{R}$.
		Hence, $d \circ E_1(\lambda) = \lambda^\alpha$. Continuity also demands that $\alpha > 0$.
		
		It only remains to consider $d \circ E_{ij} \!: \mathbb{R} \to \mathbb{R}^{++}$.
		We have $d \circ E_{ij}(0) = d(I) = 1$ and $E_{ij}(c) E_{ij}(\gamma) = E_{ij}(c + \gamma)$,
		i.e.
		\begin{equation} \label{d_add_to_mult}
			d \circ E_{ij}(c + \gamma) = d \circ E_{ij}(c) \cdot d \circ E_{ij}(\gamma) \quad \forall\ c,\gamma \in \mathbb{R}.
		\end{equation}
		Proceeding as above, $g := \ln \circ\ d \circ E_{ij} \!: \mathbb{R} \to \mathbb{R}$
		is linear, so that $g(c) = \alpha_{ij} c$ for some $\alpha_{ij} \in \mathbb{R}$,
		hence $d \circ E_{ij}(c) = e^{\alpha_{ij} c}$.
		One can verify that the following identity holds for all $i,j$:
		\begin{equation}
			E_{ji}(-1) = E_i(-1) R_{ij} E_{ji}(1) E_{ij}(-1).
		\end{equation}
		This gives $d\big(E_{ji}(-1)\big) = (\pm 1) (\pm 1) d\big(E_{ji}(1)\big) d\big(E_{ij}(-1)\big)$
		and, using \eqref{d_add_to_mult},
		\begin{equation}
		\setlength\arraycolsep{2pt}
		\begin{array}{rcl}
			d\big(E_{ij}(1)\big) 
				&=& d\big(E_{ji}(2)\big) = d\big(E_{ji}(1+1)\big) = d\big(E_{ji}(1)\big)d\big(E_{ji}(1)\big) \\[5pt]
				&=& d\big(E_{ij}(2)\big)d\big(E_{ij}(2)\big) = d\big(E_{ij}(4)\big),
		\end{array}
		\end{equation}
		which requires $\alpha_{ij}=0$.
		
		We conclude that $d$ is completely determined by $\alpha > 0$,
		where $d \circ E_1(\lambda) = \lambda^\alpha$ and $\lambda \geq 0$,
		plus whether $d$ takes negative values or not. This proves the theorem.
	\end{proof}

\end{document}